\documentclass[a4paper,12pt,leqno]{article}

\usepackage{amsmath}
\usepackage{amscd}
\usepackage[psamsfonts]{amssymb}
\usepackage{graphicx}
\usepackage{lscape}
\usepackage[french]{babel}
\usepackage[utf8]{inputenc}
\usepackage[all]{xy}
\usepackage{color}
\usepackage{footmisc}
\usepackage[toc,page]{appendix}
\usepackage{titlesec}
\usepackage{hyperref}

\newcounter{sect}
\newcounter{subsect}[sect]
\newcommand{\sect}[1]
{\refstepcounter{subsect}\refstepcounter{sect}
\renewcommand{\thesubsect}{\thesect}
{\large{\bf\rule{0ex}{10ex}\thesect}
~~#1\rule[-2ex]{0ex}{2ex}\hfill}}

\newcommand{\subsect}[1]{
\renewcommand{\thesubsect}{\arabic{sect}.\arabic{subsect}}
\refstepcounter{subsect}
{{\bf\rule{0ex}{4ex}\thesubsect}
~~#1\rule[-2ex]{0ex}{2ex}\hfill}}

\newcounter{th}[subsect]
\renewcommand{\theth}{\thesubsect.\arabic{th}}
\newenvironment{theo}[1][]{\refstepcounter{th}\rule{0ex}{4ex}%
\textsc{Th\'{e}or\`{e}me}~\theth#1.~~\em}{\em%
\rule[-1.5ex]{0ex}{1.5ex}}

\newenvironment{pro}[1][]{\refstepcounter{th}\rule{0ex}{4ex}%
\textsc{Proposition}~\theth#1.~~\em}{\em\rule[-1.5ex]{0ex}{1.5ex}}

\newenvironment{lem}[1][]{\refstepcounter{th}\rule{0ex}{4ex}%
\textsc{Lemme}~\theth#1.~~\em}{\em\rule[-1.5ex]{0ex}{1.5ex}}

\newenvironment{cor}[1][]{\refstepcounter{th}\rule{0ex}{4ex}%
\textsc{Corollaire}~\theth#1.~~\em}{\em\rule[-1.5ex]{0ex}{1.5ex}}

\newenvironment{pro-def}[1][]{\refstepcounter{th}\rule{0ex}{4ex}%
\textsc{Proposition-D\'{e}finition}~\theth#1.~~\em}{\em\rule[-1.5ex]{0ex}{1.5ex}}

\newenvironment{cor-def}[1][]{\refstepcounter{th}\rule{0ex}{4ex}%
\textsc{Corollaire-D\'{e}finition}~\theth#1.~~\em}{\em\rule[-1.5ex]{0ex}{1.5ex}}

\newenvironment{lem-def}[1][]{\refstepcounter{th}\rule{0ex}{4ex}%
\textsc{Lemme-D\'{e}finition}~\theth#1.~~\em}{\em\rule[-1.5ex]{0ex}{1.5ex}}

\newenvironment{scho}[1][]{\refstepcounter{th}\rule{0ex}{4ex}%
\textsc{Scholie}~\theth#1.~~\em}{\em\rule[-1.5ex]{0ex}{1.5ex}}

\newenvironment{scho-def}[1][]{\refstepcounter{th}\rule{0ex}{4ex}%
\textsc{Scholie-D\'{e}finition}~\theth#1.~~\em}{\em\rule[-1.5ex]{0ex}{1.5ex}}

\def\numero{$\mathrm{n}^{\text{o}}$}

\begin{document}

\begin{center}

\large
\textbf{Sur l'application de Quillen pour la cohomologie modulo $2$ de certains groupes finis}

\normalsize
\vspace{0.375cm}
\textsc{Jean Lannes}

\end{center}

\sect{Introduction}

\bigskip
\textsc{La théorie de Quillen}

\bigskip
Soit $G$ un groupe fini arbitraire. Le célèbre article de Quillen \cite{Qui}  fournit en particulier une ``approximation'' de $\mathrm{H}^{*}(G;\mathbb{F}_{2})$, la cohomologie modulo $2$ du groupe $G$. On rappelle ci-après la théorie de Quillen dans ce contexte.

\medskip
On considère la catégorie, disons $\mathcal{Q}_{G}$, suivante~:

\smallskip
-- les objets de $\mathcal{Q}_{G}$ sont les $2$-sous-groupes abéliens élémentaires $E\subset G$ ($E\approx(\mathbb{Z}/2)^{d}$ pour un certain $d$) ;

\smallskip
-- les morphismes de $\mathcal{Q}_{G}$ sont les homomorphismes de groupes $f:E\to E'$ induits par une conjugaison dans $G$.

\smallskip
On note respectivement $\mathrm{ob}\hspace{0.5pt}\mathcal{Q}_{G}$ et  $\mathrm{mor}\hspace{0.5pt}\mathcal{Q}_{G}$ l'ensemble des objets et l'ensemble des morphismes de  $\mathcal{Q}_{G}$~; on obervera que ces deux ensembles sont finis.

\medskip
Comme une conjugaison dans $G$ induit l'identité de $\mathrm{H}^{*}(G;\mathbb{F}_{2})$, le produit des homomorphismes de restriction
$$
\mathrm{H}^{*}(G;\mathbb{F}_{2})\to\prod_{E\in\mathrm{ob}\hspace{0.5pt}\mathcal{Q}_{G}}
\mathrm{H}^{*}(E;\mathbb{F}_{2})
$$
induit  un homomorphisme de $\mathbb{F}_{2}$-algèbres (commutatives) $\mathbb{N}$-graduées de $\mathrm{H}^{*}(G;\mathbb{F}_{2})$ sur la sous-algèbre de $\prod_{E\in\mathrm{ob}\hspace{0.5pt}\mathcal{Q}_{G}}\mathrm{H}^{*}(E;\mathbb{F}_{2})$ constituée des  éléments ${(x_{E})}_{E\in\mathrm{ob}\hspace{0.5pt}\mathcal{Q}_{G}}$ vérifiant $x_{E}=f^{*}x_{E'}$ pour tout $f\in\mathrm{mor}\hspace{0.5pt}\mathcal{Q}_{G}$. Nous notons $\mathrm{L}(G)$ cette sous-algèbre et
$$
\mathrm{q}_{G}:\mathrm{H}^{*}(G;\mathbb{F}_{2})\to\mathrm{L}(G)
$$
l'homomorphisme défini ci-dessus~; nous appelons $\mathrm{q}_{G}$  {\em l'application de Quillen}.

\medskip
\begin{rem}\label{limite} Soit $\mathcal{K}^{\mathbb{F}_{2}}$ la catégorie des $\mathbb{F}_{2}$-algèbres $\mathbb{N}$-graduées. Quillen observe que l'application $E\mapsto\mathrm{H}^{*}(E;\mathbb{F}_{2})$ peut être vue comme un foncteur défini sur la catégorie $\mathcal{Q}_{G}^{\mathrm{op}}$, opposée de  $\mathcal{Q}_{G}$, à valeurs dans $\mathcal{K}^{\mathbb{F}_{2}}$ et que l'objet $\mathrm{L}(G)$ de  $\mathcal{K}^{\mathbb{F}_{2}}$ n'est rien d'autre que la limite de ce foncteur~:
$$
\hspace{24pt}
\mathrm{L}(G)
\hspace{4pt}:=\hspace{4pt}
\lim_{\mathcal{Q}_{G}^{\mathrm{op}}}\mathrm{H}^{*}(E;\mathbb{F}_{2})
\hspace{23pt}.
$$
\end{rem}

\bigskip
Le théorème 6.2 de \cite{Qui} (``main theorem'') se spécialise de la façon suivante~:

\begin{theo}\label{Quillen1} Soit $G$ un groupe fini.

\medskip
{\em (a)} Tout élément du noyau de $\mathrm{q}_{G}$ est nilpotent.

\medskip
{\em (b)} Pour tout élément $x$ de $\mathrm{L}(G)$ il existe un entier naturel $r$ tel que $x^{2^{r}}$ est dans l'image de $\mathrm{q}_{G}$.

\end{theo}

\begin{rem}\label{a-bis} Pour le point (a) ci-dessus nous avons respecté la formulation de \cite{Qui}. Comme $\mathrm{H}^{*}((\mathbb{Z}/2)^{d};\mathbb{F}_{2})$ est isomorphe à une algèbre de polynômes $\mathbb{F}_{2}[u_{1},u_{2},\ldots,u_{d}]$ ($u_{1},u_{2},\ldots,u_{d}$ désignant des indéterminées de degré un) la $\mathbb{F}_{2}$-algèbre $\mathrm{L}(G)$ est réduite si bien que tout élément nilpotent de $\mathrm{H}^{*}(G;\mathbb{F}_{2})$ est dans le noyau de $\mathrm{q}_{G}$. On peut donc remplacer le point (a) en question par le suivant~:

\medskip
(a-bis) {\em Le noyau de $\mathrm{q}_{G}$ est le nilradical de $\mathrm{H}^{*}(G;\mathbb{F}_{2})$.}
\end{rem}

\vspace{0.75cm}
\textsc{Entrée en scène de l'algèbre de Steenrod}

\medskip
On observe que la cohomologie modulo $2$ d'un groupe $G$ est la cohomologie modulo $2$ d'un espace topologique à savoir son espace classifiant $\mathrm{B}G$~: $\mathrm{H}^{*}(G;\mathbb{F}_{2})=\mathrm{H}^{*}(\mathrm{B}G;\mathbb{F}_{2})$. Soit $\mathrm{A}$ l'algèbre de Steenrod modulo $2$~; $\mathrm{H}^{*}(G;\mathbb{F}_{2})$ est donc muni d'une action de $\mathrm{A}$ qui en fait un {\em $\mathrm{A}$-module instable}. Rappelons qu'un $\mathrm{A}$-module instable est un $\mathbb{F}_{2}$-espace vectoriel $\mathbb{N}$-gradué $M=\{M^{n}\}_{n\in\mathbb{N}}$ muni d'applications $\mathrm{A}^{i}\otimes M^{n}\to M^{n+i}$, qui en font un $\mathrm{A}$-module ($\mathbb{N}$-gradué), tel que l'on a $\mathrm{Sq}^{i}x=0$ pour $i>\vert x\vert$, $\vert x\vert$ désignant le degré d'un élément (homogène) de $M$ ({\em condition d'instabilité})\footnote{On pourrait supposer que $M$ est $\mathbb{Z}$-gradué~; en effet la condition d'instabilité entraîne que $x=\mathrm{Sq}^{0}x$ est nul pour $\vert x\vert<0$. La terminologie ``$\mathrm{A}$-module instable'' peut paraître étrange~: si l'on convient qu'un $\mathrm{A}$-module stable est un $\mathrm{A}$-module ne vérifiant pas la condition d'instabilité alors ce n'est rien d'autre qu'un $\mathrm{A}$-module $\mathbb{Z}$-gradué. Le mot ``stable'' dans ce contexte fait référence à la théorie de l'homotopie stable (voir par exemple \cite{Ad}).}. La catégorie dont les objets sont les $\mathrm{A}$-modules instables et dont les morphismes sont les applications $\mathrm{A}$-linéaires de degré zéro est notée $\mathcal{U}$. 
\medskip
De plus $\mathrm{H}^{*}(G;\mathbb{F}_{2})=\mathrm{H}^{*}(\mathrm{B}G;\mathbb{F}_{2})$ est une  {\em $\mathrm{A}$-algèbre instable}, c'est-à-dire que l'on a~:

\smallskip
($\mathcal{K}_{1}$)\hspace{8pt}$\mathrm{Sq}^{i}(x\smile y)=\sum_{j+k=i}\mathrm{Sq}^{j}x\smile\mathrm{Sq}^{k}y$ pour tous $x$ et $y$ dans $\mathrm{H}^{*}(G;\mathbb{F}_{2})$~;

\smallskip
($\mathcal{K}_{2}$)\hspace{8pt}$\mathrm{Sq}^{\vert x\vert}x=x\smile x$ pour tout $x$ dans $\mathrm{H}^{*}(G;\mathbb{F}_{2})$.

\smallskip
La catégorie des $\mathrm{A}$-algèbres instables est notée $\mathcal{K}$~; par définition elle est munie de deux foncteurs ``oubli''~:
$$
\hspace{24pt}
\begin{CD}
\mathcal{K}^{\mathbb{F}_{2}}@<<<\mathcal{K}@>>>\mathcal{U}
\end{CD}
\hspace{23pt}.
$$

\smallskip
La condition ($\mathcal{K}_{2}$) ci-dessus conduit aux définitions ci-dessous.

\medskip
Soit $M$ un $\mathrm{A}$-module instable~; on note $\mathrm{Sq}_{0}$ l'application, de $M$ dans $M$, $x\mapsto\mathrm{Sq}^{\vert x\vert}x$ ($\mathrm{Sq}_{0}$ n'est pas de degré zéro mais multiplie le degré par $2$). On dit qu'un élément de $M$ est {\em nilpotent} s'il est annulé par une itérée de $\mathrm{Sq}_{0}$. On dit que $M$ est {\em nilpotent} si tous ses éléments sont nilpotents. En général, soit $\mathrm{Nil}(M)$ le sous-ensemle ($\mathbb{N}$-gradué) de $M$ constitué des éléments nilpotents~; $\mathrm{Nil}(M)$ est stable par addition, la relation $\mathrm{Sq}_{0}\hspace{1pt}\mathrm{Sq}^{i}=\mathrm{Sq}^{2i}\hspace{1pt}\mathrm{Sq}_{0}$ implique qu'il est aussi stable sous l'action des opérations de Steenrod~: c'est le plus grand sous-$\mathrm{A}$-module instable  nilpotent  de $M$ qu'il est raisonnable d'appeler le {\em nilradical} de~$M$.

\bigskip
L'application $\mathrm{q}_{G}:\mathrm{H}^{*}(G;\mathbb{F}_{2})\to\mathrm{L}(G)$ est un homomorphisme de $\mathrm{A}$-algèbres instables. Le théorème \ref{Quillen1} fait intervenir l'homomorphisme sous-jacent de $\mathbb{F}_{2}$-algèbres $\mathbb{N}$-graduées~; ce qui précède montre que l'on peut le reformuler en termes de l'homomorphisme  sous-jacent de $\mathrm{A}$-modules instables~:

\begin{theo}\label{Quillen2} Soit $G$ un groupe fini~; le noyau et le conoyau de l'homomorphisme de $\mathrm{A}$-modules instables sous-jacent à $\mathrm{q}_{G}$ sont  nilpotents. 

\end{theo}

Sans surprise, on dit qu'un $\mathrm{A}$-module instable $M$ est {\em réduit} si $0$ est son seul élément nilpotent.

\medskip
La catégorie $\mathcal{U}$ est une catégorie abélienne avec assez d'injectifs (et de projectifs). On constate que  $M$ est réduit si et seulement si $\mathrm{Hom}_{\mathcal{U}}(N,M)=\mathrm{Ext}_{\mathcal{U}}^{0}(N,M)$ est nul pour tout $N$ nilpotent (pour s'en convaincre prendre $N=\mathrm{Nil}(M)$). On dit que $M$ est {\em $\mathcal{N}\hspace{-1.5pt}il$-fermé} si l'on a en outre $\mathrm{Ext}_{\mathcal{U}}^{1}(N,M)=0$ pour tout $N$ nilpotent. Un peu d'algèbre homologique dans la catégorie $\mathcal{U}$, voir \ref{nilloc-1} et \ref{dempropiso}, fournit l'énoncé suivant~:

\begin{pro}\label{iso} Soit $G$ un groupe fini~; les deux conditions suivantes sont équivalentes~:
\begin{itemize}
\item[(i)] L'application de Quillen $\mathrm{q}_{G}$  est un isomorphisme.
\item[(ii)] Le $\mathrm{A}$-module instable $\mathrm{H}^{*}(G;\mathbb{F}_{2})$ est $\mathcal{N}\hspace{-1.5pt}il$-fermé.
\end{itemize}
\end{pro}

Soit $S$ un sous-groupe d'indice impair de $G$~; le fait que l'homomorphisme de transfert $\mathrm{tr}:\mathrm{H}^{*}(S;\mathbb{F}_{2})\to\mathrm{H}^{*}(G;\mathbb{F}_{2})$ commute aux opérations de Steenrod entraîne que $\mathrm{H}^{*}(G;\mathbb{F}_{2})$ est (canoniquement)  facteur direct, comme $\mathrm{A}$-module instable, de $\mathrm{H}^{*}(S;\mathbb{F}_{2})$. La définition-même de la notion de $\mathrm{A}$-module instable $\mathcal{N}\hspace{-1.5pt}il$-fermé en termes des foncteurs $\mathrm{Ext}_{\mathcal{U}}^{k}(N,-)$, $k=0,1$ et $N$ nilpotent, implique donc~:

\begin{pro}\label{Sylow-1} Soient $G$ un groupe fini et $S\subset G$ un sous-groupe d'indice  impair. Si le $\mathrm{A}$-module instable $\mathrm{H}^{*}(S;\mathbb{F}_{2})$ est $\mathcal{N}\hspace{-1.5pt}il$-fermé alors il en est de même pour le $\mathrm{A}$-module instable $\mathrm{H}^{*}(G;\mathbb{F}_{2})$.
\end{pro}

\bigskip
Les propositions \ref{iso} et \ref{Sylow-1} conduisent aux corollaires suivants~:

\begin{cor}\label{Sylow-2} Soient $G$ un groupe fini et $S\subset G$ un sous-groupe d'indice  impair. Si le $\mathrm{A}$-module instable $\mathrm{H}^{*}(S;\mathbb{F}_{2})$ est $\mathcal{N}\hspace{-1.5pt}il$-fermé alors $\mathrm{q}_{G}$ est un isomorphisme .
\end{cor}

\bigskip
On en vient maintenant aux résultats de notre article.

\bigskip
\textsc{Trois familles de groupes finis $G$ telles que la cohomologie modulo $2$ d'un $2$-Sylow de $G$ est $\mathcal{N}\hspace{-1.5pt}il$-fermée}

\medskip
On montre dans \cite{GLZ} que si la cohomologie modulo $2$ d'un groupe fini $G$ est $\mathcal{N}\hspace{-1.5pt}il$-fermée (comme $\mathrm{A}$-module instable) alors il en est de même pour celle du produit en couronne\footnote{La notation $G\wr\mathfrak{S}_{2}$ semble plus fréquente.} $\mathfrak{S}_{2}\wr G$. On en déduit  que la cohomologie modulo $2$ d'un $2$-Sylow du groupe symétrique $\mathfrak{S}_{n}$ est $\mathcal{N}\hspace{-1.5pt}il$-fermée si l'entier $n$ est une puissance de~$2$~; le cas général en résulte en observant que le produit tensoriel de deux $\mathrm{A}$-modules instables $\mathcal{N}\hspace{-1.5pt}il$-fermés est encore $\mathcal{N}\hspace{-1.5pt}il$-fermé. Le corollaire \ref{Sylow-2} dit alors que l'application de Quillen $\mathrm{q}_{\mathfrak{S}_{n}}$ est un isomorphisme.

\medskip
Par les mêmes méthodes (mais avec un peu plus d'efforts) nous montrons dans le présent  article que la cohomologie modulo $2$ d'un $2$-Sylow du groupe alterné~$\mathfrak{A}_{n}$ est $\mathcal{N}\hspace{-1.5pt}il$-fermée.

\smallskip
Notre stratégie est de remplacer les $\mathrm{A}$-modules instables par les $\mathrm{P}$-$\mathrm{A}$-modules instables, $\mathrm{P}$ désignant la cohomologie modulo $2$ du groupe $\mathbb{Z}/2$. Rappelons la définiton de ces objets. Un $\mathrm{P}$-$\mathrm{A}$-module instable est un $\mathrm{A}$-module instable $M$ muni d'une application $\mathrm{A}$-linéaire $\mathrm{P}\otimes M\to M$ qui fait de $M$ un $\mathrm{P}$-module. Exemple~:  l'homomorphisme de $\mathrm{A}$-algèbres instables  $\mathrm{P}\to\mathrm{H}^{*}(\mathfrak{S}_{n};\mathbb{F}_{2})$,  induit par la signature, fait du $\mathrm{A}$-module instable $\mathrm{H}^{*}(\mathfrak{S}_{n};\mathbb{F}_{2})$ un $\mathrm{P}$-$\mathrm{A}$-module instable.

\smallskip
Comme dans le cas des groupes symétriques nous commençons par traiter le cas  de $\mathfrak{A}_{n}$ avec $n$ une puissance de $2$. Ce cas particulier est au coeur de notre travail~; le passage au cas général se fait par une méthode similaire à celle employée pour les groupes symétriques.

\smallskip
Comme précédemment le corollaire \ref{Sylow-2} dit que l'application de Quillen $\mathrm{q}_{\mathfrak{A}_{n}}$ est un isomorphisme. Ici un commentaire s'impose : Chad Giusti et Dev Sinha ont récemment déterminé $\mathrm{H}^{*}(\mathfrak{A}_{n};\mathbb{F}_{2})$ \cite{GS} et  montré en particulier que  $\mathrm{H}^{*}(\mathfrak{A}_{n};\mathbb{F}_{2})$ est réduit ce qui entraîne que $\mathrm{q}_{\mathfrak{A}_{n}}$ est un monomorphisme. Nous montrons que  $\mathrm{q}_{\mathfrak{A}_{n}}$ est un isomorphisme mais nous ne déterminons pas $\mathrm{L}(\mathfrak{A}_{n})$~!

\medskip
Enfin nous observons que les $2$-Sylow des groupes de Coxeter finis irréducti\-bles (le cas des groupes diédraux mis à part)  sont produits de $2$-Sylow de groupes symétriques ou alternés. La cohomologie modulo $2$ d'un $2$-Sylow d'un groupe de Coxeter fini irréductible est donc encore $\mathcal{N}\hspace{-1.5pt}il$-fermée (le cas des groupes diédraux est traité séparément) si bien qu'il en est de même pour tout  groupe de Coxeter fini $W$. Il en résulte à nouveau que l'application de Quillen $\mathrm{q}_{W}$ est un isomorphisme.

\footnotesize
\begin{center}
\textit{Crédits et remerciements}
\end{center}

\vspace{-0.2cm}
Les recherches exposées dans cet article ont commencé au début des années 2000 à l'initiative de Mark Feshbach (qui connaissait \cite{GLZ} et avait observé que les 2-Sylow des groupes de Coxeter finis irréductibles, les groupes diédraux mis à part, étaient reliés à ceux des groupes symétriques et alternés). Nous avions traité le cas du groupe alterné $\mathfrak{A}_{8}$ en utilisant la structure particulière de ``son'' 2-Sylow.  La mort prématurée de Mark en 2010 n'a malheureusement pas permis à notre collaboration d'aller jusqu'à son terme~; je dédie cet article à sa mémoire.

\smallskip
Je remercie Jean-Pierre  Serre pour les questions qu'il m'a posées sur la cohomologie modulo~$2$ des groupes alternés qui ont renouvelé mon intérêt pour ce sujet.

\smallskip
Je remercie le VIASM qui m'a invité à donner un cours sur mon travail à l'automne 2023. Je remercie en particulier Nguyen The Cuong et Pham Van Tuan qui ont été des cobayes-auditeurs efficaces.

\smallskip
Je remercie l'IRL FVMA CNRS, l'IMJ-PRG et le VIASM pour leur soutien matériel à ma visite à Hanoï.

\smallskip
Je remercie enfin le CMLS pour son hospitalité durant la rédaction de cet article.
\normalsize

\sect{Sur les $\mathrm{A}$-modules instables réduits et $\mathcal{N}\hspace{-1.5pt}il$-fermés}\label{nildefs}

\medskip
Pour le confort du lecteur nous commençons par regrouper quelques-unes des définitions apparues dans l'introduction.

\medskip
\begin{defi}\label{defnil} Soit $M$ un $\mathrm{A}$-module instable.

\smallskip
(a) On note $\mathrm{Sq}_{0}:M\to M$ l'application linéaire $x\mapsto\mathrm{Sq}^{\vert x\vert}x$, $\vert x\vert$ désignant le degré de $x$ (on observera que l'on a $\vert\mathrm{Sq}_{0}x\vert=2\vert x\vert$).

\smallskip
(b) On dit qu'un élément $x$ de $M$ est {\em nilpotent} s'il est annulé par une itérée de l'application $\mathrm{Sq}_{0}$. On dit que le $\mathrm{A}$-module instable $M$ est {\em nilpotent} si tous ses éléments sont nilpotents.

\end{defi}

\medskip
\begin{pro-def}\label{defnilrad} {\em Soit $M$ un $\mathrm{A}$-module instable.

\smallskip
On appelle {\em nilradical} de $M$ le sous-espace vectoriel $\mathbb{N}$-gradué constitué des éléments nilpotents~; on le note $\mathrm{Nil}(M)$.} Le nilradical de $M$ est stable sous l'action des opérations de Steenrod~; c'est le plus grand sous-$\mathrm{A}$-module instable nilpotent de $M$.
\end{pro-def}

\begin{pro-def}\label{defred} {\em Soit $M$ un $\mathrm{A}$-module instable.}

\smallskip
Les trois conditions suivantes sont équivalentes~:
\begin{itemize}
\item[(i)] l'application linéaire $\mathrm{Sq}_{0}:M\to M$ est injective~;
\item[(ii)] le nilradical $\mathrm{Nil}(M)$ est nul~;
\item[(iii)] on a $\mathrm{Hom}_{\mathcal{U}}(N,M)=0$ pour tout $\mathrm{A}$-module instable nilpotent $N$.
\end{itemize}

\medskip
{\em Si ces conditions sont vérifiées on dit que $M$ est {\em réduit}}.
\end{pro-def}

\medskip
\begin{defi}\label{defnf} On dit qu'un  $\mathrm{A}$-module instable $M$ est {\em $\mathcal{N}\hspace{-1.5pt}il$-fermé} si $M$ est réduit et si l'on a en outre $\mathrm{Ext}_{\mathcal{U}}^{1}(N,M)=0$  pour tout $\mathrm{A}$-module instable nilpotent $N$~; en d'autres termes si l'on a $\mathrm{Ext}_{\mathcal{U}}^{k}(N,M)=0$ pour $k=0,1$ et tout $\mathrm{A}$-module instable nilpotent $N$.
\end{defi}

\bigskip
Les trois propositions ci-après  s'obtiennent en considérant la longue suite exacte des $\mathrm{Ext}^{*}_{\mathcal{U}}(N,-)$ avec $N$ nilpotent, associée à une suite exacte courte.

\medskip
\begin{pro}\label{critnilf-0}
Soit $0\to M'\to M\to M''\to 0$ une suite exacte dans la catégorie $\mathcal{U}$ avec $M$ $\mathcal{N}\hspace{-1.5pt}il$-fermé~; les deux propriétés suivantes sont équivalentes~:
\begin{itemize}
\item[(i)] $M''$ est réduit~;
\item[(ii)] $M'$ est $\mathcal{N}\hspace{-1.5pt}il$-fermé.
\end{itemize}
\end{pro}

\begin{scho}\label{critnilf-1}
Soit $0\to M'\to M\to M''$ une suite exacte dans la catégorie $\mathcal{U}$~; si $M$ est $\mathcal{N}\hspace{-1.5pt}il$-fermé et $M''$ réduit alors $M'$~est $\mathcal{N}\hspace{-1.5pt}il$-fermé.
\end{scho}

\medskip
\textit{Démonstration.} Soit $M'''$ l'image de l'homomorphisme $M\to M''$. Comme $M''$ est réduit $M'''$ l'est aussi, si bien que l'on dispose d'une suite exacte courte $0\to M'\to M\to M'''\to 0$ avec $M$ $\mathcal{N}\hspace{-1.5pt}il$-fermé et $M'''$ réduit. On peut donc invoquer l'implication $(i)\Rightarrow(ii)$ de \ref{critnilf-0}.
\hfill$\square$

\medskip
\begin{scho}\label{invnf} Soit $M$ un $\mathrm{A}$-module instable muni d'une action (à gauche) d'un groupe $G$ commutant à l'action de $\mathrm{A}$, en clair muni d'une famille d'automorphismes $(\mathrm{a}_{g}:M\circlearrowleft)_{g\in G}$ de $M$ vérifiant  $\mathrm{a}_{gh}=\mathrm{a}_{g}\circ\mathrm{a}_{h}$ pour tous $g$ et $h$ dans $G$.

\smallskip
Si $M$ est  $\mathcal{N}\hspace{-1,5pt}il$-fermé alors le sous-$\mathrm{A}$-module instable $M^{G}$ constitués des éléments de $M$ invariants par $G$ est aussi $\mathcal{N}\hspace{-1,5pt}il$-fermé.
\end{scho}

\medskip
\textit{Démonstration.} Par définition même on a une suite exacte de $\mathrm{A}$-modules instables
$$
\hspace{24pt}
\begin{CD}
0@>>>M^{G}@>>>M@>f>>\prod_{g\in G}M_{g}
\end{CD}
\hspace{23pt},
$$
$M_{g}$ désignant une copie de $M$ et $f$ le produit des homomorphismes $\mathrm{a}_{g}-\mathrm{id}$~; on achève à l'aide de \ref{critnilf-1} en observant qu'un produit (arbitraire) de $\mathrm{A}$-modules instables $\mathcal{N}\hspace{-1,5pt}il$-fermés est $\mathcal{N}\hspace{-1,5pt}il$-fermé et \textit{a fortiori} réduit.
\hfill$\square$

\medskip
\begin{pro}\label{critnilf-2}
Soit $0\to M'\to M\to M''\to 0$ une suite exacte dans la catégorie $\mathcal{U}$ avec $M''$ $\mathcal{N}\hspace{-1.5pt}il$-fermé~; les deux propriétés suivantes sont équivalentes~:
\begin{itemize}
\item[(i)] $M'$ est $\mathcal{N}\hspace{-1.5pt}il$-fermé~;
\item[(ii)] $M$ est $\mathcal{N}\hspace{-1.5pt}il$-fermé.
\end{itemize}
\end{pro}

\begin{pro}\label{nilloc-1}
Soit $f:M\to L$ un $\mathcal{U}$-morphisme. Si $L$ est $\mathcal{N}\hspace{-1.5pt}il$-fermé et si $\ker f$ et $\mathop{\mathrm{coker}} f$ sont nilpotents alors les deux conditions sont équivalentes~:
\begin{itemize}
\item[(i)] $f$ est un isomorphisme~;
\item[(ii)] $M$ est $\mathcal{N}\hspace{-1.5pt}il$-fermé.
\end{itemize}
\end{pro}

\textit{Démonstration de $(ii)\Rightarrow(i)$.} Si $M$ est réduit il en est de même pour $\ker f$, d'où $\ker f=0$. En considérant la longue suite exacte des $\mathrm{Ext}^{*}_{\mathcal{U}}(N,-)$ avec $N$ nilpotent, associée à la suite exacte courte $0\to M\overset{f}{\to}L\to\mathop{\mathrm{coker}} f\to 0$, on constate que $\mathop{\mathrm{coker}} f$ est réduit, d'où $\mathop{\mathrm{coker}} f=0$.
\hfill$\square$

\bigskip
\textsc{Notation.} Soit $G$ un groupe fini, nous abrègerons  dans la suite de cet article la notation $\mathrm{H}^{*}(G;\mathbb{F}_{2})$ en $\mathrm{H}^{*}G$.

\bigskip
\begin{exple}\label{dempropiso} La proposition \ref{nilloc-1} fournit une  démonstration de l'énoncé \ref{iso}. En effet on a,  par définition de $\mathrm{L}(G)$, une suite exacte dans la catégorie $\mathcal{U}$
$$
0\to\mathrm{L}(G)\to\mathrm{L}^{0}(G)\to\mathrm{L}^{1}(G)
$$
avec
$$
\hspace{24pt}
\mathrm{L}^{0}(G):=\prod_{E\in\mathrm{ob}\hspace{0.5pt}\mathcal{Q}_{G}}
\mathrm{H}^{*}E
\hspace{12pt}\text{et}\hspace{12pt}
\mathrm{L}^{1}(G):=\prod_{f\in\mathrm{mor}\hspace{0.5pt}\mathcal{Q}_{G}}
\mathrm{H}^{*}(\mathrm{source\hspace{3pt}de}\hspace{3pt}f)
\hspace{23pt}.
$$
Or, à l'occasion des recherches sur la conjecture de Sullivan, il a été découvert \cite{Ca}\cite{Mi}\cite{LZens}, que $\mathrm{H}^{*}V$ est un $\mathrm{A}$-module instable injectif (réduit puisque l'on~a $\mathrm{H}^{*}V\cong\mathbb{F}_{2}[u_{1},u_{2},\ldots,u_{d}]$, $\{u_{1},u_{2},\ldots,u_{d}\}$ désignant une base de $\mathrm{H}^{1}V$), pour tout $2$-groupe abélien élémentaire $V$~; il en résulte que $\mathrm{L}(G)$ est $\mathcal{N}\hspace{-1.5pt}il$-fermé~(par exemple d'après \ref{critnilf-1}) si bien que l'on peut appliquer~\ref{nilloc-1}. En fait on verra (condition ($iii)$ de \ref{critnilf-3}) que tout $\mathrm{A}$-module instable $\mathcal{N}\hspace{-1.5pt}il$-fermé prend place dans une suite exacte analogue à celle qui définit $\mathrm{L}(G)$.
\end{exple}

\bigskip
\textit{Origine de la terminologie ``$\mathcal{N}\hspace{-1.5pt}il$-fermé''}

\footnotesize
\medskip
On reprend les hypothèses de \ref{nilloc-1}. Soit $L'$ un $\mathrm{A}$-module instable $\mathcal{N}\hspace{-1.5pt}il$-fermé, alors des arguments du même type que ceux utilisés dans la démonstration de \ref{nilloc-1} montrent que $f$ induit un isomorphisme $\mathrm{Hom}_{\mathcal{U}}(L,L')\cong\mathrm{Hom}_{\mathcal{U}}(M,L')$.  En d'autres termes, étant donné un $\mathcal{U}$-morphisme $f':M\to L'$, il existe un unique $\mathcal{U}$-morphisme $\phi: L\to L'$ tel que l'on a $f'=\phi\circ f$. Si l'on suppose en outre que $f'$ satisfait les hypothèses de la proposition \ref{nilloc-1} alors on constate que $\phi$ est un isomorphisme (échanger les rôles de $f$ et $f'$)~: le $\mathcal{U}$-morphisme $f:M\to L$ est ``unique à isomorphisme canonique près''.

\medskip
Ce qui précède est relié à la théorie de la {\em localisation dans les catégories abéliennes} développée dans \cite[Chap. III]{Ga} où la terminologie (-)-fermé est introduite. Soit $\mathcal{N}\hspace{-1.5pt}il$ la sous-catégorie pleine de $\mathcal{U}$ dont les objets sont les $\mathrm{A}$-modules instables nilpotents. Cette sous-catégorie est {\em épaisse}~: étant donnée une $\mathcal{U}$-suite exacte $0\to M'\to M\to M''\to 0$, $M$ est nilpotent si et seulement si $M'$ et $M''$ le sont. Cette propriété permet de définir la {\em catégorie quotient}  $\mathcal{U}/\mathcal{N}\hspace{-1.5pt}il$ et un foncteur $\mathbf{t}:\mathcal{U}\to\mathcal{U}/\mathcal{N}\hspace{-1.5pt}il$. De plus ce foncteur admet un adjoint à droite $\mathbf{s}:\mathcal{U}/\mathcal{N}\hspace{-1.5pt}il\to\mathcal{U}$ (on dit alors que la sous-catégorie est {\em localisante})~; l'existence de cet adjoint résulte du fait que $\mathcal{U}$ a assez d'injectifs et que tout objet de $\mathcal{U}$ contient un sous-objet maximal qui appartient à $\mathcal{N}\hspace{-1.5pt}il$. Le foncteur $\mathbf{\ell}:=\mathbf{s}\circ\mathbf{t}:\mathcal{U}\circlearrowleft$ est le {\em  foncteur localisation}~; par construction on dispose d'une transformation naturelle $\eta_{M}:M\to\mathbf{\ell}(M)$, à savoir l'unité de l'adjonction, et on constate que $\mathbf{\ell}(M)$ est $\mathcal{N}\hspace{-1.5pt}il$-fermé et que $\ker\eta_{M}$ et $\mathop{\mathrm{coker}}\eta_{M}$ sont nilpotents.
\normalsize

\medskip
La conclusion de la discussion ci-dessus est la suivante : $\mathrm{L}(G)$ est ``la'' localisation ``away from $\mathcal{N}\hspace{-1.5pt}il$'' du $\mathrm{A}$-module instable $\mathrm{H}^{*}G$ et le $\mathcal{U}$-morphisme sous-jacent à $\mathrm{q}_{G}$ s'identifie à $\eta_{\mathrm{H}^{*}G}:\mathrm{H}^{*}G\to\mathrm{L}(G)$. Pour une description ``concrète'' de la catégorie $\mathcal{U}/\mathcal{N}\hspace{-1.5pt}il$ et de l'endofonteur $\mathbf{\ell}:\mathcal{U}\circlearrowleft$, le lecteur pourra consulter \cite{HLS2}.

\bigskip
\textsc{Produit tensoriel de $\mathrm{A}$-modules instables $\mathcal{N}\hspace{-1.5pt}il$-fermés}

\medskip
Soient $M_{1}$ et $M_{2}$ deux $\mathrm{A}$-modules instables.  Le produit tensoriel des deux $\mathbb{F}_{2}$-espaces vectoriels $\mathbb{N}$-gradués sous-jacents est naturellement muni d'une structure de $\mathrm{A}$-module grâce à la structure d'algèbre de Hopf de $\mathrm{A}$\footnote{En clair~: soient $x_{1}$ un élément de $M_{1}$, $x_{2}$ un élément de $M_{2}$ et $i$ un entier naturel, on a la {\em formule de Cartan}~:
$\mathrm{Sq}^{i}(x_{1}\otimes x_{2})\hspace{2pt}=\hspace{2pt}
\sum_{i_{1}+i_{2}=i}
\mathrm{Sq}^{i_{1}}x_{1}\otimes\mathrm{Sq}^{i_{2}} x_{2}$.

\smallskip
\noindent
Profitons de cette note pour faire l'observation suivante~: Soient $K$ une $\mathrm{A}$-algèbre instable et $\varphi:K\otimes K\to K$ l'application donnée par le produit de cette algèbre~; l'axiome ($\mathcal{K}_{1}$) des $\mathrm{A}$-algèbres instables (voir introduction) signifie que $\varphi$ est un $\mathcal{U}$-morphisme \label{axiome-K1}}~; on le note $M_{1}\otimes M_{2}$, c'est encore un $\mathrm{A}$-module instable.

\medskip
On a $\mathrm{Sq}_{0}\hspace{1pt}(x_{1}\otimes x_{2})=\mathrm{Sq}_{0\hspace{1pt}}x_{1}\otimes\mathrm{Sq}_{0}\hspace{1pt}x_{2}$ pour tout $x_{1}$ dans $M_{1}$ et tout $x_{2}$ dans~$M_{2}$. Cette formule implique que si $M_{1}$ et $M_{2}$ sont réduits alors il en est de même pour $M_{1}\otimes M_{2}$. Pareillement~:

\begin{pro}\label{tensnf} Soient $M_{1}$ et $M_{2}$ deux  $\mathrm{A}$-modules instables; si $M_{1}$ et $M_{2}$ sont  $\mathcal{N}\hspace{-1.5pt}il$-fermés  alors il en est de même pour $M_{1}\otimes M_{2}$.
\end{pro}

\medskip
\textit{Démonstration.} Celle-ci utilise la théorie des $\mathcal{U}$-injectifs \cite{LZens}\cite{LS} ($\mathcal{U}$-injectif est une abréviation pour ``objet injectif de la catégorie $\mathcal{U}$ des $\mathrm{A}$-modules instables'') que nous mettons en {\oe}uvre ci-après

\begin{pro}\label{reduit} Soit $M$ un $\mathrm{A}$-module instable.  Les propriétés suivantes sont équivalentes~:
\begin{itemize}
\item[(i)]$M$ est réduit~;
\item[(ii)] il existe un $\mathcal{U}$-monomorphisme $M\hookrightarrow I$ avec $I$ un $\mathcal{U}$-injectif réduit~;
\item[(iii)] il existe une famille $(E_{\gamma})_{\gamma\in\Gamma}$ de $2$-groupes abéliens élémentaires et un $\mathcal{U}$-monomorphisme $M\hookrightarrow \bigoplus_{\gamma\in\Gamma}\mathrm{H}^{*}E_{\gamma}$.
\end{itemize}
\end{pro}

\textit{Démonstration de $(i)\Rightarrow(ii)$.} Soit $i:M\to I$ une enveloppe injective de $M$. Si $M$ est réduit alors on a $i^{-1}(\mathrm{Nil}(I))=0$ et donc $\mathrm{Nil}(I)=0$ (la définition du nilradical $\mathrm{Nil}(-)$ d'un $\mathrm{A}$-module instable est rappelée en \ref{defnilrad}) si bien que $I$ est réduit.
\hfill$\square$

\medskip
\textit{Démonstration de $(ii)\Rightarrow(iii)$.} On pose $\mathrm{P}:=\mathrm{H}^{*}\mathbb{Z}/2$. On considère l’ensemble des facteurs directs indécomposables de $\mathrm{P}^{\otimes m}$, $m$ parcourant $\mathbb{N}$ (par convention $\mathrm{P}^{\otimes m}=\mathbb{F}_{2}$ pour $m=0$), et on choisit un sous-ensemble $\mathcal{L}$ de cet ensemble tel que chaque classe d’isomorphisme de ces facteurs directs ait dans $\mathcal{L}$ un représentant et un seul ; le théorème \cite[6.2.1]{LS} dit qu'il existe une (unique) famille de cardinaux $(a_{L})_{L\in\mathcal{L}}$ telle que l'on a $I\cong\bigoplus_{L\in\mathcal{L}}L^{\oplus a_{L}}$.
\hfill$\square$

\medskip
\textit{Démonstration de $(iii)\Rightarrow(i)$.} Cette implication est évidente puisque les $\mathrm{H}^{*}E_{\gamma}$ sont réduits comme $\mathrm{A}$-modules instables.
\hfill$\square$

\begin{cor}\label{critnilf-3} Soit $M$ un $\mathrm{A}$-module instable.  Les propriétés suivantes sont équivalentes~:
\begin{itemize}
\item[(i)]$M$ est $\mathcal{N}\hspace{-1.5pt}il$-fermé~;
\item[(ii)] il existe une suite exacte dans $\mathcal{U}$ de la forme $0\to M\to I^{0}\to I^{1}$ avec $I^{0}$ et $I^{1}$ des $\mathcal{U}$-injectifs réduits~;
\item[(iii)] il existe deux familles $(E_{\gamma})_{\gamma\in\Gamma}$ et $(F_{\delta})_{\delta\in\Delta}$ de $2$-groupes abéliens élémentai\-res et une suite exacte dans $\mathcal{U}$ de la forme
$$
\hspace{24pt}
0\to M\to\bigoplus_{\gamma\in\Gamma}\mathrm{H}^{*}E_{\gamma}\to \bigoplus_{\delta\in\Gamma}\mathrm{H}^{*}F_{\delta}
\hspace{23pt}.
$$
\end{itemize}
\end{cor}

\textit{Démonstration de $(i)\Rightarrow(ii)$. }Si $M$ est $\mathcal{N}\hspace{-1.5pt}il$-fermé, il est \textit{a fortiori} réduit. La proposition \ref{reduit} dit qu'il existe un monomorphisme $M\hookrightarrow I^{0}$ avec $I^{0}$ un $\mathcal{U}$-injectif réduit~; soit $Q$ le conoyau de ce monomorphisme, l'implication $(ii)\Rightarrow(i)$ de \ref{critnilf-0} montre que $Q$ est réduit. Il existe donc un monomorphisme $Q\hookrightarrow I^{1}$ avec $I^{1}$ un $\mathcal{U}$-injectif réduit.
\hfill$\square$

\medskip
\textit{Démonstration de $(ii)\Rightarrow(i)$.} Puisque la catégorie $\mathcal{U}$ a assez d'injectifs, la suite exacte  $0\to M\to I^{0}\to I^{1}$ se prolonge en une résolution injective
$$
\hspace{24pt}
0\to M\to I^{0}\to I^{1}\to I^{2}\to I^{3}\to\ldots
\hspace{23pt}.
$$
Soit $N$ un $\mathrm{A}$-module instable nilpotent~; $\mathrm{Hom}_{\mathcal{U}}(N,I^{k})=0$ pour $k=0,1$ implique $\mathrm{Ext}_{\mathcal{U}}^{k}(N,M)=0$ pour $k=0,1$.
\hfill$\square$

\medskip
\textit{Démonstration de $(i)\Rightarrow(iii)$.} Elle est analogue à celle de $(i)\Rightarrow(ii)$, l'implication  $(i)\Rightarrow(ii)$ de la proposition \ref{reduit} étant remplacé par l'implication  $(i)\Rightarrow(iii)$ de cette même proposition.
\hfill$\square$

\medskip
\textit{Démonstration de $(iii)\Rightarrow(i)$.} C'est un cas particulier de $(ii)\Rightarrow(i)$ car la cohomologie modulo $2$ d'un groupe abélien élémentaire est un $\mathcal{U}$-injectif réduit.
\hfill$\square$

\bigskip
\textit{Fin de la démonstration de la proposition \ref{tensnf}}

\medskip
On utilise par exemple l'équivalence $(i)\iff(iii)$ du corollaire \ref{critnilf-3}. Soient $M_{1}$ et $M_{2}$ deux $\mathrm{A}$-modules instables $\mathcal{N}\hspace{-1.5pt}il$-fermés et $0\to M_{i}\to I_{i}^{0}\to I_{i}^{1}$, $i=1,2$, deux suites exactes de $\mathrm{A}$-modules instables avec $I_{i}^{k}$ des sommes directes de cohomologie modulo $2$ de $2$-groupes abéliens élémentaires  pour $i=1,2$ et $k=0,1$. On a alors une suite exacte de $\mathrm{A}$-modules instables
$$
\hspace{24pt}
0\to M_{1}\otimes M_{2}\to I_{1,2}^{0}:= I_{1}^{0}\otimes I_{2}^{0}\to I_{1,2}^{1}:=( I_{1}^{0}\otimes I_{2}^{1})\oplus  (I_{1}^{0}\otimes I_{2}^{0})
\hspace{24pt}
$$
et $I_{1,2}^{0}$ et $I_{1,2}^{1}$ sont aussi des sommes directes de cohomologie modulo $2$ de $2$-groupes abéliens élémentaires.
\hfill$\square$

\bigskip
\textsc{Compléments}

\bigskip
\textit{Suspension d'un $\mathrm{A}$-module instable}

\medskip
Soient $M=(M^{n})_{n\in\mathbb{Z}}$ un $\mathrm{A}$-module $\mathbb{Z}$-gradué et $k$ un élément de $\mathbb{Z}$. On note $\Sigma^{k} M$ le $\mathrm{A}$-module $\mathbb{Z}$-gradué $(M^{n-k})_{n\in\mathbb{Z}}$ et $\Sigma^{k} x$ l'élément de degré $n$ de $\Sigma^{k} M$ correspondant à un élément de degré $n-k$ de $M$~; $\Sigma^{k} M$ s'appelle la {\em $k$-suspension} de~$M$, on abrège $\Sigma^{1}M$ en $\Sigma M$ et {\em $1$-suspension} en {\em suspension}.

\medskip
La proposition suivante est évidente~:

\begin{pro}\label{susp} Soit $M$ un $\mathrm{A}$-module instable.

\medskip
{\em (a)} Si $k$ est un élément de $\mathbb{N}$ alors le $\mathrm{A}$-module $\mathbb{Z}$-gradué $\Sigma^{k}M$ est encore un $\mathrm{A}$-module~instable.

\medskip
{\em (b)} Les trois conditions suivantes sont équivalentes~:
\begin{itemize}
\item[(i)] $M$ est une suspension (en clair $M$ est la suspension d'un un $\mathrm{A}$-module instable)~;
\item[(ii)]  le $\mathrm{A}$-module $\mathbb{Z}$-gradué $\Sigma^{-1}M$ est un $\mathrm{A}$-module instable~;
\item[(iii)] l'application $\mathrm{Sq}_{0}:M\to M$ est triviale.
\end{itemize}

\medskip
{\em (c)} Si $M$ est une suspension alors on a $\mathrm{Hom}_{\mathcal{U}}(M,M')=0$ quand le $\mathrm{A}$-module instable $M'$ est réduit.

\end{pro}

\bigskip
\textit{Foncteur $\mathrm{T}_{V}$ et $\mathrm{A}$-modules instables $\mathcal{N}\hspace{-1.5pt}il$-fermés}

\medskip
Soit $V$ un $2$-groupe abélien élémentaire~; on étudie dans \cite{LaDurham}\cite{LaT} le foncteur $\mathrm{T}_{V}:\mathcal{U}\to\mathcal{U}$ adjoint à droite du foncteur $\mathcal{U}\to\mathcal{U},M\mapsto\mathrm{H}^{*}V\otimes M$. On montre en particulier que $\mathrm{T}_{V}$ est exact (cette propriété est implicite dans \cite{LZens}).

\begin{pro}\label{Tnilf} Soit $M$ un $\mathrm{A}$-modules instable~; si $M$ est $\mathcal{N}\hspace{-1.5pt}il$-fermé alors il en est de même pour $\mathrm{T}_{V}M$.
\end{pro}

\medskip
\textit{Démonstration.} Soit $E$ un $2$-groupe abélien élémentaire~; le $\mathrm{A}$-module instable  $\mathrm{T}_{V}\mathrm{H}^{*}E$ est isomorphe à une somme directe de copies de $\mathrm{H}^{*}E$ indexées par l'ensemble $\mathrm{Hom}(V,E)$ (voir \cite[4.2]{LaDurham}). Comme $\mathrm{T}_{V}$ est exact et commute aux sommes directes, on conclut en invoquant l'équivalence $(i)\iff(iii)$ de~\ref{critnilf-3}.
\hfill$\square$

\begin{pro}\label{Gphi} Soient $G$ un groupe fini, $V$ un $2$-groupe abélien élémentaire et  $\phi:V\to G$ un homomorphisme de groupes~; soit $G_{\phi}\subset G$ le sous-groupe de~$G$ centralisateur de $\phi(V)$. Si le $\mathrm{A}$-module instable $\mathrm{H}^{*}G$ est $\mathcal{N}\hspace{-1.5pt}il$-fermé, alors il en est de même pour  le $\mathrm{A}$-module instable $\mathrm{H}^{*}G_{\phi}$.
\end{pro}

\medskip
\textit{Démonstration.}  On pose $\mathrm{Rep}(V,G):=G\backslash\mathrm{Hom}(V,G)$, $G$ agissant à gauche sur l'ensemble $\mathrm{Hom}(V,G)$ par conjugaison au but. On montre (voir la remarque qui suit \cite[3.4.6]{LaT}) que l'on a un isomorphisme canonique de $\mathrm{A}$-modules instables
$$
\hspace{24pt}
\mathrm{T}_{V}\hspace{1pt}\mathrm{H}^{*}G
\hspace{4pt}\cong\hspace{4pt}
\prod_{\phi}\mathrm{H}^{*}G_{\phi}
\hspace{23pt},
$$
$\phi$ décrivant un système de représentants de $\mathrm{Rep}(V,G)$ dans $\mathrm{Hom}(V,G)$ (ce résultat est une généralisation de  \cite[4.2]{LaDurham}). En particulier $\mathrm{H}^{*}G_{\phi}$ est un facteur direct de $\mathrm{T}_{V}\hspace{1pt}\mathrm{H}^{*}G$ qui est $\mathcal{N}\hspace{-1.5pt}il$-fermé d'après \ref{Tnilf}.
\hfill$\square$

\sect{Groupes symétriques}\label{gs}

On montre dans cette section que la cohomologie modulo  $2$ d'un $2$-sous-groupe de Sylow d'un groupe symétrique est $\mathcal{N}\hspace{-1.5pt}il$-fermée~; ce résultat apparaît déjà dans \cite{GLZ}. Notre rédaction est bien plus détaillée que celle de \cite{GLZ} mais sa trame est assez semblable~; en fait le rôle principal de cette section est de servir de modèle à la stratégie que nous adopterons dans la section \ref{Alt} où nous traiterons des groupes alternés.

\medskip
On note $\mathfrak{S}_{n}$($n\in\mathbb{N}$) le groupe des permutations de l'ensemble $\{1,2,\ldots, n\}$ (vide pour $n=0$).

\medskip
Soit $G$ un groupe fini~; on abrège {\em $2$-sous-groupe de Sylow de $G$} en {\em $2$-Sylow de~$G$} (abréviation fréquente dans la littérature).

\subsect{``Le'' $2$-Sylow de $\mathfrak{S}_{n}$}\label{2-Sylow-1}

\smallskip
Le contenu de cette sous-section est fort classique.

\medskip
On considère tout d'abord le cas où $n$ est une puissance de $2$.

\medskip
On note $\mathrm{S}_{2^{m}}$ ($m\in\mathbb{N}$), le groupe défini par récurrence de la façon suivante~:
$$
\mathrm{S}_{2^{0}}=\mathfrak{S}_{1}
\hspace{12pt}\text{et}\hspace{12pt}
\mathrm{S}_{2^{m}}=\mathfrak{S}_{2}\wr\mathrm{S}_{2^{m-1}}:=\mathfrak{S}_{2}\ltimes(\mathrm{S}_{2^{m-1}}\times\mathrm{S}_{2^{m-1}})
$$
($\mathfrak{S}_{2}$ agissant à gauche sur le groupe $\mathrm{S}_{2^{m-1}}\times\mathfrak{S}_{2^{m-1}}$ par permutation des deux facteurs). On constate que $\mathrm{S}_{2^{m}}$ est un $2$-groupe de cardinal $2^{(2^{m}-1)}$.

\medskip
Soit $m\geq 1$ un entier~; la considération de la bijection
$$
\{1,2\}\times\{1,2,\ldots,2^{m-1}\}\to\{1,2,\ldots,2^{m}\}
\hspace{12pt},\hspace{12pt}
(i,k)\mapsto k+(i-1)\hspace{2pt}2^{m-1}
$$
conduit à la définition d'un monomorphisme $\mathfrak{S}_{2}\wr\mathfrak{S}_{2^{m-1}}\hookrightarrow\mathfrak{S}_{2^{m}}$ et donc, par récurrence sur $m$, à l'identification de $\mathrm{S}_{2^{m}}$ avec un sous-groupe de $\mathfrak{S}_{2^{m}}$. On constate, par exemple à l'aide du lemme \ref{alpha} ci-après (bien connu), que ce sous-groupe est un $2$-Sylow.

\medskip
On traite ensuite le cas général.

\medskip
On écrit $n=2^{m_{1}}+2^{m_{2}}+\ldots +2^{m_{r}}$ avec  $m_{1}>m_{2}>\ldots>m_{r}\geq 0$ et on pose $\mathrm{S}_{n}:=\mathrm{S}_{2^{{m}_{1}}}\times\mathrm{S}_{2^{{m}_{2}}}\times\ldots\times\mathrm{S}_{2^{{m}_{r}}}$~; la considération de la bijection évidente
$$
\coprod_{1\leq i\leq r}
\{1,2,\ldots,2^{m_{i}}\}
\hspace{4pt}\cong\hspace{4pt}
\{1,2,\ldots,n\}
$$
donne un monomorphisme $\mathfrak{S}_{2^{{m}_{1}}}\times\mathfrak{S}_{2^{{m}_{2}}}\times\ldots\times\mathfrak{S}_{2^{{m}_{r}}}\hookrightarrow\mathfrak{S}_{n}$ et conduit donc  à l'identification de $\mathrm{S}_{n}$ avec un sous-groupe de $\mathfrak{S}_{n}$. On se convainc que ce sous-groupe est un $2$-Sylow en invoquant à nouveau le lemme \ref{alpha}.

\begin{lem}\label{alpha} Soit $n\geq 0$ un entier ~; soit $\alpha(n)$ le nombre de chiffres $1$ dans l'écriture de $n$ en base $2$. Alors la valuation $2$-adique de $n!$ est $n-\alpha(n)$.
\end{lem}

\medskip
\textit{Démonstration.} On note respectivement $[x]$ et $\mathrm{v}_{2}(y)$ la partie entière d'un nombre réel $x$ et la valuation $2$-adique d'un nombre rationnel non nul $y$.

\medskip
Soit $F^{k}$, $k\in\mathbb{N}$, le sous-ensemble de $\{1,2,\ldots, n\}$ constitué des entiers $p$ avec $\mathrm{v}_{2}(p)\geq k$~; on a donc une filtration décroissante
$$
\{1,2,\ldots, n\}=F^{0}\supset F^{1}\supset\ldots\supset F^{k}\supset\ldots
$$
($F^{k}$ est vide pour $k$ assez grand). Comme le cardinal de $F^{k}$ est $[\frac{n}{2^{k}}]$, on a~:
$$
\hspace{24pt}
\mathrm{v}_{2}(n!)\hspace{4pt}=\hspace{4pt}\sum_{k=0}^{\infty}\hspace{2pt}k\hspace{2pt}(\hspace{2pt}[\frac{n}{2^{k}}]-[\frac{n}{2^{k+1}}]\hspace{2pt})\hspace{4pt}=\hspace{4pt}\sum_{k=1}^{\infty}\hspace{2pt}[\frac{n}{2^{k}}]
\hspace{24pt}
$$
(les sommes ci-dessus sont en fait finies).

\medskip
D'autre part le $k$-ième chiffre de l'écriture en base $2$ de $n$ est
$[\frac{n}{2^{k}}]-2\hspace{1pt}[\frac{n}{2^{k+1}}]$,
si bien que l'on a~;
$$
\hspace{24pt}
\alpha(n)\hspace{4pt}=\hspace{4pt}\sum_{k=0}^{\infty}\hspace{2pt}(\hspace{2pt}[\frac{n}{2^{k}}]-2\hspace{1pt}[\frac{n}{2^{k+1}}]\hspace{2pt})\hspace{4pt}=\hspace{4pt}n-\sum_{k=1}^{\infty}\hspace{2pt}[\frac{n}{2^{k}}]
\hspace{23pt}.
$$
\hfill$\square$

\subsect{Constructions quadratiques}\label{cq-1}

\medskip
\footnotesize
Nous aurons à considérer dans cette sous-section ``le'' classifiant $\mathrm{B}G$ d'un groupe fini $G$~; pour fixer les idées il sera commode de disposer d'un modèle ``fonctoriel'' en $G$.

\smallskip
Soit $G$ un groupe fini, on note $\mathrm{E}G$ le groupe simplicial défini par $\mathrm{E}_{n}G:=G^{\{0,1,\ldots,n\}}$, les faces et dégénérescences étant les applications évidentes~; on note encore $\mathrm{E}G$ la réalisation géométrique de cet ensemble simplicial. L'espace $\mathrm{E}G$ est contractile  et muni d'une action (à droite) de $G$ qui est (topologiquement) libre. On note $\mathrm{B}G$ le quotient $\mathrm{E}G/G$ ($\mathrm{E}G\to\mathrm{B}G$ est un revêtement galoisien de groupe $G$). Pour une généralisation de cette construction le lecteur pourra consulter \cite[\S 3]{Seg} (dans cette référence $G$ est un groupe topologique).

\normalsize

\medskip
Soit $G$ un groupe fini~; on rappelle que la notation $\mathfrak{S}_{2}\wr G$ désigne le produit semi-direct $\mathfrak{S}_{2}\ltimes (G\times G)$, $\mathfrak{S}_{2}$ agissant à gauche sur $G\times G$ par permutation des deux facteurs ($\mathfrak{S}_{2}\wr G$ est  le { \em produit en couronne} de $\mathfrak{S}_{2}$ et $G$, il est plus souvent noté $G\wr\mathfrak{S}_{2}$).

\medskip
On constate que l'on a
$$
\hspace{24pt}
\mathrm{B}\hspace{1pt}(\mathfrak{S}_{2}\wr G)
\hspace{4pt}=\hspace{4pt}
\mathrm{E}\mathfrak{S}_{2}\times_{\mathfrak{S}_{2}} (\mathrm{B}G\times\mathrm{B}G)
\hspace{23pt},
$$
le groupe $\mathfrak{S}_{2}$ agissant à gauche sur l'espace $\mathrm{B}G\times\mathrm{B}G$ par permutation des deux facteurs. Plus généralement~:

\bigskip
\begin{defi}\label{cq}Soit $X$ un espace topologique~; on pose
$$
\hspace{24pt}
\mathfrak{S}_{2}X
\hspace{4pt}:=\hspace{4pt}
\mathrm{E}\mathfrak{S}_{2}\times_{\mathfrak{S}_{2}} (X\times X)
\hspace{23pt},
$$
le groupe $\mathfrak{S}_{2}$ agissant à gauche sur l'espace $X\times X$ par permutation des deux facteurs. L'espace $\mathfrak{S}_{2}X$ est appelé la {\em construction quadratique} sur $X$.
\end{defi}

\bigskip
On note $\pi:\mathrm{E}\mathfrak{S}_{2}\times (X\times X)\to\mathfrak{S}_{2}X$ l'application de passage au quotient~; $\pi$~est un revêtement double.

\bigskip
\textsc{Calcul de $\mathrm{H}^{*}(\mathfrak{S}_{2}X;\mathbb{F}_{2})$}\label{cq-2}

\medskip
Ce calcul est aussi classique (voir par exemple \cite[\S 3]{Mg}, \cite[Chap.IV, \S 2]{Vo})~; il est intimement relié à la définition des opérations de Steenrod.

\medskip
On note $\mathrm{C}_{\bullet}X$ le complexe des chaînes singulières d'un espace topologique~$X$ et $\varepsilon:\mathrm{C}_{\bullet}X\to\mathbb{Z}$ son augmentation. L'action de $\mathfrak{S}_{2}$ sur l'espace $\mathrm{E}\mathfrak{S}_{2}$ fait de $\mathrm{C}_{\bullet}\hspace{1pt}\mathrm{E}\mathfrak{S}_{2}$ un complexe de $\mathbb{Z}[\mathfrak{S}_{2}]$-modules à droite et $\varepsilon:\mathrm{C}_{\bullet}\hspace{1pt}\mathrm{E}\mathfrak{S}_{2}\to\mathbb{Z}$ est une résolution libre de $\mathbb{Z}$, vu comme un $\mathbb{Z}[\mathfrak{S}_{2}]$-module à droite trivial. Pour allèger la notation on pose $\mathrm{W}:=\mathrm{C}_{\bullet}\hspace{1pt}\mathrm{E}\mathfrak{S}_{2}$.

\medskip
Le théorème d'Eilenberg-Zilber montre que l'on dispose d'une équivalence d'homotopie de $\mathbb{Z}[\mathfrak{S}_{2}]$-complexes de chaînes (fonctorielle en $X$)~:
$$
\hspace{24pt}
\mathrm{C}_{\bullet}\hspace{1pt}(\mathrm{E}\mathfrak{S}_{2}\times(X\times X))
\hspace{4pt}\cong\hspace{4pt}
\mathrm{W}\otimes_{\mathbb{Z}}\mathrm{C}_{\bullet}\hspace{1pt}(X\times X)
\hspace{23pt},
$$
le groupe $\mathfrak{S}_{2}$ agissant sur $\mathrm{C}_{\bullet}\hspace{1pt}(X\times X)$ par échange des facteurs. Le théorème d'Eilenberg-Zilber  dit encore que l'on dispose
d'une équivalence d'homotopie, disons $\mathrm{ez}:\mathrm{C}_{\bullet}\hspace{1pt}(X\times X)\cong\mathrm{C}_{\bullet}X\otimes\hspace{1pt}\mathrm{C}_{\bullet}X$ (fonctorielle en $X$ et unique à homotopie fonctorielle près). Le groupe $\mathfrak{S}_{2}$ opère aussi sur $\mathrm{C}_{\bullet}X\otimes\hspace{1pt}\mathrm{C}_{\bullet}X$ par échange des facteurs, mais il est bien connu que $\mathrm{ez}$ ne peut être $\mathfrak{S}_{2}$-équivariante car si elle l'était les opérations de Steenrod $\mathrm{Sq}^{i}$ seraient triviales pour $i>0$. Cependant~:

\begin{pro}\label{ma} Soit $X$ un espace topologique. On dispose d'une équiva\-lence d'homotopie de $\mathbb{Z}[\mathfrak{S}_{2}]$-complexes de chaînes~
$$
\hspace{24pt}
\mathrm{W}\otimes_{\mathbb{Z}}\mathrm{C}_{\bullet}\hspace{1pt}(X\times X)
\hspace{4pt}\cong\hspace{4pt}
\mathrm{W}\otimes_{\mathbb{Z}}(\mathrm{C}_{\bullet}X\otimes\mathrm{C}_{\bullet}X)
\hspace{23pt}.
$$
De plus cette équivalence d'homotopie est unique à homotopie fonctorielle près.
\end{pro}

\begin{cor}\label{Ccq} Soit $X$ un espace topologique. On dispose d'une équiva\-lence d'homotopie fonctorielle en $X$~:
$$
\hspace{24pt}
\mathrm{C}_{\bullet}\hspace{1pt}\mathfrak{S}_{2}X
\hspace{4pt}\cong\hspace{4pt}
\mathrm{W}\otimes_{\mathbb{Z}[\mathfrak{S}_{2}]}(\mathrm{C}_{\bullet}X\otimes\mathrm{C}_{\bullet}X)
\hspace{23pt},
$$
le groupe $\mathfrak{S}_{2}$ agissant sur $\mathrm{C}_{\bullet}X\otimes\mathrm{C}_{\bullet}X$ par échange des facteurs. De plus cette équivalence d'homotopie est unique à homotopie fonctorielle près.\footnote{Si l'on suppose que $X$ est un $\mathrm{CW}$-complexe et que $\mathrm{C}_{\bullet}$ désigne le complexe des chaînes cellulaires alors les complexes $\mathrm{C}_{\bullet}\hspace{1pt}\mathfrak{S}_{2}X$ et $\mathrm{W}\otimes_{\mathbb{Z}[\mathfrak{S}_{2}]}(\mathrm{C}_{\bullet}X\otimes\mathrm{C}_{\bullet}X)$ sont isomorphes~; ce point de vue est notamment adopté dans \cite{Mg} et \cite{Vo}.}
\end{cor}

\medskip
\textit{Démonstration.} La théorie des revêtements  montre que le complexe de chaînes $\mathrm{C}_{\bullet}\hspace{1pt}\mathfrak{S}_{2}X$ est le quotient de $\mathrm{C}_{\bullet}\hspace{1pt}(\mathrm{E}\mathfrak{S}_{2}\times(X\times X))$ par l'action de $\mathfrak{S}_{2}$.
$\hfill\square$

\medskip
Le lecteur pourra trouver une démonstration détaillée de la proposition \ref{ma} dans \cite[Chap. I]{Z78} qui emploie la méthode des modèles acycliques. Nous donnons ci-après, en petits caractères, un aperçu de cette démonstration.

\footnotesize
\medskip
Soit $\mathcal{T}$ la catégorie des espaces topologiques et $\Theta$ l'endofoncteur $(X_{1},X_{2})\mapsto(X_{2},X_{1})$ de la catégorie $\mathcal{T}\times\mathcal{T}$. On constate que $\Theta\circ\Theta$ est le foncteur identité de $\mathcal{T}\times\mathcal{T}$(en  d'autes termes que la catégorie $\mathcal{T}\times\mathcal{T}$ est munie d'une action du groupe  $\mathfrak{S}_{2}$)~; ceci permet de définir une catégorie $\mathfrak{S}_{2}\wr\mathcal{T}$. Quelques détails~:

\medskip
-- Les objets de $\mathfrak{S}_{2}\wr\mathcal{T}$ sont les mêmes que ceux de $\mathcal{T}\times\mathcal{T}$.

\medskip
-- Soient $Y$ et $Z$ deux objets de  $\mathcal{T}\times\mathcal{T}$,  on pose
$$
\hspace{24pt}
\mathrm{Hom}_{\mathfrak{S}_{2}\wr\mathcal{T}}\hspace{2pt}(Y,Z)
\hspace{4pt}:=\hspace{4pt}
\mathrm{Hom}_{\mathcal{T}\times\mathcal{T}}\hspace{2pt}(Y,Z)
\coprod
\mathrm{Hom}_{\mathcal{T}\times\mathcal{T}}\hspace{2pt}(\Theta Y,Z)
\hspace{23pt};
$$
le lecteur devinera sans peine la règle de composition des morphismes de $\mathfrak{S}_{2}\wr\mathcal{T}$ et obser\-vera, chemin faisant, que l'application $\mathrm{Hom}_{\mathfrak{S}_{2}\wr\mathcal{T}}\hspace{2pt}(Y,Z)\to\mathfrak{S}_{2}$ qui envoie le premier terme sur $\mathrm{id}$ et le second sur  la transposition $\tau_{1,2}$ induit un foncteur $\mathrm{p}:\mathfrak{S}_{2}\wr\mathcal{T}\to\underline{\mathfrak{S}_{2}}$,  $\underline{\mathfrak{S}_{2}}$ désignant la catégorie associée au groupe $\mathfrak{S}_{2}$\footnote{La catégorie $\underline{\mathfrak{S}_{2}}$ possède un seul objet, disons $*$, et le monoïde $\mathrm{Hom}_{\underline{\mathfrak{S}_{2}}}(*,*)$ est le groupe $\mathfrak{S}_{2}$.}

\medskip
-- L'inclusion  $\mathrm{Hom}_{\mathcal{T}\times\mathcal{T}}\hspace{2pt}(Y,Z)\hookrightarrow\mathrm{Hom}_{\mathfrak{S}_{2}\wr\mathcal{T}}\hspace{2pt}(Y,Z)$ induit un foncteur, qui est l'identité sur les objets, de $\mathcal{T}\times\mathcal{T}$ dans  $\mathfrak{S}_{2}\wr\mathcal{T}$~; on le note $\mathrm{i}$.

\smallskip
-- L'inclusion  $\mathrm{Hom}_{\mathcal{T}\times\mathcal{T}}\hspace{2pt}(\Theta Y,\Theta Y)\hookrightarrow\mathrm{Hom}_{\mathfrak{S}_{2}\wr\mathcal{T}}\hspace{2pt}(Y,\Theta Y)$ envoie l'identité de $\Theta Y$ dans la catégorie  $\mathcal{T}\times\mathcal{T}$ sur un élément que l'on note $\iota(Y)$. On constate que l'application
$$
\mathrm{Hom}_{\mathcal{T}\times\mathcal{T}}\hspace{2pt}(\Theta Y,Z)\to\mathrm{Hom}_{\mathfrak{S}_{2}\wr\mathcal{T}}\hspace{2pt}(Y,Z)\hspace{12pt},\hspace{12pt}
f \mapsto\mathrm{i}(f)\circ\iota(Y)
$$
est l'inclusion naturelle te que le composé $\iota(\Theta Y)\circ\iota(Y)$ est l'identité de $Y$.

\medskip
Après avoir mis en place le formalisme ci-dessus, on considère les foncteurs
$$
\hspace{4pt}
\mathrm{F}_{1}(X_{1},X_{2})=\mathrm{W}\otimes_{\mathbb{Z}}\mathrm{C}_{\bullet}\hspace{1pt}(X_{1}\times X_{2})
\hspace{12pt},\hspace{12pt}
\mathrm{F}_{2}(X_{1},X_{2})=\mathrm{W}\otimes_{\mathbb{Z}}(\mathrm{C}_{\bullet}\hspace{1pt}(X_{1})\otimes\mathrm{C}_{\bullet}\hspace{1pt}(X_{2}))
\hspace{11pt}.
$$
définis sur $\mathcal{T}\times\mathcal{T}$ et à valeurs dans la catégorie dans la catégorie des $\mathbb{Z}[\mathfrak{S}_{2}]$-complexes de chaînes (l'action de $\mathfrak{S}_{2}$ sur $\mathrm{C}_{\bullet}\hspace{1pt}(X_{1}\times X_{2})$ et $\mathrm{C}_{\bullet}\hspace{1pt}(X_{1})\otimes\mathrm{C}_{\bullet}\hspace{1pt}(X_{2})$ étant triviale).

\medskip
On montre que ces foncteurs se prolongent canoniquement en des foncteurs, disons $\widehat{\mathrm{F}}_{1}$ et $\widehat{\mathrm{F}}_{2}$, définis sur la catégorie $\mathfrak{S}_{2}\wr\mathcal{T}$. La valeur de $\widehat{\mathrm{F}}_{1}$  (resp. $\widehat{\mathrm{F}}_{2}$) sur $\iota(X_{1},X_{2})$ est le produit tensoriel de l'action de $\mathfrak{S}_{2}$ sur $\mathrm{W}$ et de l'isomorphisme canonique $\mathrm{C}_{\bullet}\hspace{1pt}(X_{1}\times X_{2})\cong\mathrm{C}_{\bullet}\hspace{1pt}(X_{2}\times X_{1})$ (resp. $\mathrm{C}_{\bullet}\hspace{1pt}(X_{1})\otimes\mathrm{C}_{\bullet}\hspace{1pt}(X_{2})\cong\mathrm{C}_{\bullet}\hspace{1pt}(X_{2})\otimes\mathrm{C}_{\bullet}\hspace{1pt}(X_{1})$).

\medskip
On vérifie que les foncteurs $\widehat{\mathrm{F}}_{1}$ et $\widehat{\mathrm{F}}_{2}$, sont libres sur les ``modèles'' $\{(\Delta^{n_{1}},\Delta^{n_{2}})\}$, $(n_{1},n_{2})$ parcourant $\mathbb{N}\times\mathbb{N}$ ($\Delta^{n}$ désigne ici le $n$-simplexe standard),  et que leurs versions augmentées (par le foncteur constant de valeur $\mathbb{Z}$) sont acycliques sur ces mêmes modèles. Le théorème des modèles acycliques fournit alors l'énoncé suivant~:

\begin{pro}\label{78} Soient $X_{1}$ et $X_{2}$ deux espaces topologiques. On dispose d'une équivalence d'homotopie fonctorielle en $(X_{1},X_{2})$, vu comme un objet de la catégorie $\mathfrak{S}_{2}\wr\mathcal{T}$,
$$
\hspace{24pt}
\widehat{\mathrm{ez}}_{\hspace{0.75pt}(X_{1},X_{2})}
\hspace{4pt}:\hspace{4pt}
\mathrm{W}\otimes_{\mathbb{Z}}\mathrm{C}_{\bullet}\hspace{1pt}(X_{1}\times X_{2})
\hspace{4pt}\cong\hspace{4pt}
\mathrm{W}\otimes_{\mathbb{Z}}(\mathrm{C}_{\bullet}\hspace{1pt}(X_{1})\otimes\mathrm{C}_{\bullet}\hspace{1pt}(X_{2}))
\hspace{23pt}.
$$
De plus cette équivalence d'homotopie est unique à homotopie fonctorielle (au même sens que ci-dessus) près.
\end{pro}

\medskip
\begin{rem}\label{ez}  La catégorie $\mathcal{T}\times\mathcal{T}$ s'identifie à une sous-catégorie (non pleine) de $\mathfrak{S}_{2}\wr\mathcal{T}$~;  toujours d'après le théorème des modèles acycliques, la ``restriction'' de la transformation naturelle $\widehat{\mathrm{ez}}$ à $\mathcal{T}\times\mathcal{T}$ est homotope, (fonctoriellement en $\mathcal{T}\times\mathcal{T}$) à la transformation naturelle $1_{\mathrm{W}}\otimes_{\mathbb{Z}}\mathrm{ez}$ ($\mathrm{ez}$ désignant la transformation naturelle d'Eilenberg-Zilber).
\end{rem}

\medskip
Soit $X$ un espace topologique~; il découle de la proposition  \ref{78} que les équivalences d'homotopie
$$
\hspace{24pt}
\widehat{\mathrm{ez}}_{\hspace{0.75pt}(X,X)}
\hspace{4pt}:\hspace{4pt}
\mathrm{W}\otimes_{\mathbb{Z}}\mathrm{C}_{\bullet}\hspace{1pt}(X\times X)
\hspace{4pt}\cong\hspace{4pt}
\mathrm{W}\otimes_{\mathbb{Z}}(\mathrm{C}_{\bullet}\hspace{1pt}(X)\otimes\mathrm{C}_{\bullet}\hspace{1pt}(X))
\hspace{24pt}
$$
induisent une transformation naturelle de foncteurs définis sur la catégorie des espaces topologiques  et à valeurs dans la catégorie des $\mathbb{Z}[\mathfrak{S}_{2}]$-complexes de chaînes. La proposition \ref{ma} en résulte.

\medskip
\begin{rem}
La catégorie $\underline{\mathfrak{S}_{2}}\times\mathcal{T}$ s'identifie à une sous-catégorie (non pleine) de $\mathfrak{S}_{2}\wr\mathcal{T}$. La transformation naturelle évoquée ci-dessus est simplement la restriction de $\widehat{\mathrm{ez}}$ à $\underline{\mathfrak{S}_{2}}\times\mathcal{T}$.
\end{rem}

\normalsize

\bigskip
La proposition \ref{Ccq} implique~:

\begin{cor}\label{cqevh} Soit $X$ un espace topologique. On a un isomorphisme canonique de $\mathbb{F}_{2}$-espaces vectoriels gradués
$$
\hspace{24pt}
\mathrm{H}_{*}(\mathfrak{S}_{2}X;\mathbb{F}_{2})
\hspace{4pt}\cong\hspace{4pt}
\mathrm{H}_{*}(\mathfrak{S}_{2};\mathrm{H}_{*}(X;\mathbb{F}_{2})\otimes\mathrm{H}_{*}(X;\mathbb{F}_{2}))
\hspace{23pt},
$$
le groupe $\mathfrak{S}_{2}$ agissant sur $\mathrm{H}_{*}(X;\mathbb{F}_{2})\otimes\mathrm{H}_{*}(X;\mathbb{F}_{2})$ par échange des facteurs.
\end{cor}

\medskip
Précisons la notation. Soient $G$ un groupe et $M_{*}=(M_{n})_{n\in\mathbb{N}}$ un $\mathbb{Z}[G]$-module gradué~; $\mathrm{H}_{*}(G;M_{*})$ désigne le groupe abélien gradué dont le $n$-ième terme est $\bigoplus_{p+q=n}\mathrm{H}_{p}(G;M_{q})$. On définit $\mathrm{H}^{*}(G;M^{*})$ \textit{mutatis mutandis}.

\bigskip
\textit{Démonstration.} On note $\mathrm{H}_{\bullet}(X;\mathbb{F}_{2})$ le complexe de chaînes à différentielle nulle
$$
\hspace{24pt}
\mathrm{H}_{0}(X;\mathbb{F}_{2})\overset{0}{\leftarrow}\mathrm{H}_{1}(X;\mathbb{F}_{2})\overset{0}{\leftarrow}\mathrm{H}_{2}(X;\mathbb{F}_{2})\overset{0}{\leftarrow}\ldots\overset{0}{\leftarrow}\mathrm{H}_{n}(X;\mathbb{F}_{2})\overset{0}{\leftarrow}\ldots
\hspace{23pt}.
$$
Comme $\mathrm{C}_{\bullet}(X;\mathbb{F}_{2})$ est un complexe de chaînes sur un corps il existe une équivalence d'homotopie $h:\mathrm{C}_{\bullet}(X;\mathbb{F}_{2})\to\mathrm{H}_{\bullet}(X;\mathbb{F}_{2})$ induisant l'identité en homologie~; de plus un tel $h$ est unique à homotopie près. Ces deux propriétés (fort classiques) résultent par exemple de l'équivalence $(iii)\Leftrightarrow(iv)$ de \cite[\S 2, \numero\hspace{3pt}5, Proposition 6]{Boh}. Précisons un peu. Posons $\mathrm{C}_{\bullet}(X;\mathbb{F}_{2})=\mathrm{C}_{\bullet}$ et $\mathrm{H}_{\bullet}(X;\mathbb{F}_{2})=\mathrm{H}_{\bullet}$. La condition $(iv)$ évoquée ci-dessus dit en particulier que $\mathrm{C}_{\bullet}$ est isomorphe à un complexe de la forme $\mathrm{H}_{\bullet}\oplus\mathrm{D}_{\bullet}$ avec $\mathrm{D}_{\bullet}$ homotope à $0$. Pour un complexe de cette forme, l'existence de $h$ et son unicité à homotopie près sont évidentes.

\medskip
On invoque alors un célèbre lemme de Steenrod \cite[Lemma 5.2]{St} dont nous rappelons l'énoncé dans le contexte qui nous intéresse~:

\begin{lem}\label{2-St} Soient $f_{0},f_{1}:C_{\bullet}\to D_{\bullet}$ deux homomorphismes de complexes de chaînes. Si $f_{0}$ et $f_{1}$ sont homotopes alors il en est de même pour
\footnotesize
$$
\hspace{2pt}
1_{\mathrm{W}}\otimes_{\mathbb{Z}[\mathfrak{S}_{2}]}(f_{0}\otimes f_{0})\hspace{1pt},\hspace{1pt}1_{\mathrm{W}}\otimes_{\mathbb{Z}[\mathfrak{S}_{2}]}(f_{1}\otimes f_{1}):\mathrm{W}\otimes_{\mathbb{Z}[\mathfrak{S}_{2}]}(C_{\bullet}\otimes C_{\bullet})\to\mathrm{W}\otimes_{\mathbb{Z}[\mathfrak{S}_{2}]}(D_{\bullet}\otimes D_{\bullet})
\hspace{1pt}.
$$
\normalsize
\end{lem}

\begin{scho}\label{2-St-scho}Soit $f:C_{\bullet}\to D_{\bullet}$ un homomorphisme de complexes de chaînes. Si $f$ est une équivalence d'homotopie alors il en est de même pour
$$
\hspace{18pt}
1_{\mathrm{W}}\otimes_{\mathbb{Z}[\mathfrak{S}_{2}]}(f\otimes f)\hspace{2pt}:\hspace{2pt}\mathrm{W}\otimes_{\mathbb{Z}[\mathfrak{S}_{2}]}(C_{\bullet}\otimes C_{\bullet})\to\mathrm{W}\otimes_{\mathbb{Z}[\mathfrak{S}_{2}]}(D_{\bullet}\otimes D_{\bullet})
\hspace{17pt}.
$$
\end{scho}

Le lemme de Steenrod et son scholie montrent que l'homomorphisme de complexes de chaînes
\footnotesize
$$
\begin{CD}
\mathrm{W}\otimes_{\mathbb{Z}[\mathfrak{S}_{2}]}(\mathrm{C}_{\bullet}(X;\mathbb{F}_{2})\otimes\mathrm{C}_{\bullet}(X;\mathbb{F}_{2}))
@>1_{\mathrm{W}}\otimes_{\mathbb{Z}[\mathfrak{S}_{2}]}(h\otimes h)>>
\mathrm{W}\otimes_{\mathbb{Z}[\mathfrak{S}_{2}]}(\mathrm{H}_{\bullet}(X;\mathbb{F}_{2})\otimes\mathrm{H}_{\bullet}(X;\mathbb{F}_{2}))
\end{CD}
$$
\normalsize
est une équivalence d'homotopie et que la classe d'homotopie de cet homomorphisme est indépendante du choix de $h$.
\hfill$\square$

\medskip
\begin{rem}Il est implicite dans \cite{St} que l'on dispose de versions équivariantes de \ref{2-St} et \ref{2-St-scho} dans lesquelles les complexes sont remplacés par des $\mathbb{Z}[\mathfrak{S}_{2}]$-complexes et le produit tensoriel $\otimes_{\mathbb{Z}[\mathfrak{S}_{2}]}$ par le produit tensoriel~$\otimes_{\mathbb{Z}}$. La présence de $\mathrm{W}$ dans la théorie de Steenrod est nécessaire~: Soit $\mathrm{E}_{\bullet}$ le complexe défini par $\mathrm{E}_{n}=\mathbb{F}_{2}$ pour $n=0,1$, $\mathrm{E}_{n}=0$ pour $n\not=0,1$ et $\mathrm{d}_{1}=1$~; $\mathrm{E}_{\bullet}$ est homotope à $0$ mais $\mathrm{E}_{\bullet}\otimes\mathrm{E}_{\bullet}$ n'est pas homotope à $0$ en tant que $\mathbb{Z}[\mathfrak{S}_{2}]$-complexe.
\end{rem}

\bigskip
Dualement~:

\begin{cor}\label{cqev} Soit $X$ un espace topologique avec $\mathrm{H}^{*}(X;\mathbb{F}_{2})$ de dimension finie en chaque degré. On a un isomorphisme canonique de $\mathbb{F}_{2}$-espaces vectoriels gradués
$$
\hspace{24pt}
\mathrm{H}^{*}(\mathfrak{S}_{2}X;\mathbb{F}_{2})
\hspace{4pt}\cong\hspace{4pt}
\mathrm{H}^{*}(\mathfrak{S}_{2};\mathrm{H}^{*}(X;\mathbb{F}_{2})\otimes\mathrm{H}^{*}(X;\mathbb{F}_{2}))
\hspace{23pt},
$$
le groupe $\mathfrak{S}_{2}$ agissant sur $\mathrm{H}^{*}(X;\mathbb{F}_{2})\otimes\mathrm{H}^{*}(X;\mathbb{F}_{2})$ par échange des facteurs.
\end{cor}

\medskip
(L'hypothèse de finitude est juste là pour assurer que le $\mathbb{F}_{2}$-espace vectoriel $\mathrm{H}^{s}(X;\mathbb{F}_{2})\otimes\mathrm{H}^{t}(X;\mathbb{F}_{2})$ est le dual de $\mathrm{H}_{s}(X;\mathbb{F}_{2})\otimes\mathrm{H}_{t}(X;\mathbb{F}_{2})$ pour tout couple d'entiers $(s,t)$.)

\bigskip
\textsc{Notation.} A partir de maintenant la cohomologie que nous considèrerons sera la cohomologie à coefficients dans $\mathbb{F}_{2}$,  aussi nous abrègerons la notation $\mathrm{H}^{*}(-;\mathbb{F}_{2})$ en $\mathrm{H}^{*}(-)$ ou $\mathrm{H}^{*}-$.

\bigskip
\begin{scho}\label{cqB} Soit $X$ un espace topologique avec $\mathrm{H}^{n}X$ de dimension finie pour tout $n$. Soit $B_{n}$ une base (totalement) ordonnée de $\mathrm{H}^{n}X$~; on pose $B:=\coprod_{n\in\mathbb{N}}B_{n}$.

\smallskip
On munit $B$ de la relation d'ordre (total) qui prolonge celles des $B_{n}$ et qui vérifie en outre $x<y$ pour $x\in B_{n}$, $y\in B_{p}$ et $n<p$.

\smallskip
On a un isomorphisme canonique de $\mathbb{F}_{2}$-espaces vectoriels gradués
$$
\hspace{24pt}
\mathrm{H}^{*}\mathfrak{S}_{2}X
\hspace{4pt}\cong\hspace{4pt}
\bigoplus_{(x,y)\in B\times B\hspace{1pt},\hspace{1pt} x<y}
\mathbb{F}_{2}\hspace{2pt}.\hspace{2pt}(x\otimes y+y\otimes x)
\hspace{6pt}\oplus\hspace{6pt}
\bigoplus_{x\in B}
\mathrm{H}^{*}\mathfrak{S}_{2}\hspace{2pt}.\hspace{2pt}x\otimes x
\hspace{23pt}.
$$
\end{scho}

\medskip
\textit{Démonstration.} Le $\mathbb{F}_{2}[\mathfrak{S}_{2}]$-module ($\mathbb{N}$-gradué) $\mathrm{H}^{*}X\otimes\mathrm{H}^{*}X$ est somme directe des sous-modules suivants~:

\smallskip
-- le sous-module engendré par $x\otimes y$ et $y\otimes x$ avec $(x,y)\in B\times B\hspace{1pt},\hspace{1pt} x<y$~;

\smallskip
-- le sous-module engendré par $x\otimes x$ avec $x\in B$.

\smallskip
Oublions la graduation~; le premier est isomorphe à $\mathbb{F}_{2}[\mathfrak{S}_{2}]$, le second à  $\mathbb{F}_{2}$ (muni de l'action triviale de $\mathfrak{S}_{2}$). La formule pour $\mathrm{H}^{*}\mathfrak{S}_{2}X$ en résulte.
\hfill$\square$

\bigskip
Nous nous proposons maintenant de ``polir'' la formule du scholie \ref{cqB} en retravaillant la définition des deux termes du second membre.

\medskip
1) On considère le revêtement double $\pi:\mathrm{E}\mathfrak{S}_{2}\times (X\times X)\to\mathfrak{S}_{2}X$ et l'homomorphisme de transfert qui lui est associé~:
$$
\hspace{24pt}
\mathrm{tr}:\mathrm{H}^{*}(\mathrm{E}\mathfrak{S}_{2}\times (X\times X))=\mathrm{H}^{*}X\otimes\mathrm{H}^{*}X\to\mathrm{H}^{*}\mathfrak{S}_{2}X
\hspace{23pt}.
$$
Soient $z_{1}$ et  $z_{2}$ deux éléments de $\mathrm{H}^{*}X$, on constate (invoquer \ref{ez}) que l'on a $\mathrm{tr}(z_{1}\otimes z_{2})=z_{1}\otimes z_{2}+z_{2}\otimes z_{1}$ (et en particulier $\mathrm{tr}(z_{1}\otimes z_{2})=0$ pour $z_{1}=z_{2}$). D'où l'égalité~:
$$
\bigoplus_{(x,y)\in B\times B\hspace{1pt},\hspace{1pt} x<y}
\mathbb{F}_{2}\hspace{2pt}.\hspace{2pt}(x\otimes y+y\otimes x)
\hspace{4pt}=\hspace{4pt}
\mathop{\mathrm{im}}\mathrm{tr}
$$
(la notation $\mathop{\mathrm{im}}\mathrm{tr}$ désigne ci-dessus l'image de l'homomorphisme $\mathrm{tr}$).

\medskip
2)  Soient $A$ un groupe abélien et $n$ un entier~; on note $A[n]$ le complexe de chaînes dont le $n$-ième terme est $A$ et dont tous les autres sont nuls.

\smallskip
Soient $z$ un élément de $\mathrm{H}^{n}X$ et $\tilde{z}:\mathrm{C}_{n}X\to\mathbb{F}_{2}$ un cocycle représentant $z$~; on rappelle que $\tilde{z}$ s'identifie à une morphisme de complexes $\mathrm{C}_{\bullet}X\to\mathbb{F}_{2}[n]$ et $z$ à la classe d'homotopie de ce morphisme.

\smallskip
Le morphisme de complexes
$$
\begin{CD}
\mathrm{W}\otimes\mathrm{C}_{\bullet}X\otimes\mathrm{C}_{\bullet}X
@>\varepsilon\hspace{1pt}\otimes\hspace{1pt}\tilde{z}\hspace{1pt}\otimes\hspace{1pt}\tilde{z}>>
\mathbb{Z}[0]\otimes\mathbb{F}_{2}[n]\otimes\mathbb{F}_{2}[n]=\mathbb{F}_{2}[2n]
\end{CD}
$$
se factorise à travers $\mathrm{W}\otimes_{\mathbb{Z}[\mathfrak{S}_{2}]}(\mathrm{C}_{\bullet}X\otimes\mathrm{C}_{\bullet}X)$. Compte tenu de \ref{Ccq}, il définit un élément de $\mathrm{H}^{2n}\hspace{0.5pt}\mathfrak{S}_{2}X$~; le lemme ci-dessous est encore dû à Steenrod.

\begin{lem-def}\label{P2} L'élément de $\mathrm{H}^{2n}\hspace{0.5pt}\mathfrak{S}_{2}X$ introduit ci-dessus est indépendant du choix du cocycle $\tilde{z}$ représentant $z$~; on le note $\mathrm{P}_{\hspace{-2pt}2}\hspace{1pt}z$ et on l'appelle la {\em deuxième puissance de Steenrod} de $z$.
\end{lem-def}

\medskip
\textit{Démonstration.} Conséquence de \ref{2-St} avec $C_{\bullet}=\mathrm{C}_{\bullet}X$ et $D_{\bullet}=\mathbb{F}_{2}[n]$.
\hfill$\square$

\medskip
Soient $z_{1}$ et $z_{2}$ deux  éléments de $\mathrm{H}^{n}X$, on constate que l'on a $\mathrm{P}_{\hspace{-2pt}2}\hspace{1pt}(z_{1}+z_{2})=\mathrm{P}_{\hspace{-2pt}2}\hspace{1pt}z_{1}+\mathrm{P}_{\hspace{-2pt}2}\hspace{1pt}z_{2}+\mathrm{tr}(z_{1}\otimes z_{2})$~:  l'application $\mathrm{P}_{\hspace{-2pt}2}:\mathrm{H}^{n}X\to\mathrm{H}^{2n}\hspace{0.5pt}\mathfrak{S}_{2}X$ est quadratique.

\medskip
La définition même de $\mathrm{P}_{\hspace{-2pt}2}$ montre que le terme $\bigoplus_{x\in B}\mathrm{H}^{*}\mathfrak{S}_{2}\hspace{2pt}.\hspace{2pt}x\otimes x$ peut être réécrit $\bigoplus_{x\in B}\mathrm{H}^{*}\mathfrak{S}_{2}\hspace{2pt}.\hspace{2pt}\mathrm{P}_{\hspace{-2pt}2}\hspace{1pt}x$.

\bigskip
\textsc{La diagonale de Steenrod et le foncteur $\mathrm{R}_{1}$ de Singer}

\medskip
L'application $\mathrm{id}\times\delta:\mathrm{E}\mathfrak{S}_{2}\times X\to\mathrm{E}\mathfrak{S}_{2}\times (X\times X)$, $\delta$ désignant la diagonale de $X$, induit par passage au quotient une application $\Delta:\mathrm{B}\mathfrak{S}_{2}\times X\to\mathfrak{S}_{2}X$ que l'on appelle la {\em diagonale de Steenrod}.

\smallskip
La contemplation du diagramme de revêtements doubles
$$
\begin{CD}
\mathrm{E}\mathfrak{S}_{2}\times X
@>\mathrm{id}\times\delta>>
\mathrm{E}\mathfrak{S}_{2}\times (X\times X) \\
\hspace{12pt}@VVV  @VV\pi V \\
\hspace{5pt}\mathrm{B}\mathfrak{S}_{2}\times X
@>\Delta>>
\mathfrak{S}_{2}X
\end{CD}
$$
montre que l'homomorphisme composé
$$
\begin{CD}
\mathrm{H}^{*}X\otimes\mathrm{H}^{*}X
@>\mathrm{tr}>>
\mathrm{H}^{*}\mathfrak{S}_{2}X
@>\Delta^{*}>>
\mathrm{H}^{*}(\mathrm{B}\mathfrak{S}_{2}\times X)
\end{CD}
$$
est nul (observer que l'homomorphisme $\mathrm{tr}:\mathrm{H}^{0}\hspace{1pt}\mathrm{B}\mathfrak{S}_{2}\to:\mathrm{H}^{0}\hspace{1pt}\mathrm{E}\mathfrak{S}_{2}$ est nul), d'où l'égalité $\Delta^{*}(\mathop{\mathrm{im}}\mathrm{tr})=0$.

\medskip
On a $\mathrm{H}^{*}\mathfrak{S}_{2}=\mathrm{H}^{*}\mathbb{Z}/2=\mathbb{F}_{2}[u]$,  $u$ désignant l'élément non nul de $\mathrm{H}^{1}\mathfrak{S}_{2}$~; on a donc $\mathrm{H}^{*}(\mathrm{B}\mathfrak{S}_{2}\times X)=\mathbb{F}_{2}[u]\otimes\mathrm{H}^{*}X$. Soit $z$ un élément de $\mathrm{H}^{n}X$, on rappelle que l'on peut définir les opérations de Steenrod par la formule ci-dessous~:
$$
\hspace{24pt}
\Delta^{*}\mathrm{P}_{\hspace{-2pt}2}\hspace{1pt}z
\hspace{4pt}=\hspace{4pt}
\sum_{i=0}^{n}u^{n-i}\otimes\mathrm{Sq}^{i}z
\hspace{23pt}.
$$
On note que l'expression au second membre de  cette égalité ne fait intervenir que la structure de $\mathrm{A}$-module instable de $\mathrm{H}^{*}X$~; cette observation conduit à la définition ci-après.

\medskip
\begin{defi}\label{St1} Soient $M$ un $\mathrm{A}$-module instable et $z$ un élément (homogène) de $M$~; on note $\mathrm{St}_{1}\hspace{1pt}z$ l'élément
$\sum_{i=0}^{\vert z\vert}u^{\vert z\vert-i}\otimes\mathrm{Sq}^{i}z$ de $\mathbb{F}_{2}[u]\otimes M$.\footnote{L'indice $1$ est là parce qu'il existe des $\mathrm{St}_{s}\hspace{1pt}z$ avec $s\in\mathbb{N}$, voir par exemple \cite{LZdest}.}

\smallskip
On  observera que l'application $\mathbb{F}_{2}$-linéaire
$$
\mathrm{St}_{1}:
M\to\mathbb{F}_{2}[u]\otimes M\hspace{6pt},\hspace{6pt}z\mapsto\mathrm{St}_{1}\hspace{1pt}z
$$
``multiplie le degré par $2$''.
\end{defi}

\medskip
\begin{pro-def}\label{R1}Soit $M$ un $\mathrm{A}$-module instable~; on note $\mathrm{R}_{1}M$ le sous-$\mathbb{F}_{2}[u]$-module de $\mathbb{F}_{2}[u]\otimes M$ engendré par $\mathrm{St}_{1}M$.

\medskip
{\em (a)} Le module $\mathrm{R}_{1}M$ est un sous-$\mathrm{A}$-module de $\mathbb{F}_{2}[u]\otimes M$, en particulier $\mathrm{R}_{1}M$ est un $\mathrm{A}$-module instable.

\medskip
{\em (b)} La correspondance $M\mapsto\mathrm{R}_{1}M$ s'étend en un endofoncteur $\mathrm{R}_{1}:\mathcal{U}\circlearrowleft$.

\medskip
{\em (c)} Le foncteur $\mathrm{R}_{1}$ préserve les sommes directes.

\medskip
{\em (d)} Soit $\mathcal{O}:\mathcal{K}\to\mathcal{U}$ le foncteur oubli. Si $M$ est une $\mathrm{A}$-algèbre instable alors $\mathrm{R}_{1}\mathcal{O}M$ est une  sous-$\mathrm{A}$-algèbre instable de la $\mathrm{A}$-algèbre instable $\mathbb{F}_{2}[u]\otimes M$. La correspondance $M\mapsto\mathrm{R}_{1}M$ s'étend en un endofoncteur $\mathrm{R}_{1}:\mathcal{K}\circlearrowleft$ tel que le diagramme de foncteurs
$$
\begin{CD}
\mathcal{K}@>\mathrm{R}_{1}>>\mathcal{K} \\
@V\mathcal{O}VV @V\mathcal{O}VV \\
\mathcal{U}@>\mathrm{R}_{1}>>\mathcal{U}
\end{CD}
$$
est commutatif.

\medskip
{\em (e)} Soit $B\subset M$ une base (au sens gradué) du $\mathbb{F}_{2}$-espace vectoriel $\mathbb{N}$-gradué sous-jacent à $M$, alors $\mathrm{St}_{1}B$ est une base (au sens gradué) du $\mathbb{F}_{2}[u]$-module $\mathbb{N}$-gradué sous-jacent à $\mathrm{R}_{1}M$.

\medskip
{\em (f)} Soit $E$ un $\mathbb{F}_{2}$-espace vectoriel $\mathbb{N}$-gradué~;  on note $\Phi E$ le $\mathbb{F}_{2}$-espace vectoriel $\mathbb{N}$-gradué défini par
$$
(\Phi E)^{n}=
\begin{cases}
E^{\frac{n}{2}} & \text{pour $n$ pair,} \\
0 &  \text{pour $n$ impair.}
\end{cases}
$$
Soit $\mathcal{E}$ la catégorie des $\mathbb{F}_{2}$-espaces vectoriels $\mathbb{N}$-gradués~; on a un $\mathcal{E}$-isomorphisme naturel en $M$~:
$$
\hspace{24pt}
\mathcal{O\hspace{1pt}}\mathrm{R}_{1}M
\hspace{4pt}\cong\hspace{4pt}
\mathbb{F}_{2}[u]\otimes\Phi\hspace{1pt}\mathcal{O}M
\hspace{23pt},
$$
$\mathcal{O}$ désignant cette fois le foncteur oubli de $\mathcal{U}$ vers $\mathcal{E}$.

\medskip
{\em (g)} Le foncteur $\mathrm{R}_{1}$ est exact.

\end{pro-def}

\medskip
\textit{Démonstration.} Pour les points (a) et (e) nous renvoyons à \cite{LZdest}. Les points (b) et (c) sont évidents. Le point (f) résulte du point  (e) et du fait que le degré de $\mathrm{St}_{1}\hspace{1pt}z$ est le double de celui de $z$. Le point (f) implique le point (g) (une $\mathcal{U}$-suite est exacte si et seulement si la $\mathcal{E}$-suite sous-jacente est exacte).

Soient $M$ une $\mathrm{A}$-algèbre instable, $z_{1},z_{2}$ deux éléments de $M$ et $i$ un entier naturel~;  l'égalité
$\mathrm{Sq}^{i}(z_{1}z_{2})=\sum_{j+k=i}\mathrm{Sq}^{j}z_{1}\hspace{2pt}\mathrm{Sq}^{k}z_{2}$ est équivalente à l'égalité $\mathrm{St}_{1}(z_{1}z_{2})=\mathrm{St}_{1}\hspace{1pt}z_{1}\hspace{2pt}\mathrm{St}_{1}\hspace{1pt}z_{2}$, d'où le point (d).
\hfill$\square$

\bigskip
\textsc{Les endofoncteurs $\Phi:\mathcal{U}\circlearrowleft$ et $\Phi:\mathcal{K}\circlearrowleft$}

\medskip
On a défini dans le point (f) de \ref{R1} un endofoncteur ``double''  $\Phi:\mathcal{E}\circlearrowleft$~; on va définir ci-dessous des endofoncteurs $\Phi:\mathcal{U}\circlearrowleft$ et $\Phi:\mathcal{K}\circlearrowleft$ compatibles (en un sens évident) avec les foncteurs oubli $\mathcal{K}\to\mathcal{U}\to\mathcal{E}$.

\medskip
Soit $M$ un $\mathbb{F}_{2}$-espace vectoriel $\mathbb{N}$-gradué~; on note $\sigma:M\otimes M\circlearrowleft $ l'automorphisme involutif $x\otimes y\mapsto y\otimes x$ définissant l'action  de $\mathfrak{S}_{2}$ sur $M\otimes M$. On considère le $0$-ième groupe de cohomologie de Tate $\widehat{\mathrm{H}}^{0}(\mathfrak{S}_{2};M\otimes M)$  (``les invariants divisés par les normes'')~:
$$
\hspace{24pt}
\widehat{\mathrm{H}}^{0}(\mathfrak{S}_{2};M\otimes M)
\hspace{4pt}:=\hspace{4pt}
\ker(1-\sigma)/\mathop{\mathrm{im}}(1+\sigma)
\hspace{23pt}.
$$
On constate que l'application
$$
M\to\widehat{\mathrm{H}}^{0}(\mathfrak{S}_{2};M\otimes M)
\hspace{12pt},\hspace{12pt}
x\mapsto\text{classe de\hspace{4pt}} x\otimes x
$$
est un isomorphisme ``multipliant le degré par 2'', en d'autres termes que l'on a un $\mathcal{E}$-isomorphisme naturel $\Phi M\cong \widehat{\mathrm{H}}^{0}(\mathfrak{S}_{2};M\otimes M)$.

\medskip
Si $M$ est un $\mathrm{A}$-module instable alors  $\sigma:M\otimes M\circlearrowleft$ est un $\mathcal{U}$-automorphisme si bien que $\Phi M\cong \widehat{\mathrm{H}}^{0}(\mathfrak{S}_{2};M\otimes M)$ est muni d'une structure naturelle de  $\mathrm{A}$-module instable. Il n'est pas difficile d'expliciter l'action des opérations de Steenrod sur $\Phi M$~:
$$
\mathrm{Sq}^{i}\hspace{1pt}\Phi x
=
\begin{cases}
\Phi\hspace{1pt}\mathrm{Sq}^{\frac{i}{2}}x & \text{pour $i$ pair,} \\
0 &  \text{pour $n$ impair,}
\end{cases}
$$
$\Phi x$ désignant la classe de $x\otimes x$.

\medskip
Si $M$ est une $\mathrm{A}$-algèbre instable alors $\ker(1-\sigma)=(M\otimes M)^{\mathfrak{S}_{2}}$ est une sous- $\mathrm{A}$-algèbre instable de $M\otimes M$ et  $\mathop{\mathrm{im}}(1+\sigma)$ et un idéal de $(M\otimes M)^{\mathfrak{S}_{2}}$ stable sous l'action $\mathrm{A}$ si bien que $\Phi M$ est muni d'une $\mathcal{K}$-structure naturelle. On constate que le produit de $\Phi M$ est simplement le ``double'' de celui de $M$~: $\Phi x\hspace{1pt}\Phi y=\Phi (xy)$.

\medskip
\begin{rem}\label{lambda2} Soit $M$ un $\mathbb{F}_{2}$-espace vectoriel $\mathbb{N}$-gradué ou un $\mathrm{A}$-module instable, on constate que $\mathop{\mathrm{im}}(1+\sigma)$ est naturellement isomorphe au quotient de $M\otimes M$ par le sous-objet engendré par  les $x\otimes x$, $x$ parcourant $M$~; il est donc raisonnable de poser $\mathop{\mathrm{im}}(1+\sigma):=\Lambda^{2}M$.
\end{rem}

\bigskip
On revient à présent sur l'application
$$
\hspace{24pt}
\Delta^{*}: \mathrm{H}^{*}\mathfrak{S}_{2}X\to\mathrm{H}^{*}(\mathrm{B}\mathfrak{S}_{2}\times X)=\mathbb{F}_{2}[u]\otimes\mathrm{H}^{*}X
\hspace{23pt}.
$$

\begin{pro}\label{cqse1} Soit $X$ un espace topologique avec $\mathrm{H}^{*}X$ de dimension finie en chaque degré~; on a les deux égalités suivantes~:

\medskip
{\em (a)}\hspace{8pt} $\mathop{\mathrm{im}}\Delta^{*}=\mathrm{R}_{1}\mathrm{H}^{*}X$,

\medskip
{\em (b)}\hspace{8pt} $\ker\Delta^{*}=\mathop{\mathrm{im}}\mathrm{tr}$.
\end{pro}

\bigskip 
\textit{Démonstration.} L'égalité $\Delta^{*}\mathrm{P}_{\hspace{-2pt}2}\hspace{1pt}x=\mathrm{St}_{1}\hspace{1pt}x$ et le fait que $\Delta^{*}$ est $\mathbb{F}_{2}[u]$-linéaire montrent que l'on a dans $\mathbb{F}_{2}[u]\otimes\mathrm{H}^{*}X$ l'inclusion $\mathrm{R}_{1}\mathrm{H}^{*}X\subset\mathop{\mathrm{im}}\Delta^{*}$.

\smallskip
Le scholie \ref{cqB} et le point (f) de \ref{R1} montrent que les $\mathbb{F}_{2}$-espaces vectoriels $\mathbb{N}$-gradués sous-jacents à $\mathrm{H}^{*}\mathfrak{S}_{2}X/\mathop{\mathrm{im}}\mathrm{tr}$ et $\mathrm{R}_{1}\mathrm{H}^{*}X$ sont tous deux isomorphes au $\mathbb{F}_{2}$-espace vectoriel $\mathbb{N}$-gradué $\mathbb{F}_{2}[u]\otimes\Phi\mathrm{H}^{*}X$. Cette observation implique l'égalité
$$
\hspace{24pt}
\dim\hspace{1pt}(\mathrm{H}^{*}\mathfrak{S}_{2}X/\mathop{\mathrm{im}}\mathrm{tr})^{n}
\hspace{4pt}=\hspace{4pt}
\dim\hspace{1pt}(\mathrm{R}_{1}\mathrm{H}^{*}X)^{n}
\hspace{23pt},
$$
pour tout $n$ dans $\mathbb{N}$, la notation $(-)^{n}$ désignant ici l'espace des éléments de degré $n$ d'un $\mathbb{F}_{2}$-espace vectoriel $\mathbb{N}$-gradué.

\smallskip
L'égalité $\Delta^{*}(\mathop{\mathrm{im}}\mathrm{tr})=0$ dit que $\Delta^{*}$ se factorise à travers $\mathrm{H}^{*}\mathfrak{S}_{2}X/\mathop{\mathrm{im}}\mathrm{tr}$, ce qui compte tenu de l'égalité ci-dessus implique l'inégalité
$$
\hspace{24pt}
\dim\hspace{1pt}
(\mathop{\mathrm{im}}\Delta^{*})^{n}
\hspace{4pt}\leq\hspace{4pt}
\dim\hspace{1pt}(\mathrm{R}_{1}\mathrm{H}^{*}X)^{n}
\hspace{23pt}.
$$
Les informations ci-dessus permettent de se convaincre que $\Delta^{*}$ induit un isomorphisme de  $\mathrm{H}^{*}\mathfrak{S}_{2}X/\mathop{\mathrm{im}}\mathrm{tr}$ sur $\mathrm{R}_{1}\mathrm{H}^{*}X$, d'où (a) et (b).
\hfill$\square$

\bigskip
\begin{rem}\label{ejP2} La suite exacte de Gysin d'un revêtement double (voir section \ref{Gysin}) montre que $\mathop{\mathrm{im}}\mathrm{tr}$ est un $\mathbb{F}_{2}[u]$-module {\em via} l'augmentation $\varepsilon: \mathbb{F}_{2}[u]\to~\hspace{-4pt}\mathbb{F}_{2}$.
\end{rem}

\bigskip
La proposition \ref{cqse1} peut se reformuler ainsi~:
                
\begin{scho}\label{cqse2} Soit $X$ un espace topologique avec $\mathrm{H}^{*}X$ de dimension finie en chaque degré~;  on a une suite exacte, naturelle en $X$, de $\mathrm{A}$- modules instables
$$
0\to\Lambda^{2}\mathrm{H}^{*}X\to\mathrm{H}^{*}\mathfrak{S}_{2}X\to\mathrm{R}_{1}\mathrm{H}^{*}X\to 0
\hspace{23pt},
$$
les deuxième et troisième  flèches étant respectivement induites par les homomorphismes $\mathrm{tr}:\mathrm{H}^{*}X\otimes\mathrm{H}^{*}X\to\mathrm{H}^{*}\mathfrak{S}_{2}X$ et $\Delta^{*}:\mathrm{H}^{*}\mathfrak{S}_{2}X\to\mathrm{H}^{*}\mathrm{B}\mathfrak{S}_{2}\otimes\mathrm{H}^{*}X$.

\medskip
\em{ (La notation $\Lambda^{2}$ est introduite en \ref{lambda2}.)}
\end{scho}

\medskip
Dans la suite exacte ci-dessus les deux termes de part et d'autre de $\mathrm{H}^{*}\mathfrak{S}_{2}X$ s'exprime fonctoriellement en $\mathrm{H}^{*}X$~; en va montrer qu'il en est de même pour $\mathrm{H}^{*}\mathfrak{S}_{2}X$ en précisant l'extension. On considère pour cela le diagramme commutatif suivant
$$
\begin{CD}
0@>>>\Lambda^{2}\mathrm{H}^{*}X@>>>\mathrm{H}^{*}\mathfrak{S}_{2}X
@>\Delta^{*}>>\mathrm{R}_{1}\mathrm{H}^{*}X
@>>>0 \\
& & @V\mathrm{id}VV  @V\pi^{*}VV \\
0@>>>\Lambda^{2}\mathrm{H}^{*}X@>>>(\mathrm{H}^{*}\otimes\mathrm{H}^{*}X)^{\mathfrak{S}_{2}}@>\nu>>\Phi\mathrm{H}^{*}X@>>>0
\end{CD}
$$
dans lequel~:

\smallskip
-- Les deux lignes sont exactes.

\smallskip
-- La deuxième flèche de la ligne du bas est l'inclusion $\ker(1-\sigma)\hookrightarrow\mathop{\mathrm{im}}(1+\sigma)$.

\smallskip
-- La flèche notée $\pi^{*}$ est celle induite par l'application $\mathrm{E}\mathfrak{S}_{2}\times (X\times X)\to\mathfrak{S}_{2}X$.

\smallskip
-- La flèche $\nu$ est la surjection $\mathrm{H}^{0}(\mathfrak{S}_{2};\mathrm{H}^{*}X\otimes \mathrm{H}^{*}X)\to
\widehat{\mathrm{H}}^{0}(\mathfrak{S}_{2};\mathrm{H}^{*}X\otimes \mathrm{H}^{*}X)$.

\medskip
Puisque les deux lignes sont exactes on peut compléter ce diagramme commutatif par une flèche, disons $f:\mathrm{R}_{1}\mathrm{H}^{*}X\to\Phi\mathrm{H}^{*}X$, uniquement déterminée. Soient $x$ un élément de $\mathrm{H}^{*}X$ et $j$ un entier naturel, les égalités
$$
\Delta^{*}\mathrm{P}_{\hspace{-2pt}2}\hspace{1pt}x=\mathrm{St}_{1}\hspace{1pt}x
\hspace{8pt},\hspace{8pt}
\pi^{*}(u^{j}\hspace{1pt}\mathrm{P}_{\hspace{-2pt}2}\hspace{1pt}x)=
\begin{cases}
x\otimes x & \text{pour $j=0$} \\ 0  & \text{pour $j>0$}
\end{cases}
\hspace{8pt}\text{et}\hspace{8pt}
\nu(x\otimes x)=\Phi x
$$
entraînent que l'on a
$$
f(u^{j}\hspace{1pt}\mathrm{St}_{1}\hspace{1pt}x)
\hspace{4pt}=\hspace{4pt}
\begin{cases}
\Phi x & \text{pour $j=0$}, \\ 0  & \text{pour $j>0$}.
\end{cases}
$$
Or on montre (voir par exemple \cite[4.2.6]{LZdest}) qu'il existe une unique transformation naturelle $\rho:\mathrm{R}_{1}\to\Phi$ de $\mathcal{U}$-endofoncteurs vérifiant
$$
\rho_{M}\hspace{1pt}(u^{j}\hspace{1pt}\mathrm{St}_{1}\hspace{1pt}x)
\hspace{4pt}=\hspace{4pt}
\begin{cases}
\Phi x & \text{pour $j=0$} \\ 0  & \text{pour $j>0$}
\end{cases}
$$
pour tout $\mathrm{A}$-module instable $M$ et tout élément $x$ de $M$~; on constate aussi que $\rho$ s'étend en une transformation naturelle, toujours notée  $\rho:\mathrm{R}_{1}\to\Phi$, de $\mathcal{K}$-endofoncteurs.

\medskip
On a donc obtenu le diagramme commutatif suivant~:
$$
\begin{CD}
\hspace{12pt}0@>>>\Lambda^{2}\mathrm{H}^{*}X
@>>>\mathrm{H}^{*}\mathfrak{S}_{2}X
@>\Delta^{*}>>\mathrm{R}_{1}\mathrm{H}^{*}X@>>>0\hspace{12pt} \\
& & @V\mathrm{id}VV  @V\pi^{*}VV  @V\rho VV  \\
\hspace{12pt}0@>>>\Lambda^{2}\mathrm{H}^{*}X@>>>(\mathrm{H}^{*}\otimes\mathrm{H}^{*}X)^{\mathfrak{S}_{2}}@>\nu>>\Phi\mathrm{H}^{*}X@>>>0\hspace{11pt}.
\end{CD}
$$
On observe que puisque la flèche verticale de gauche est l'identité, le carré de gauche est cartésien. Cette observation conduit aux définitions ci-après~:

\bigskip
\begin{defi}\label{cqalg} Soit $M$ un $\mathrm{A}$-module instable (resp.une $\mathrm{A}$-algèbre instable)~; on pose
$$
\hspace{24pt}
\mathfrak{S}_{2}M:=
\lim\hspace{2pt}(\hspace{2pt}(M\otimes M)^{\mathfrak{S}_{2}}\overset{\nu}{\longrightarrow}\Phi M\overset{\rho}{\longleftarrow}\mathrm{R}_{1}M\hspace{2pt})
\hspace{23pt},
$$
limite (produit fibré) dans la catégorie $\mathcal{U}$ (resp.  $\mathcal{K}$). Les deux endofoncteurs $\mathfrak{S}_{2}:\mathcal{U}\circlearrowleft$ et $\mathfrak{S}_{2}:\mathcal{K}\circlearrowleft$ ainsi définis commutent (en un sens évident) avec le foncteur oubli $\mathcal{K}\to\mathcal{U}$.
\end{defi}

\bigskip
Nous avons tout fait  pour avoir~:

\begin{pro}\label{cqfinal} Soit $X$ un espace topologique avec $\mathrm{H}^{*}X$ de dimension finie en chaque degré~;  on a un isomorphisme de $\mathrm{A}$-algèbres instables, et {\em a fortiori} de $\mathrm{A}$-modules instables, naturel en~$X$~:
$$
\hspace{24pt}
\mathrm{H}^{*}\mathfrak{S}_{2}X
\hspace{4pt}\cong\hspace{4pt}
\mathfrak{S}_{2}\mathrm{H}^{*}X
\hspace{23pt}.
$$                 
\end{pro}

Enoncé qui implique le suivant, concernant la cohomologie modulo $2$ des groupes finis~:

\begin{cor}\label{cqgroupe} Soit $G$ un groupe fini~;  on a un isomorphisme de $\mathrm{A}$-algèbres instables, et {\em a fortiori} de $\mathrm{A}$-modules instables, naturel en $G$~:
$$
\hspace{24pt}
\mathrm{H}^{*}(\mathfrak{S}_{2}\wr G)
\hspace{4pt}\cong\hspace{4pt}
\mathfrak{S}_{2}\mathrm{H}^{*}G
\hspace{23pt}.
$$                 
\end{cor}

\medskip
\textit{Démonstration.} Elle résulte des points suivants~:

\smallskip
--\hspace{8pt}$\mathrm{H}^{*}G$ est de dimension finie en chaque degré~;

\smallskip
--\hspace{8pt}$\mathrm{H}^{*}G$ est la cohomologie modulo $2$ de l'espace classifiant $\mathrm{B}G$~;

\smallskip
--\hspace{8pt} l'espace $\mathrm{B}(\mathfrak{S}_{2}\wr G)$ coïncide avec la construction quadratique $\mathfrak{S}_{2}\mathrm{B}G$.

\subsect{Retour à $\mathrm{H}^{*}\mathrm{S}_{n}$}\label{Sn}

Soit $n$ un entier naturel~; on rappelle que $\mathrm{S}_{n}$ désigne le $2$-Sylow du groupe symétrique $\mathfrak{S}_{n}$ décrit dans la sous-section \ref{2-Sylow-1} 
dont nous reprenons les notations.

\medskip
\begin{theo}\label{symnf} Le $\mathrm{A}$-module instable
$\mathrm{H}^{*}\mathrm{S}_{n}$ est $\mathcal{N}\hspace{-1,5pt}il$-fermé pour tout entier naturel $n$.
\end{theo}

\bigskip
\textit{Démonstration.} On a $\mathrm{S}_{n}=\mathrm{S}_{2^{m_{1}}}\times\mathrm{S}_{2^{m_{1}}}\times\ldots\times\mathrm{S}_{2^{m_{r}}}$ et donc
$$
\hspace{24pt}
\mathrm{H}^{*}\mathrm{S}_{n}
\hspace{4pt}\cong\hspace{4pt}
\mathrm{H}^{*}\mathrm{S}_{2^{m_{1}}}\otimes\mathrm{H}^{*}\mathrm{S}_{2^{m_{2}}}\otimes\ldots\otimes\mathrm{H}^{*}\mathrm{S}_{2^{m_{r}}}
\hspace{23pt}.
$$
Cet isomorphisme et la proposition \ref{tensnf} montrent qu'il suffit  de démontrer le théorème \ref{symnf} pour $n=2^{m}$ avec $m\in\mathbb{N}$. Comme l'on a $\mathrm{S}_{2^{m+1}}=\mathfrak{S}_{2}\wr\mathrm{S}_{2^{m}}$,  le corollaire \ref{cqgroupe} dit que l'on a un isomorphisme de  $\mathrm{A}$-modules instables  $\mathrm{H}^{*}\mathrm{S}_{2^{m+1}}=\mathfrak{S}_{2}\mathrm{H}^{*}\mathrm{S}_{2^{m}}$~; on achève par récurrence  sur $m$ à l'aide de la proposition \ref{quadnilf} ci-après ($\mathrm{H}^{*}\mathrm{S}_{2^{0}}=\mathbb{F}_{2}$ est  $\mathcal{N}\hspace{-1,5pt}il$-fermé !).
\hfill$\square$

\medskip
\begin{cor}\label{quillensym}
L'application de Quillen $\mathrm{q}_{\mathfrak{S}_{n}}:\mathrm{H}^{*}(\mathfrak{S}_{n}\;\mathbb{F}_{2})\to\mathrm{L}(\mathfrak{S}_{n})$ est un isomorphisme pour tout entier naturel $n$.
\end{cor}

\medskip
\textit{Démonstration.} C'est une illustration du corollaire \ref{Sylow-2}.
\hfill$\square$

\bigskip
\begin{pro}\label{quadnilf} Soit $M$ un $\mathrm{A}$-module instable~; si $M$ est  $\mathcal{N}\hspace{-1.5pt}il$-fermé  alors il en est de même pour $\mathfrak{S}_{2}M$.
\end{pro}

\medskip
\textit{Démonstration.} Par définition même de l'endofoncteur  $\mathfrak{S}_{2}:\mathcal{U}\circlearrowleft$, on a la suite exacte suivante dans la catégorie $\mathcal{U}$:
$$
\hspace{24pt}
0\to\mathfrak{S}_{2}M\to(M\otimes M)^{\mathfrak{S}_{2}}\oplus\mathrm{R}_{1}M\to\Phi M
\hspace{23pt}.
$$
On invoque le scholie \ref{critnilf-1}.~:

\smallskip
-- Comme $M$ est  $\mathcal{N}\hspace{-1.5pt}il$-fermé, il est \textit{a fortiori} réduit, c'est donc aussi le cas pour~$\Phi M$.

\smallskip
-- Le  $\mathrm{A}$-module instable $(M\otimes M)^{\mathfrak{S}_{2}}$ est $ \mathcal{N}\hspace{-1.5pt}il$-fermé d'après \ref{tensnf} et \ref{invnf}.

\smallskip
 -- Le  $\mathrm{A}$-module instable $\mathrm{R}_{1}M$ est aussi  $\mathcal{N}\hspace{-1.5pt}il$-fermé, voir ci-dessous.
  \hfill$\square$
  
  \medskip
\begin{pro}\label{R1nilf} Soit $M$ un $\mathrm{A}$-module instable~; si $M$ est  $\mathcal{N}\hspace{-1.5pt}il$-fermé  alors il en est de même pour $\mathrm{R}_{1}M$.
\end{pro}

\medskip
\textit{Démonstration.} La condition $(iii)$ de \ref{critnilf-3} dit qu'il existe une $\mathcal{U}$-suite exacte de la forme $0\to M\to I^{0}\to I^{1}$ avec $I^{0}$ et $I^{1}$ somme directe de cohomologies modulo $2$ de $2$-groupes abéliens élémentaires. Puisque l'endofonteur $\mathrm{R}_{1}$ est exact on a encore une $\mathcal{U}$-suite exacte  $0\to \mathrm{R}_{1}M\to \mathrm{R}_{1}I^{0}\to \mathrm{R}_{1}I^{1}$~;  compte tenu de \ref{critnilf-1} et du fait que $ \mathrm{R}_{1}$ préserve les somme directes, on est ramené à montrer que $\mathrm{R}_{1}\mathrm{H}^{*}V$ est $ \mathcal{N}\hspace{-1.5pt}il$-fermé pour  tout $2$-groupe abélien élémentaire $V$. On montre dans \cite{LZdest} (Lemme 4.4.6.2) que $\mathrm{R}_{1}\mathrm{H}^{*}V$ est isomorphe au sous-$\mathrm{A}$-module instable d'invariants $(\mathrm{H}^{*}(\mathbb{Z}/2\oplus V))^{\Gamma}$,  $\Gamma$ désignant le sous-groupe de $\mathrm{Aut}(\mathbb{Z}/2\oplus V)$ constitué des automorphismes qui sont l'identité sur $V$ ($\Gamma$ agit à droite sur $\mathrm{H}^{*}(\mathbb{Z}/2\oplus V)$). On achève en invoquant à nouveau  \ref{invnf}.
\hfill$\square$

\sect{Variations sur la suite exacte de Gysin d'un revêtement double}\label{Gysin}

\bigskip
Soit $p:Y\to X$ un revêtement double~; soit $e$  sa classe caractéristique ($e$~est un élément de $\mathrm{H}^{1}(X;\mathbb{F}_{2})=\mathrm{H}^{1}(X;\mathbb{Z}/2)$, $p$ est un revêtement galoisien de groupe de Galois $\mathbb{Z}/2$).

\medskip
Comme précédemment, on abrège ci-après  la notation $\mathrm{H}^{*}(-;\mathbb{F}_{2})$ en $\mathrm{H}^{*}(-)$ ou simplement $\mathrm{H}^{*}-$.

\bigskip
On considère la longue suite exacte de Gysin
$$
\begin{CD}
\hspace{6pt}\ldots\to\mathrm{H}^{n-1}X@>e\smile>>\mathrm{H}^{n}X@>p^{*}>>\mathrm{H}^{n}Y@>\mathrm{tr}>>\mathrm{H}^{n}X@>e\smile>>\mathrm{H}^{n+1}X\to\ldots\hspace{4pt}.
\end{CD}
$$
On note respectivement  $\mathop{\mathrm{im}} e$ et $\ker e$  l'image et le noyau, dans $\mathrm{H}^{*}X$, de la multiplication par $e$.

\bigskip
\begin{pro}\label{A-stable} Dans $\mathrm{H}^{*}X$,  $\mathop{\mathrm{im}} e$ et $\ker e$ sont stables sous l'action de l'algèbre de Steenrod.
\end{pro}

\bigskip
\textit{Démonstration.} Dans le cas de $\mathop{\mathrm{im}} e$ ceci résulte immédiatement de la formule
$$
\hspace{24pt}
\mathrm{Sq}^{i}(e\hspace{1pt}x)
\hspace{4pt}=\hspace{4pt}
e\hspace{1pt}\mathrm{Sq}^{i}x+e^{2}\hspace{1pt}\mathrm{Sq}^{i-1}x
\hspace{23pt}.
$$
Cette même formule permet de régler le cas de $\ker e$ par récurrence sur $i$.
\hfill$\square$

\bigskip
On pose $\mathop{\mathrm{coker}} e:=\mathrm{H}^{*}X\slash\mathop{\mathrm{im}} e$~; d'après ce qui précède $\mathop{\mathrm{coker}} e$ est muni d'une structure canonique de $\mathrm{A}$-module instable.

\bigskip
\begin{scho}\label{A-Gysin} La suite exacte courte de Gysin
$$
\begin{CD}
0@>>>\mathop{\mathrm{coker}} e@>p^{*}>>\mathrm{H}^{*}Y@>\mathrm{tr}>>
\ker e @>>> 0
\end{CD}
$$
est une suite exacte dans la catégorie $\mathcal{U}$.
\end{scho}

\medskip
\begin{rem}\label{sdP} Si $\mathrm{H}^{*}X$ est de dimension finie en chaque degré, alors il en est de même pour $\mathrm{H}^{*}Y$ et les séries de Poincaré de $\mathrm{H}^{*}X$ et $\ker e$ déterminent celle de $\mathrm{H}^{*}Y$. Précisons~: soit $\mathrm{S}(E;t)\in\mathbb{N}[[t]]$ la série de Poincaré d'un  $\mathbb{F}_{2}$-espace vectoriel $\mathbb{N}$-gradué $E$ de dimension finie en chaque degré; on a~:
$$
\hspace{24pt}
\mathrm{S}(\mathrm{H}^{*}Y;t)
\hspace{4pt}=\hspace{4pt}
(1-t)\hspace{1pt}\mathrm{S}(\mathrm{H}^{*}X;t)+(1+t)\hspace{1pt}\mathrm{S}(\ker e\hspace{1pt};t)
\hspace{23pt}.
$$
\end{rem}

\begin{pro}\label{kere} Si $\mathrm{H}^{*}X$ est $\mathcal{N}\hspace{-1.5pt}il$-fermé alors il en est de même pour $\ker e$.
\end{pro}

\textit{Démonstration.} Compte tenu de \ref{critnilf-0}, il faut montrer que la $\mathrm{A}$-algèbre instable quotient $\mathrm{H}^{*}X\slash\ker e$ est réduite, c'est-à-dire que si $x$ est un élément de $\mathrm{H}^{*}X$ avec $e\hspace{1pt}x^{2}=0$ alors on a aussi $e\hspace{1pt}x=0$. Or on a les implications $e\hspace{1pt}x^{2}=0\Rightarrow(e\hspace{1pt}x)^{2}=0\Rightarrow e\hspace{1pt}x=0$~; la première est triviale, la seconde tient à l'hypothèse faite sur $\mathrm{H}^{*}X$~: un $\mathrm{A}$-module instable $\mathcal{N}\hspace{-1.5pt}il$-fermé est \textit{a fortiori} réduit.
\hfill$\square$

\bigskip
\begin{cor}\label{cokere} Si $\mathrm{H}^{*}X$ est $\mathcal{N}\hspace{-1.5pt}il$-fermé les deux propriétés suivantes sont équi\-valentes~:
\begin{itemize}
\item[(i)] $\mathrm{H}^{*}Y$ est $\mathcal{N}\hspace{-1.5pt}il$-fermé~;
\item[(ii)] $\mathop{\mathrm{coker}} e$ est $\mathcal{N}\hspace{-1.5pt}il$-fermé.
\end{itemize}
\end{cor}

\textit{Démonstration.} Compte tenu du scholie \ref{A-Gysin} et de la proposition ci-dessus, on peut invoquer la proposition \ref{critnilf-2}.
\hfill$\square$

\sect{Sur les $\mathrm{P}$-$\mathrm{A}$-modules et $\mathrm{P}$-$\mathrm{A}$-algèbres instables}\label{P-A-inst}

\subsect{Les notions de $\mathrm{P}$-$\mathrm{A}$-modules et $\mathrm{P}$-$\mathrm{A}$-algèbres instables}\label{P-A-notions}

\medskip
Rappelons que $\mathrm{P}$ désigne la cohomologie modulo $2$ du groupe $\mathbb{Z}/2$~; $\mathrm{P}$ est une $\mathrm{A}$-algèbre instable.

\bigskip
\begin{defi}\label{def-P-A-1} Un {\em $\mathrm{P}$-$\mathrm{A}$-module instable} est un $\mathrm{A}$-module instable $M$ muni d'une application $\mathrm{A}$-linéaire $\mathrm{P}\otimes M\to M$ qui fait de $M$ un $\mathrm{P}$-module ($\mathbb{N}$-gradué). Un homomorphisme  de $\mathrm{P}$-$\mathrm{A}$-modules instables est un homomorphisme de $\mathbb{F}_{2}$-espaces vectoriels $\mathbb{N}$-gradués  à la fois $\mathrm{A}$-linéaire et  $\mathrm{P}$-linéaire. La catégorie des  $\mathrm{P}$-$\mathrm{A}$-modules instables est notée $\mathrm{P}$-$\mathcal{U}$. 
\end{defi}

\bigskip
\begin{exple}\label{P-tens-M} Produit tensoriel de $\mathrm{P}$ et d'un $\mathrm{A}$-module instable.

\smallskip
Soit $M$ un $\mathrm{A}$-module instable ; le $\mathcal{U}$-morphisme
$$
\begin{CD}
\mathrm{P}\otimes(\mathrm{P}\otimes M)=(\mathrm{P}\otimes\mathrm{P})\otimes M@>\varphi\hspace{1pt}\otimes\hspace{1pt}1>>\mathrm{P}\otimes M
\end{CD}
$$
($\varphi$ désignant le $\mathcal{U}$-morphisme donné par le produit de $\mathrm{P}$, voir note \footref{axiome-K1}) fait de $\mathrm{P}\otimes M$ un $\mathrm{P}$-$\mathrm{A}$-module instable. Le foncteur $\mathcal{U}\to\mathrm{P}\text{-}\mathcal{U}$ ainsi défini est adjoint à gauche du foncteur oubli $\mathrm{P}\text{-}\mathcal{U}\to\mathcal{U}$ :
$$
\mathrm{Hom}_{\mathcal{U}}(M,N)
\hspace{4pt}=\hspace{4pt}
\mathrm{Hom}_{\mathcal{U}}(\mathrm{P}\otimes M,N)
$$
pour tout $\mathrm{A}$-module instable $M$ et tout $\mathrm{P}$-$\mathrm{A}$-module instable $N$
\end{exple}

\bigskip
\begin{defi}\label{def-P-A-2} Une {\em $\mathrm{P}$-$\mathrm{A}$-algèbre instable} est une $\mathrm{A}$-algèbre instable~$K$ munie d'un $\mathcal{K}$-morphisme $f:\mathrm{P}\to K$. Soient $(K_{0}\hspace{1pt};f_{0})$ et $(K_{1}\hspace{1pt};f_{1})$ deux $\mathrm{P}$-$\mathrm{A}$-algèbres instables~; un homomorphisme  de $\mathrm{P}$-$\mathrm{A}$-algèbres instables de $(K_{0}\hspace{1pt};f_{0})$ dans $(K_{1}\hspace{1pt};f_{1})$ est la donnée d'un $\mathcal{K}$-morphisme $g:K_{0}\to K_{1}$ avec $f_{1}=g\circ f_{0}$. La catégorie des  $\mathrm{P}$-$\mathrm{A}$-modules instables est notée $\mathrm{P}\backslash\mathcal{K}$. 
\end{defi}

\bigskip
On dispose d'un ``foncteur oubli'' évident $\mathrm{P}\backslash\mathcal{K}\to\mathrm{P}\text{-}\mathcal{U}$~: la donnée d'un $\mathcal{K}$-morphisme $\mathrm{P}\to K$ fait du $\mathrm{A}$-module instable sous-jacent à $K$ un $\mathrm{P}$-$\mathrm{A}$-module instable.

\bigskip
Plus généralement, on définit, \textit{mutatis mutandis}, les catégories $L\text{-}\mathcal{U}$ et $L\backslash\mathcal{K}$ pour toute $\mathrm{A}$-algèbre instable $L$.

\bigskip
On précise ci-après la structure de $\mathrm{P}$ et on démystifie la catégorie $\mathrm{P}\backslash\mathcal{K}$.

\medskip
 La $\mathbb{F}_{2}$-algèbre $\mathbb{N}$-graduée sous-jacente à $\mathrm{P}$ est isomorphe à l'algèbre de polynômes $\mathbb{F}_{2}[u]$, $u$ désignant une indéterminée de degré $1$. La~structure de $\mathrm{A}$-algèbre instable de $\mathrm{P}$ est uniquement déterminée par cette information. Rappelons  pourquoi. Soient $K$ une $\mathrm{A}$-algèbre instable et $x$ un élément de~$K^{1}$, les axiomes que vérifie une $\mathrm{A}$-algèbre instable impliquent $\mathrm{Sq}^{i}\hspace{1pt}x^{n}=\binom{n}{i}\hspace{1pt}x^{n+i}$ pour tous entiers naturels $i$ et $n$ (formule due à Henri Cartan).

\medskip
La $\mathrm{A}$-algèbre instable $\mathrm{P}$ est en outre un co-groupe dans la catégorie $\mathcal{K}$. En~clair, en plus des morphismes $\varphi:\mathrm{P}\otimes\mathrm{P}\to\mathrm{P}$ (produit) et $\eta:\mathbb{F}_{2}\to\mathrm{P}$ (unité), on dispose de morphismes $\psi: \mathrm{P}\to\mathrm{P}\otimes\mathrm{P}$ (diagonale) et $\varepsilon:\mathrm{P}\to\mathbb{F}_{2}$ (augmentation) qui font de la $\mathbb{F}_{2}$-algèbre $\mathbb{N}$-graduée sous-jacente à $\mathrm{P}$ une algèbre de Hopf bicommutative~;  la diagonale $\psi$ est le $\mathcal{K}$-morphisme induit par l'homomorphisme de groupes $\mathbb{Z}/2\times\mathbb{Z}/2\to\mathbb{Z}/2, (x_{1},x_{2})\mapsto x_{1}+x_{2}$.

\medskip
 La formule de Cartan montre que l'application
$$
\hspace{24pt}
\mathrm{Hom}_{\mathcal{K}}(\mathrm{P},K)\to K^{1}
\hspace{12pt},\hspace{12pt}
f\mapsto f(u)
\hspace{24pt}
$$
est bijective. D'après ce qui précède le foncteur $K\mapsto\mathrm{Hom}_{\mathcal{K}}(\mathrm{P},K)$ peut être vu comme un foncteur à valeurs dans la catégorie des groupes abéliens et la bijection ci-dessus est un isomorphisme de groupes (naturel en $K$).

\medskip
La discussion ci-dessus  montre qu'un objet de $\mathrm{P}\backslash\mathcal{K}$ peut être vu comme un couple $(K;e)$, $K$ étant une $\mathrm{A}$-algèbre instable et~$e$ un élément de $K^{1}$,  et qu'un morphisme de $(K;e)$ dans $(K';e')$ est un $\mathcal{K}$-morphisme $g:K\to K'$ avec $g(e)=e'$.

\bigskip
On revient maintenant sur la donnée initiale de la section \ref{Gysin}~: 

\medskip
Soit $p:Y\to X$ un revêtement double~; soit $e\in\mathrm{H}^{1}X$  sa classe caractéristique. Compte tenu de ce qui précède, la donnée de $e\in\mathrm{H}^{1}(X;\mathbb{F}_{2})$ est  équivalente à celle d'un homomorphisme de $\mathrm{A}$-algèbres instables $\mathrm{P}\to\mathrm{H}^{*}X$~; on observera que si le revêtement $p$ est classifié par une application $c:X\to\mathrm{B}(\mathbb{Z}/2)$ alors cet homomorphisme s'identifie à $c^{*}:\mathrm{H}^{*}\mathrm{B}(\mathbb{Z}/2)\to\mathrm{H}^{*}X$. Reformulons~:  La donnée de $e$ fait de $\mathrm{H}^{*}X$ une $\mathrm{P}$-$\mathrm{A}$-algèbre instable et du même coup un $\mathrm{P}$-$\mathrm{A}$-module  instable.

\bigskip
On dispose d'un foncteur ``oubli'' évident $\mathrm{P}\text{-}\mathcal{U}\to\mathcal{U}$. Ce foncteur possède une ``section'' tout aussi évidente~:

\bigskip
\begin{defi}\label{theta} Soit $M$ un $A$-module instable, l'augmentation $\mathrm{P}\to\mathbb{F}_{2}$ permet de faire de $M$ un $\mathrm{P}$-$\mathrm{A}$-module instable, un tel  $\mathrm{P}$-$\mathrm{A}$-module instable, sera dit {\em$\mathrm{P}$-trivial}. On note $\theta:\mathcal{U}\to\mathrm{P}\text{-}\mathcal{U}$ le foncteur ainsi défini.
\end{defi}

\medskip
Le foncteur $\theta$ possède un adjoint à droite défini ci-dessous.

\bigskip
\begin{pro-def}\label{tau} Soit $M$ un $\mathrm{P}$-$\mathrm{A}$-module instable. On pose 
$$
\hspace{24pt}
\tau M
\hspace{4pt}:=\hspace{4pt}
\{\hspace{2pt}x\hspace{2pt};\hspace{2pt} x\in M \hspace{4pt}\text{et}\hspace{4pt}u\hspace{1pt}x=0\hspace{2pt}\}
\hspace{23pt};
$$
$\tau M$ est stable sous l'action de $\mathrm{A}$, c'est le plus grand sous-$\mathrm{P}$-$\mathrm{A}$-module instable $\mathrm{P}$-trivial de $M$, on l'appelle la {\em partie triviale} de $M$. Le foncteur $M\mapsto\tau M$ sera suivant le contexte considéré comme un endofoncteue de $\mathrm{P}\text{-}\mathcal{U}$ ou un foncteur $\mathrm{P}\text{-}\mathcal{U}\to\mathcal{U}$ (adjoint à droite de $\theta$).
\end{pro-def}

\medskip
\textit{Démonstration.} La démonstration de la stabilité de $\tau M$ sous l'action de $\mathrm{A}$ est la même que celle du $\ker e$ de la proposition \ref{A-stable}.
\hfill$\square$

\subsect{Produit tensoriel de $\mathrm{P}$-$\mathrm{A}$-modules instables}\label{P-A-prod-1}

Soient $M_{1}$ et $M_{2}$ deux $\mathrm{P}$-$\mathrm{A}$-modules instables. Le produit tensoriel $M_{1}\otimes M_{2}$ des deux $\mathrm{A}$-modules instables sous-jacents est naturellement un $L$-$\mathrm{A}$-module instable avec  $L=\mathrm{P}\otimes\mathrm{P}$ et tout aussi naturellement un $\mathrm{P}$-$\mathrm{A}$-module instable grâce au $\mathcal{K}$-morphisme $\psi:\mathrm{P}\to\mathrm{P}\otimes\mathrm{P}$ introduit en \ref{P-A-notions}. 
Le $\mathrm{P}$-$\mathrm{A}$-module instable ainsi défini est appelé le {\em produit tensoriel de $M_{1}$ et $M_{2}$}.

\medskip
Comme la loi de composition interne d'un groupe est associative le produit tensoriel de $\mathrm{P}$-$\mathrm{A}$-modules instables est associatif (en un sens que le lecteur devinera aisément). Soit $(M_{1},M_{2},\ldots M_{r})$ une suite finie de $\mathrm{P}$-$\mathrm{A}$-modules instables~; le $\mathrm{P}$-$\mathrm{A}$-module instable $M_{1}\otimes M_{2}\otimes\ldots\otimes M_{r}$ peut être défini \textit{via}  le $\mathcal{K}$-morphisme $\psi_{r}:\mathrm{P}\to\mathrm{P}\otimes\mathrm{P}\otimes\ldots\otimes\mathrm{P}$ induit par l'homomorphisme de groupes
$$
\hspace{8pt}
\mathbb{Z}/2\times\mathbb{Z}/2\times\ldots\times\mathbb{Z}/2\to\mathbb{Z}/2
\hspace{6pt},\hspace{6pt}(x_{1},x_{2},\ldots,x_{r})\mapsto x_{1}+x_{2}+\ldots+x_{r}
\hspace{7pt}.
$$

\subsect{Produit tensoriel sur $\mathrm{P}$ de deux $\mathrm{P}$-$\mathrm{A}$-modules instables}\label{P-A-prod-2}

\medskip
On note $\mathcal{E}$ la catégorie des  $\mathbb{F}_{2}$-espaces vectoriels $\mathbb{N}$-gradués et $\mathrm{P}\text{-}\mathcal{E}$ la catégorie des $\mathbb{F}_{2}$-espaces vectoriels $\mathbb{N}$-gradués munis d'une structure de $\mathrm{P}$-module (au sens gradué). Il est clair que l'on dispose d'un foncteur oubli $\mathrm{P}\text{-}\mathcal{U}\to\mathrm{P}\text{-}\mathcal{E}$.

\medskip
En vrac quelques observations et définitions~:

\medskip
(1) Puisque l'algèbre $\mathrm{P}$ est commutative, les notions de $\mathrm{P}$-module à gauche et à droite coïncident.

\medskip
(2) Soient  $M_{1}$ et $M_{2}$ deux objets de $\mathrm{P}$-$\mathcal{E}$~; par définition $M_{1}\otimes_{\mathrm{P}}M_{2}$ est le $\mathcal{E}$-objet coégalisateur des  deux applications $M_{1}\otimes\mathrm{P}\otimes M_{2}\to M_{1}\otimes M_{2}$ suivantes~:
$$
\hspace{24pt}
x_{1}\otimes a\otimes x_{2}\mapsto ax_{1}\otimes x_{2}
\hspace{12pt}\text{et}\hspace{12pt}
x_{1}\otimes a\otimes x_{2}\mapsto x_{1}\otimes ax_{2}
\hspace{23pt}.
$$

\medskip
(3) Soient  $M_{1}$ et $M_{2}$ deux objets de  $\mathrm{P}$-$\mathcal{U}$~; puisque les deux applications ci-dessus sont $\mathrm{A}$-linéaires,  la $\mathcal{E}$-structure de  $M_{1}\otimes_{\mathrm{P}}M_{2}$ peut être naturellement enrichie en une $\mathcal{U}$-structure.

\medskip
(4) A nouveau puisque $\mathrm{P}$ est commutative, le foncteur $-\otimes_{\mathrm{P}}-:\mathrm{P}\text{-}\mathcal{E}\times\mathrm{P}\text{-}\mathcal{E}\to\mathcal{E}$ (resp. $-\otimes_{\mathrm{P}}-:\mathrm{P}\text{-}\mathcal{U}\times\mathrm{P}\text{-}\mathcal{U}\to\mathcal{U}$) peut être ``relevé'' en un foncteur à valeurs dans $\mathrm{P}\text{-}\mathcal{E}$ (resp. $\mathrm{P}\text{-}\mathcal{U}$).

\medskip
(5) Le diagramme de foncteurs
$$
\begin{CD}
\mathrm{P}\text{-}\mathcal{U}\times\mathrm{P}\text{-}\mathcal{U}
@>-\otimes_{\mathrm{P}}-:>>
\mathrm{P}\text{-}\mathcal{U} \\
@V\mathrm{oubli}\times\mathrm{oubli}VV @V\mathrm{oubli}VV \\
\mathrm{P}\text{-}\mathcal{E}\times\mathrm{P}\text{-}\mathcal{E}
@>-\otimes_{\mathrm{P}}-:>>
\mathrm{P}\text{-}\mathcal{E}
\end{CD}
$$
est commutatif.

\subsect{Foncteurs $\mathrm{Tor}_{k}^{\mathrm{P}}(-,-)$ dans la catégorie $\mathrm{P}\text{-}\mathcal{U}\times\mathrm{P}\text{-}\mathcal{U}$}\label{Tor-P}

Quelques observations et définitions (toujours en vrac)~:

\bigskip
(1) Les foncteurs $\mathrm{P}\text{-}\mathcal{U}\to\mathrm{P}\text{-}\mathcal{U}, M_{2}\mapsto M_{1}\otimes_{\mathrm{P}}M_{2}\hspace{4pt}\text{et}\hspace{4pt}M_{1}\mapsto M_{1}\otimes_{\mathrm{P}}M_{2}$ sont exacts à droite.

\medskip
(2) La catégorie $\mathrm{P}\text{-}\mathcal{U}$ a assez de projectifs.

\smallskip
On résume  ci-après  la théorie des projectifs de $\mathrm{P}\text{-}\mathcal{U}$ (théorie formelle, on peut y remplacer $\mathrm{P}$ par une $\mathrm{A}$-algèbre instable $L$, avec $L^{0}=\mathbb{F}_{2}$, arbitraire).

\smallskip
Soit $n$ un entier naturel. On note $\mathrm{F}(n)$ le $\mathrm{A}$-module instable librement engendré par un générateur de degré $n$~;  $\mathrm{F}(n)$ représente le foncteur, défini sur $\mathcal{U}$ et à valeurs dans la catégorie des $\mathbb{F}_{2}$-espaces vectoriels, $M\mapsto M^{n}$~: $\mathrm{Hom}_{\mathcal{U}}(\mathrm{F}(n),M)=M^{n}$.

\smallskip
Cette égalité implique la suivante $\mathrm{Hom}_{\mathcal{U}}(\mathrm{P}\otimes\mathrm{F}(n),M)=M^{n}$ (voir \ref{P-tens-M}) qui montre que $\mathrm{P}\otimes\mathrm{F}(n)$ est un projectif de $\mathrm{P}\text{-}\mathcal{U}$. Soient $M$ un  objet de  $\mathrm{P}\text{-}\mathcal{U}$ et $B$ une base du $\mathrm{P}$-module gradué  sous-jacent à $M$, le $\mathrm{P}\text{-}\mathcal{U}$-homomorphisme $\bigoplus_{x\in B}\mathrm{P}\otimes\mathrm{F}(\vert x\vert)\to M$ est par construction un épimorphisme~:  $\mathrm{P}\text{-}\mathcal{U}$ a assez de projectifs.
Du coup tout projectif de $\mathrm{P}\text{-}\mathcal{U}$ est facteur direct d'une telle somme directe de $\mathrm{P}\otimes\mathrm{F}(n)$'s ce qui entraîne en particulier qu'il est libre comme $\mathrm{P}$-module (voir \ref{P-lib-2}). On dispose en fait d'un résultat plus précis~: 

\medskip
(3) Tout projectif de $\mathrm{P}\text{-}\mathcal{U}$ est somme directe de $\mathrm{P}\otimes\mathrm{F}(n)$'s.\footnote{Soit $\widetilde{\mathrm{A}}$ l'idéal d'augmentation de $\mathrm{A}$. Soient $M$ un projectif de $\mathrm{P}\text{-}\mathcal{U}$, $\bar{\bar{B}}$ une base (homogène) du $\mathbb{F}_{2}$-espace vectoriel $\mathbb{N}$-gradué $(M/uM)/\widetilde{\mathrm{A}}\hspace{1pt}(M/uM)$ et $B$ un ``relèvement'' de $\bar{\bar{B}}$ dans $M$~; le $\mathrm{P}\text{-}\mathcal{U}$-homomorphisme $\bigoplus_{x\in B}\mathrm{P}\otimes\mathrm{F}(\vert x\vert)\to M$ est un isomorphisme.}

\medskip
(4)  On définit les  foncteurs  $\mathrm{Tor}^{\mathrm{P}}_{k}(-,-)$ comme les dérivés à droite des foncteurs produit tensoriel.

\medskip
(5) Dériver par rapport au premier ou second argument ``donne le même résultat''.

\medskip
(6) Les foncteurs $\mathrm{Tor}^{\mathrm{P}}_{k}(-,-)$ ``commutent'' avec le foncteur oubli $\mathrm{P}\text{-}\mathcal{U}\to\mathrm{P}\text{-}\mathcal{E}$ (puisque tout projectif de $\mathrm{P}\text{-}\mathcal{U}$ est libre comme $\mathrm{P}\text{-}\mathbb{F}_{2}$-espace vectoriel $\mathbb{N}$-gradué).

\medskip
(7) On a $\mathrm{Tor}^{\mathrm{P}}_{k}(M_{1},M_{2})=0$ pour $k>1$ (conséquence de l'observation (6) et du scholie \ref{projdim}).

\medskip
Voici une autre conséquence de l'observation (6)~:

\begin{pro}\label{res-P-lib} Soient $M_{1}$ et $M_{2}$ deux $\mathrm{P}$-$\mathrm{A}$-modules instables.

\smallskip
Si $M_{1}\leftarrow C_{\bullet}$ est une résolution dans la catégorie $\mathrm{P}\text{-}\mathcal{U}$ qui est libre  comme résolution dans la catégorie $\mathrm{P}\text{-}\mathcal{E}$ alors on a un $\mathrm{P}\text{-}\mathcal{U}$-isomorphisme
$$
\hspace{24pt}
\mathrm{Tor}^{\mathrm{P}}_{k}(M_{1},M_{2})
\hspace{4pt}\cong\hspace{4pt}
\mathrm{H}_{k}(C_{\bullet}\otimes_{\mathrm{P}}M_{2})
\hspace{24pt}
$$
\end{pro}

\textit{Démonstration.} Un $\mathrm{P}$-$\mathrm{A}$-module instable qui est libre comme $\mathrm{P}$-module est $-\otimes_{\mathrm{P}}M_{2}$-acyclique (en fait le seul cas non-trivial est $k=1$).
\hfill$\square$

 \bigskip
 Exemple d'application de la proposition ci-dessus~:

\begin{pro}\label{sigmatau} Soit $M$ un $\mathrm{P}$-$\mathrm{A}$-module instable~; on a un isomorphisme canonique de $\mathrm{P}$-$\mathrm{A}$-modules instables
$$
\hspace{24pt}
\mathrm{Tor}_{1}^{\mathrm{P}}\hspace{1pt}(\mathbb{F}_{2},M)
\hspace{4pt}\cong\hspace{4pt}
\Sigma\hspace{2pt}\tau M
\hspace{23pt}.
$$
\end{pro}

\textit{Démonstration.} On considère l'augmentation $\varepsilon: \mathrm{P}\to\mathbb{F}_{2}$~; on pose $\widetilde{\mathrm{P}}:=\ker\varepsilon$ et on note $\mathrm{i}$ l'inclusion de $\widetilde{\mathrm{P}}$ dans $\mathrm{P}$. La suite exacte $0\to\widetilde{\mathrm{P}}\overset{\mathrm{i}}{\to}\mathrm{P}\to\mathbb{F}_{2}\to 0$ fournit une résolution de $\mathbb{F}_{2}$ dans la catégorie $\mathrm{P}\text{-}\mathcal{U}$ qui est une résolution libre dans la catégorie $\mathrm{P}\text{-}\mathcal{E}$. Comme $\widetilde{\mathrm{P}}$ est un $\mathrm{P}$-module libre de base $\{u\}$,
l'application $M\to\widetilde{\mathrm{P}}\otimes_{\mathrm{P}}M, x\mapsto u\otimes_{\mathrm{P}}x$ est un isomorphisme, de degré $1$,  de $\mathbb{F}_{2}$-espaces vectoriels $\mathbb{N}$-gradués. L'application $\mathrm{i}\otimes_{\mathrm{P}}M$ s'identifie à l'application $\widetilde{\mathrm{P}}\otimes_{\mathrm{P}}M\to\mathrm{P}\otimes_{\mathrm{P}}M=M\hspace{1pt},\hspace{1pt}u\otimes_{\mathrm{P}}x\mapsto u\hspace{1pt}x$ et son noyau au sous-$\mathrm{P}$-$\mathrm{A}$-module instable $\widetilde{\mathrm{P}}\otimes_{\mathrm{P}}\tau M$ de $\widetilde{\mathrm{P}}\otimes_{\mathrm{P}}M$. On a un isomorphisme canonique de $\mathrm{P}$-$\mathrm{A}$-modules instables $\mathrm{Tor}_{1}^{\mathrm{P}}\hspace{1pt}(\mathbb{F}_{2},M)\cong\widetilde{\mathrm{P}}\otimes_{\mathrm{P}}\tau M$. Il reste à expliciter la structure de $\mathrm{A}$-module instable de $\widetilde{\mathrm{P}}\otimes_{\mathrm{P}}\tau M$.

\smallskip
Soient $x$ un élément de $\tau M$ et $j$ un entier naturel~; on a~:
$$
\hspace{24pt}
\mathrm{Sq}^{j}\hspace{1pt}(u\otimes_{\mathrm{P}}x)
\hspace{4pt}=\hspace{4pt}
\sum_{k=0}^{j}
\hspace{4pt}
\mathrm{Sq}^{k}u\otimes_{\mathrm{P}}Sq^{j-k}x
\hspace{23pt}.
$$
Pour $k\geq 1$ le terme $\mathrm{Sq}^{k}u\otimes_{\mathrm{P}}Sq^{j-k}x$ est de la forme $uv_{k}\otimes_{\mathrm{P}}Sq^{j-k}x$ avec $v_{k}\in\widetilde{\mathrm{P}}$ soit encore $u\otimes_{\mathrm{P}}v_{k}Sq^{j-k}x$~; puisque $\tau M$ est stable sous l'action de~$\mathrm{A}$ ce terme est nul. On a donc $\mathrm{Sq}^{j}\hspace{1pt}(u\otimes_{\mathrm{P}}x)=u\otimes_{\mathrm{P}}\mathrm{Sq}^{j}x$ formule qui montre que le $\mathrm{A}$-module instable $\widetilde{\mathrm{P}}\otimes_{\mathrm{P}}\tau M$ s'identifie à la suspension $\Sigma\hspace{2pt}\tau M$ du $\mathrm{A}$-module instable $\tau M$.
\hfill$\square$

\medskip
\begin{rem}\label{A-Gysin-bis} Les deux $\mathrm{A}$-modules instables  $\ker e$ et $\mathop{\mathrm{coker}} e$ introduits en section~\ref{Gysin} s'identifient respectivement à $\mathbb{F}_{2}\otimes_{\mathrm{P}}\mathrm{H}^{*}X$ et $\Sigma^{-1}\mathrm{Tor}_{1}^{\mathrm{P}}\hspace{1pt}(\mathbb{F}_{2},\mathrm{H}^{*}X)$. La suite exacte de  du scholie~\ref{A-Gysin} peut donc  être réécrite sous la forme
$$
\hspace{24pt}
0\to
\mathrm{Tor}_{0}^{\mathrm{P}}\hspace{1pt}(\mathbb{F}_{2},\mathrm{H}^{*}X)\to
\mathrm{H}^{*}Y\to\Sigma^{-1}\mathrm{Tor}_{1}^{\mathrm{P}}\hspace{1pt}(\mathbb{F}_{2},\mathrm{H}^{*}X)
\to 0
\hspace{23pt}.
$$
\end{rem}

\subsect{Changement de paradigme}\label{chdp}

\begin{pro}\label{chdp-1} Soient $M_{1}$ et $M_{2}$ deux $\mathrm{P}$-$\mathrm{A}$-modules instables~; on a un isomorphisme de $\mathrm{A}$-modules instables
$$
\hspace{24pt}
M_{1}\otimes_{\mathrm{P}}M_{2}
\hspace{4pt}\cong\hspace{4pt}
\mathbb{F}_{2}\otimes_{\mathrm{P}}(M_{1}\otimes M_{2})
\hspace{23pt},
$$
naturel en $M_{1}$ et $M_{2}$.

\smallskip
{\em (Le $\mathrm{P}$-$\mathrm{A}$-module instable $M_{1}\otimes M_{2}$ qui apparaît ci-dessus est défini en \ref{P-A-prod-1}).}
\end{pro}

\medskip
\textit{Démonstration.} Soient $x_{1}$ un élément de $M_{1}$ et  $x_{2}$ un élément de $M_{2}$, par définition même de la structure de $\mathrm{P}$-$\mathrm{A}$-module instable de $M_{1}\otimes M_{2}$ on a $u\hspace{1pt}(x_{1}\otimes x_{2})=u\hspace{1pt}x_{1}\otimes x_{2}+x_{1}\otimes u\hspace{1pt}x_{2}$ ($\psi(u)=u\otimes 1+1\otimes u$)~; la proposition résulte simplement de cette égalité. Précisons un peu~:

\smallskip
Soit $R$ le sous-$\mathbb{F}_{2}$-espace vectoriel $\mathbb{N}$-gradué de $M_{1}\otimes M_{2}$ engendré par les $u\hspace{1pt}x_{1}\otimes x_{2}+x_{1}\otimes u\hspace{1pt}x_{2}$, $x_{1}$ parcourant $M_{1}$ et $x_{2}$ parcourant~$M_{2}$.

\smallskip
-- L'égalité $R=u\hspace{1pt}(M_{1}\otimes M_{2})$ montre que $R$ est stable sous l'action de $\mathrm{A}$ et que le $\mathrm{A}$-module instable sous-jacent à $\mathbb{F}_{2}\otimes_{\mathrm{P}}(M_{1}\otimes M_{2})$ est le quotient $(M_{1}\otimes M_{2})/R$.

\smallskip
-- On a dans $M_{1}\otimes M_{2}$ la congruence $u\hspace{1pt}x_{1}\otimes x_{2}\equiv x_{1}\otimes u\hspace{1pt}x_{2} \pmod{R}$ et par récurrence $u^{n}\hspace{1pt}x_{1}\otimes x_{2}\equiv x_{1}\otimes u^{n}\hspace{1pt}x_{2} \pmod{R}$, pour tout $n$ dans $\mathbb{N}$~; $R$ contient donc tous les $u^{n}\hspace{1pt}x_{1}\otimes x_{2}+x_{1}\otimes u^{n}\hspace{1pt}x_{2}$. Cette observation montre que $R$ est l'image de la différence des deux homomorphismes $M_{1}\otimes\mathrm{P}\otimes M_{2}\to M_{1}\otimes M_{2}$ qui apparaissent dans la définition de $M_{1}\otimes_{\mathrm{P}}M_{2}$~;  ceci confirme que $R$ est stable sous l'action de $\mathrm{A}$ et montre que le $\mathrm{A}$-module instable sous-jacent à $M_{1}\otimes_{\mathrm{P}}M_{2}$ est aussi  le quotient $(M_{1}\otimes M_{2})/R$.
\hfill$\square$

\medskip
\begin{rem} Les $\mathrm{A}$-modules instables
$M_{1}\otimes_{\mathrm{P}}M_{2}$ et $\mathbb{F}_{2}\otimes_{\mathrm{P}}(M_{1}\otimes M_{2})$ possèdent tous deux une structure naturelle de $\mathrm{P}$-$\mathrm{A}$-modules instable (voir le point (4) de \ref{P-A-prod-2}). Celle du second est $\mathrm{P}$-triviale (voir \ref{theta}), celle du premier ne l'est pas en général (prendre par exemple $M_{1}=M_{2}=\mathrm{P}$).
\end{rem}

\medskip
\begin{cor}\label{chdp-2} Soient $M_{1}$ et $M_{2}$ deux $\mathrm{P}$-$\mathrm{A}$-modules instables~; on a un isomorphisme de $\mathrm{A}$-modules instables
$$
\hspace{24pt}
\mathrm{Tor}^{\mathrm{P}}_{1}(M_{1},M_{2})
\hspace{4pt}\cong\hspace{4pt}
\mathrm{Tor}^{\mathrm{P}}_{1}(\mathbb{F}_{2},M_{1}\otimes M_{2})
\hspace{23pt},
$$
naturel en $M_{1}$ et $M_{2}$.
\end{cor}

\medskip
\textit{Démonstration.}  Soit $M_{1}\leftarrow C_{\bullet}$ est une résolution projective dans la catégorie $\mathrm{P}\text{-}\mathcal{U}$.

\begin{lem}\label{res-P-lib-1} Le complexe augmenté $M_{1}\otimes M_{2}\leftarrow C_{\bullet}\otimes M_{2}$  est une résolution de $M_{1}\otimes M_{2}$ dans la catégorie $\mathrm{P}\text{-}\mathcal{U}$ qui est libre dans la catégorie $\mathrm{P}\text{-}\mathcal{E}$.
\end{lem}

\medskip
Supposons ce lemme démontré. La proposition \ref{res-P-lib} dit que l'on a un isomorphisme de $\mathrm{P}$-$\mathrm{A}$-modules instables $\mathrm{Tor}^{\mathrm{P}}_{1}(M_{1},M_{2})\cong\mathrm{H}_{1}(C_{\bullet}\otimes_{\mathrm{P}}M_{2})$ et la proposition \ref{chdp-1}  que l'on a un isomorphisme de $\mathrm{A}$-modules instables $\mathrm{H}_{1}(C_{\bullet}\otimes_{\mathrm{P}}M_{2})\cong\mathrm{H}_{1}(\mathbb{F}_{2}\otimes_{\mathrm{P}}(C_{\bullet}\otimes M_{2}))$. On achève en invoquant à nouveau la proposition \ref{res-P-lib}.
\hfill$\square$

\medskip
\textit{Démonstration du lemme \ref{res-P-lib-1}.} Soit $k$ un entier naturel~; compte tenu du point (3) (ou (2)) de \ref{Tor-P} on peut supposer que $C_{k}$ est somme directe de $\mathrm{P}\otimes\mathrm{F}(n)$'s. Il suffit donc de vérifier que le produit tensoriel de $\mathrm{P}$-$\mathrm{A}$-modules instables $(\mathrm{P}\otimes\mathrm{F}(n))\otimes M_{2}$ est un $\mathrm{P}$-module libre. Comme le $\mathrm{P}$-$\mathrm{A}$-module instable $\mathrm{P}\otimes\mathrm{F}(n)$ coïncide avec le produit tensoriel de $\mathrm{P}$-$\mathrm{A}$-modules instables $\mathrm{P}\otimes\hspace{1pt}\theta\hspace{1pt}\mathrm{F}(n)$, on a
$$
\hspace{12pt}
(\mathrm{P}\otimes\mathrm{F}(n))\otimes M_{2}
\hspace{2pt}=\hspace{2pt}
(\mathrm{P}\otimes\theta\hspace{1pt}\mathrm{F}(n))\otimes M_{2}
\hspace{2pt}\cong\hspace{2pt}
\mathrm{P}\otimes(\theta\hspace{1pt}\mathrm{F}(n)\otimes M_{2})
\hspace{11pt}.
$$
On conclut à l'aide du lemme suivant~:

\begin{lem}\label{hopf} Soit $M$ un $\mathrm{P}$-$\mathrm{A}$-module instable~; les deux  produits tensoriels de $\mathrm{P}$-$\mathrm{A}$-modules instables $\mathrm{P}\otimes M$ et  $\mathrm{P}\otimes \theta\mathcal{O}M$, $\mathcal{O}$ désignant ici le foncteur oubli $\mathrm{P}\text{-}\mathcal{U}\to\mathcal{U}$, sont canoniquement isomorphes.
\end{lem}

\medskip
\textit{Démonstration} L'argument est classique~: on utilise le fait que $\mathrm{P}$ est une algèbre de Hopf avec antipode (ici l'identité). L'isomorphisme canonique en question est donné par la composition
$$
\begin{CD}
\hspace{12pt}
\mathrm{P}\otimes M@>\psi\hspace{1pt}\otimes\hspace{1pt}1_{M}>>
(\mathrm{P}\otimes\mathrm{P})\otimes M
@>=>>
\mathrm{P}\otimes(\mathrm{P}\otimes M)
@>1_{\mathrm{P}\hspace{1pt}}\otimes\hspace{1pt}\mathrm{a}>>
\mathrm{P}\otimes M
\hspace{11pt},
\end{CD}
$$
$\mathrm{a}:\mathrm{P}\otimes M\to M$ désignant l'action de $\mathrm{P}$ sur $M$.
\hfill$\square$

\medskip
Les énoncés \ref{sigmatau} et \ref{chdp-2} conduisent au lemme suivant que nous utiliserons de façon répétitive par la suite~:

\begin{lem}\label{suspred} Soient $0\to M_{2}\to M'_{2}\to M''_{2}\to 0$ une suite exacte de\linebreak $\mathrm{P}$-$\mathrm{A}$-modules instables et $M_{1}$ un $\mathrm{P}$-$\mathrm{A}$-module instable~; si le $\mathrm{A}$-module instable sous-jacent à $M_{1}\otimes_{\mathrm{P}}M_{2}$ est réduit alors les suites de $\mathrm{P}$-$\mathrm{A}$-modules instables
$$
0\to M_{1}\otimes_{\mathrm{P}}M_{2}
\to M_{1}\otimes_{\mathrm{P}}M'_{2}\to M_{1}\otimes_{\mathrm{P}}M''_{2}\to 0
$$
et 
$$
0\to \mathrm{Tor}^{\mathrm{P}}_{1}(M_{1},M_{2})
\to \mathrm{Tor}^{\mathrm{P}}_{1}(M_{1},M'_{2})
\to \mathrm{Tor}^{\mathrm{P}}_{1}(M_{1},M''_{2})\to 0
$$

sont exactes. On a le même énoncé, {\em mutatis mutandis}, avec le produi tensoriel à droite par $M_{1}$.
\end{lem}

\medskip
\textit{Démonstration.} D'après \ref{chdp-2} et \ref{sigmatau} le $\mathrm{A}$-module instable  sous jacent à $\mathrm{Tor}^{\mathrm{P}}_{1}(M_{1},M''_{2})$ est une suspension; le point (c) de \ref{susp} dit que le connectant $\partial:\mathrm{Tor}^{\mathrm{P}}_{1}(M_{1},M''_{2})\to M_{1}\otimes_{\mathrm{P}}M_{2}$ est trivial.
\hfill$\square$

\subsect{Sur les $\mathrm{P}$-$\mathrm{A}$-modules instables qui sont réduits comme $\mathrm{A}$-modules instables}

\begin{pro}\label{P-A-reduit-1} Soit $M$ un $\mathrm{P}$-$\mathrm{A}$-module instable tel que le $\mathrm{A}$-module instable sous-jacent est réduit.

\medskip
{\em (a)} Soient $x$ un élément de $M$ et $n\geq 1$ un entier~;  les deux conditions suivantes sont équivalentes~:
\begin{itemize}
\item[(i)] $u^{n}\hspace{1pt}x=0$~;
\item[(ii)] $u\hspace{1pt}x=0$.
\end{itemize}
(En d'autres termes, $\tau M$ est la $\mathrm{P}$-torsion du $\mathrm{P}$-module sous-jacent à $M$.)

\medskip
{\em (b)} Le $\mathrm{P}$-module sous-jacent à $M/\tau M$ est libre.

\medskip
{\em (c)} Le $\mathrm{A}$-module instable sous-jacent à $M/\tau M$ est réduit.
\end{pro}

\medskip
\textit{Démonstration du point (a).} La seule implication à vérifier est  $(i)\Rightarrow(ii)$. La formule 
$
\mathrm{Sq}^{i}\hspace{1pt}u^{n}\hspace{1pt}x
=\sum_{j=0}^{i}\hspace{2pt}
{n\choose j}\hspace{1pt}u^{n+j}\hspace{1pt}\mathrm{Sq}^{i-j}x
$
montre par récurrence sur l'entier~$i$ que l'on a $u^{n}\hspace{1pt}\mathrm{Sq}^{i}x=0$ pour tout $i$, et donc en particulier $u^{n}\hspace{1pt}\mathrm{Sq}_{0}\hspace{1pt}x=0$. On écrit $n=2m-\delta$ avec $m\in\mathbb{N}-\{0\}$ et  $\delta\in\{0,1\}$(en d'autres termes, on effectue la division euclidienne de $n+1$ par $2$)~; on a par construction  $u^{2m}\hspace{1pt}\mathrm{Sq}_{0}\hspace{1pt}x=0$ soit encore $\mathrm{Sq}_{0}(u^{m}\hspace{1pt}x)=0$ et donc $u^{m}\hspace{1pt}x=0$ puisque $M$ est réduit. Or on constate que l'on a $m<n$ pour $n\not=1$.
\hfill$\square$

\medskip
\textit{Remarque.} Dans le cas où $M$ est une $\mathrm{P}$-$\mathrm{A}$-algèbre instable,  la démonstration est moins alambiquée~: $u^{n}x=0\Rightarrow(ux)^{n}=0$ et $(ux)^{n}=0\Leftrightarrow ux=~0$ car l'hypothèse ``le $\mathrm{A}$-module instable sous-jacent à $M$ est réduit équivaut à l'hypothèse ``la $\mathbb{F}_{2}$-algèbre   commutative $\mathbb{N}$-graduée sous-jacente à $M$ est réduite''.

\medskip
Le point (b) est conséquence du point (a) et du lemme \ref{P-lib-2} ci-après.
\hfill$\square$

\medskip
\textit{Démonstration du point (c).} Soient $x$ un élément de $M$ et $y$ sa classe dans $M/\tau M$. L'égalité $\mathrm{Sq}_{0}\hspace{1pt}y=0$ équivaut à $u\hspace{1pt}\mathrm{Sq}_{0}\hspace{1pt}x=0$~; si cette dernière égalité est satisfaite, on a \textit{a fortiori} $\mathrm{Sq}_{0}(u\hspace{1pt}x)=0$ ($\mathrm{Sq}_{0}(u\hspace{1pt}x)=\mathrm{Sq}_{0}u\hspace{2pt}\mathrm{Sq}_{0}\hspace{1pt}x=u^{2}\hspace{1pt}\mathrm{Sq}_{0}\hspace{1pt}x$), donc $u\hspace{1pt}x=0$ puisque $M$ est réduit c'est-à-dire $y=0$.
\hfill$\square$

\begin{lem}\label{P-lib-2}Soit $R$ un $\mathrm{P}$-espace vectoriel $\mathbb{N}$-gradué~; les deux conditions suivantes sont équivalentes~:
\begin{itemize}
\item[(i)]  $R$ est un $\mathrm{P}$-module sans torsion~;
\item[(ii)] $R$ est  un $\mathrm{P}$-module libre.
\end{itemize}
\end{lem}

\textit{Démonstration de $(i)\Rightarrow(ii)$.} Soient $\bar{B}$ une base (homogène) du $\mathbb{F}_{2}$-espace vectoriel $\mathbb{N}$-gradué $\mathbb{F}_{2}\otimes_{\mathrm{P}}R$ et $B$ un relèvement (homogène) de $\bar{B}$ dans $R$, On~montre que $B$ est une base de $R$, c'est-à-dire que le $(\mathrm{P}\text{-}\mathcal{E})$-morphisme canonique $\beta:L:=\bigoplus_{x\in B}\mathrm{P}\hspace{1pt}x\to R$ est un isomorphisme. Soit $C$ le conoyau de~$\beta$~; par construction  $\mathbb{F}_{2}\otimes_{\mathrm{P}}\beta$ est un isomorphisme, il en résulte $\mathbb{F}_{2}\otimes_{\mathrm{P}}C=0$, le lemme de Nakayama $\mathbb{N}$-gradué (trivial) implique $C=0$. Soit $N$ le noyau de $\beta$~; on considère la $(\mathrm{P}\text{-}\mathcal{E})$-suite exacte $0\to N\to L\to  R\to 0$, celle-ci induit une suite exacte
$$
\hspace{24pt}
\mathrm{Tor}^{\mathrm{P}}_{1}(\mathbb{F}_{2},R)\to\mathbb{F}_{2}\otimes_{\mathrm{P}}N\to\mathbb{F}_{2}\otimes_{\mathrm{P}}L\overset{\cong}{\to}\mathbb{F}_{2}\otimes_{\mathrm{P}}R\to 0
\hspace{23pt}.
$$
On a par hypothèse $\mathrm{Tor}^{\mathrm{P}}_{1}(\mathbb{F}_{2},R)=0$ (invoquer la $(\mathrm{P}\text{-}\mathcal{E})$-version, immédiate, de \ref{sigmatau}) d'où $\mathbb{F}_{2}\otimes_{\mathrm{P}}N=0$ et donc $N=0$.
\hfill$\square$

\begin{scho}\label{projdim} Tour objet de $\mathrm{P}\text{-}\mathcal{E}$ admet une résolution libre de longueur $1$.
\end{scho}

\bigskip
La proposition \ref{P-A-reduit-2} ci-dessous implique la proposition  \ref{P-A-reduit-1}. Cependant la démonstration de \ref{P-A-reduit-2} est plus sophistiquée que celle de  \ref{P-A-reduit-1} car elle utilise la classification des injectifs de la catégorie $\mathrm{P}$-$\mathcal{U}$.

\begin{pro}\label{P-A-reduit-2} Soit $M$ un $\mathrm{P}$-$\mathrm{A}$-module instable~; les deux conditions suivantes sont équivalentes~:
\begin{itemize}
\item[(i)]  $M$ est réduit comme $\mathrm{A}$-module instable~;
\item[(ii)] $M$ se plonge dans un $\mathrm{P}$-$\mathrm{A}$-module instable de la forme $\theta J\oplus(\mathrm{P}\otimes I)$ avec $I$ et $J$ deux  $\mathrm{A}$-modules instables injectifs réduits.
\end{itemize}
\end{pro}

\textit{Démonstration de $(i)\Rightarrow(ii)$}. Soit $i:M\to E$ une enveloppe injective de~$M$ (dans la catégorie $\mathrm{P}$-$\mathcal{U}$). Soit $N$ le nilradical du $\mathrm{A}$-module instable sous-jacent à $E$~; on constate que $N$ est stable par multiplication par $\mathrm{P}$, en d'autres termes que $N$ est un sous-objet de $E$ dans la catégorie  $\mathrm{P}$-$\mathcal{U}$. Par hypothèse on a $i^{-1}(N)=0$ donc $N=0$~:  le $\mathrm{A}$-module instable sous-jacent à $E$ est réduit.

\smallskip
Soit $V$  un $2$-groupe abélien élémentaire, on rappelle ci-dessous la classification des $\mathrm{H}^{*}V$-$\mathcal{U}$-injectifs (objets injectifs dans la catégorie $\mathrm{H}^{*}V$-$\mathcal{U}$) obtenue dans~\cite{LZsmith}~: 

\smallskip
On exhibe (voir \cite[3.2]{LZsmith}) un système de représentants pour les classes d'isomorphisme de $\mathrm{H}^{*}V$-$\mathcal{U}$-injectifs indécomposables
$$
\hspace{24pt}
\{\hspace{1pt}\mathrm{E}(V,W,L,n)\hspace{1pt}\}_{(W,L,n)\hspace{1pt}\in\hspace{1pt}\mathcal{W}\times\mathcal{L}\times\mathbb{N}}
\hspace{23pt},
$$
$\mathcal{W}$ désignant l'ensemble des sous-groupes de $V$ et $\mathcal{L}$ l'ensemble de $\mathcal{U}$-injectifs réduits introduit dans la démonstration de l'implication $(ii)\Rightarrow(iii)$ de \ref{reduit}. Il~en résulte formellement que  pour tout $\mathrm{H}^{*}V$-$\mathcal{U}$-injectif $I$ il existe une unique famille de cardinaux $(a_{W,L,n})_{(W,L,n)\hspace{1pt}\in\hspace{1pt}\mathcal{W}\times\mathcal{L}\times\mathbb{N}}$ telle que $I$ est isomorphe à la somme directe
$$
\hspace{24pt}
\bigoplus_{(L,W,n)}\mathrm{E}(V,W,L,n)^{\oplus\hspace{1pt}a_{W,L,n}}
\hspace{23pt}.
$$

\smallskip
On constate que le $\mathrm{A}$-module instable sous-jacent à  $\mathrm{E}(V,W,L,n)$ est réduit si et seulement si l'on a $n=0$. On en déduit, en prenant $V=\mathbb{Z}/2$, qu'il existe deux familles de cardinaux $(a_{L})_{L\in\mathcal{L}}$ et $(b_{L})_{L\in\mathcal{L}}$ telles que le $\mathrm{P}$-$\mathrm{A}$-module instable $E$ est isomorphe à la somme directe de deux termes
$$
\hspace{24pt}
\theta\hspace{1pt}(\bigoplus_{L\in\mathcal{L}} L^{\oplus\hspace{1pt}a_{L}})
\hspace{24pt}\text{et}\hspace{24pt}
 \mathrm{P}\otimes(\bigoplus_{L\in\mathcal{L}}L^{\oplus\hspace{1pt}b_{L}})
\hspace{23pt}.
$$
\hfill$\square$

\medskip
\begin{rem}\label{ab initio} Soit $J$ un $\mathcal{U}$-injectif réduit, on peut se convaincre \textit{ab initio} de ce que $\theta J$ est un $(\mathrm{P}\text{-}\mathcal{U})$-injectif à l'aide de \ref{sigmatau}~:

\smallskip
Soit $M$ un $\mathrm{P}$-$\mathrm{A}$-module instable~; on a $\mathrm{Hom}_{\mathrm{P}-\mathcal{U}}(M,\theta J)=\mathrm{Hom}_{\mathcal{U}}(\mathbb{F}_{2}\otimes_{\mathrm{P}}M,J)$. Soit $0\to M'\to M\to M''\to 0$ une  $(\mathrm{P}\text{-}\mathcal{U})$-suite exacte. On conclut en appliquant le foncteur exact $\mathrm{Hom}_{\mathcal{U}}(-,J)$ à la $\mathcal{U}$-suite exacte
$$
\mathrm{Tor}^{\mathrm{P}}_{1}(\mathbb{F}_{2},M'')\to\mathbb{F}_{2}\otimes_{\mathrm{P}}M'\to\mathbb{F}_{2}\otimes_{\mathrm{P}}M\to\mathbb{F}_{2}\otimes_{\mathrm{P}}M''\to 0
$$
et en observant que l'on a $\mathrm{Hom}_{\mathcal{U}}(\mathrm{Tor}^{\mathrm{P}}_{1}(\mathbb{F}_{2},M''),J)=0$ compte tenu de \ref{sigmatau} et du point (c) de \ref{susp}.
\end{rem}

\medskip
\begin{cor}\label{P-A-nilf} Soit $M$ un $\mathrm{P}$-$\mathrm{A}$-module instable~; les deux conditions suivantes sont équivalentes~:
\begin{itemize}
\item[(i)]  $M$ est  $\mathcal{N}\hspace{-1,5pt}il$-fermé comme $\mathrm{A}$-module instable~;
\item[(ii)] il existe une suite exacte $0\to M\to \theta J^{0}\oplus(\mathrm{P}\otimes I^{0})\to \theta J^{1}\oplus(\mathrm{P}\otimes I^{1})$ dans la catégorie $\mathrm{P}\text{-}\mathcal{U}$ avec $J^{0}$, $J^{1}$,  $I^{0}$ et $I^{1}$ des $\mathrm{A}$-modules instables injectifs réduits.
\end{itemize}
\end{cor}

\textit{Démonstration de $(i)\Rightarrow(ii)$.} Analogue à celle de l'implication $(i)\Rightarrow(ii)$ du corollaire \ref{critnilf-3}.
\hfill$\square$

\bigskip
On ajoute maintenant à la condition $(i)$ de \ref{P-A-reduit-2} la condition ``$M$ est libre comme $\mathrm{P}$-module''.

\bigskip
\begin{pro}\label{P-A-reduit-P-lib} Soit $M$ un $\mathrm{P}$-$\mathrm{A}$-module instable~; les deux conditions suivantes sont équivalentes~:
\begin{itemize}
\item[(i)]  $M$ est libre comme $\mathrm{P}$-module et réduit comme $\mathrm{A}$-module instable~;
\item[(ii)] $M$ se plonge dans un $\mathrm{P}$-$\mathrm{A}$-module instable de la forme $\mathrm{P}\otimes I$ avec $I$ un $\mathrm{A}$-module instable injectif réduit.
\end{itemize}
\end{pro}

\textit{Démonstration de $(i)\Rightarrow(ii)$.} Soit $i:M\to E$ une enveloppe injective de~$M$ (dans la catégorie $\mathrm{P}$-$\mathcal{U}$). On a vu lors de la démonstration de l'implication  $(i)\Rightarrow(ii)$ de \ref{P-A-reduit-2} que $E$ est de la forme $\theta J\oplus\mathrm{P}\otimes I$ avec $J$ et $I$ des $\mathcal{U}$-injectifs réduits~; si $M$ est un $\mathrm{P}$-module libre on a $i^{-1}(\theta J)=0$ et donc $\theta J=0$.
\hfill$\square$

\bigskip
\begin{rem}
Soit $p:Y\to X$ un revêtement double~; on rappelle que la donnée de $p$ fait de $\mathrm{H}^{*}X$ un $\mathrm{P}$-$\mathrm{A}$-module instable. L'énoncé \ref{cokere} (reformulé en tenant compte de \ref{A-Gysin-bis}) montre que si $\mathrm{H}^{*}X$ et $\mathbb{F}_{2}\otimes_{\mathrm{P}}\mathrm{H}^{*}X$ sont $\mathcal{N}\hspace{-1.5pt}il$-fermés alors $\mathrm{H}^{*}Y$ l'est aussi.
\end{rem}

\medskip
Cette remarque conduit à la définition \textit{ad hoc} ci-desous~:

\medskip
\begin{defi} Un $\mathrm{P}$-$\mathrm{A}$-module instable $M$ sera dit {\em $\mathrm{P}$-$\mathcal{N}\hspace{-1.5pt}il$-fermé} si les $\mathrm{A}$-modules instables sous-jacents à $M$ et  $\mathbb{F}_{2}\otimes_{\mathrm{P}}M$ sont tous deux $\mathcal{N}\hspace{-1.5pt}il$-fermés.
\end{defi}

\medskip
Nous avons tout fait pour avoir~:

\begin{pro}\label{rev-P-nilf} Soit $p:Y\to X$ un revêtement double~; si le $\mathrm{P}$-$\mathrm{A}$-module instable $\mathrm{H}^{*}X$ est $\mathrm{P}$-$\mathcal{N}\hspace{-1.5pt}il$-fermé alors  le $\mathrm{A}$-module instable $\mathrm{H}^{*}Y$est $\mathcal{N}\hspace{-1.5pt}il$-fermé.
\end{pro}

\medskip
\textsc{Conventions.} Soit $M$ un $\mathrm{P}$-$\mathrm{A}$-module instable. Nous  dirons que $M$ est réduit (resp. $\mathcal{N}\hspace{-1.5pt}il$-fermé) si  le $\mathrm{A}$-module instable sous-jacent à $M$ est réduit (resp. $\mathcal{N}\hspace{-1.5pt}il$-fermé). Nous dirons que $M$ est $\mathrm{P}$-libre s'il est libre comme  $\mathrm{P}$-module $\mathbb{N}$-gradué. Ces conventions devraient alléger notre rédaction.

\medskip
\begin{pro}\label{tens-P-nilf} Soient $M_{1}$ et $M_{2}$ deux $\mathrm{P}$-$\mathrm{A}$-modules instables. Si $M_{1}$ et $M_{2}$ sont  $\mathrm{P}$-$\mathcal{N}\hspace{-1.5pt}il$-fermés alors il en est de même pour leur produit tensoriel $M_{1}\otimes M_{2}$.
\end{pro}

\medskip
\textit{Démonstration.} La proposition  \ref{tensnf} dit que $M_{1}\otimes M_{2}$ est $\mathcal{N}\hspace{-1.5pt}il$-fermé. Il reste à montrer que $\mathbb{F}_{2}\otimes_{\mathrm{P}}(M_{1}\otimes M_{2})$ est aussi $\mathcal{N}\hspace{-1.5pt}il$-fermé, soit encore, compte tenu de \ref{chdp-1} à vérifier l'énoncé suivant~:

\begin{pro}\label{tens-P-nilf-1} Soient $M_{1}$ et $M_{2}$ deux $\mathrm{P}$-$\mathrm{A}$-modules instables. Si $M_{1}$ et $M_{2}$ sont  $\mathrm{P}$-$\mathcal{N}\hspace{-1.5pt}il$-fermés alors $M_{1}\otimes_{\mathrm{P}} M_{2}$ est $\mathcal{N}\hspace{-1.5pt}il$-fermé.
\end{pro}

\medskip
Evidemment cet énoncé implique le suivant~:

\begin{pro}\label{tens-P-nilf-2} Soit $M$ un $\mathrm{P}$-$\mathrm{A}$-module instable. Si $M$ est  $\mathrm{P}$-$\mathcal{N}\hspace{-1.5pt}il$-fermé alors $M\otimes_{\mathrm{P}}M$ est $\mathcal{N}\hspace{-1.5pt}il$-fermé.
\end{pro}

\medskip
Réciproquement \ref{tens-P-nilf-2} implique \ref{tens-P-nilf-1}~: prendre $M=M_{1}\oplus M_{2}$ dans \ref{tens-P-nilf-2} et observer que $M_{1}\otimes_{\mathrm{P}} M_{2}$ est facteur direct dans $(M_{1}\oplus M_{2})\otimes_{\mathrm{P}}(M_{1}\oplus M_{2})$. Nous choisissons de démontrer \ref{tens-P-nilf-2} afin de limiter un peu le nombre de notations.

\bigskip
\textit{Démonstration de \ref{tens-P-nilf-2}.} On pose $K:=\tau M$ et $L:=M/K$. Puisque $M$ est $\mathcal{N}\hspace{-1,5pt}il$-fermé  il est \textit{a fortiori} réduit~; les points (b) et (c) de la proposition \ref{P-A-reduit-1} disent respectivement que $L$ est $\mathrm{P}$-libre (d'où le choix de la notation  !) et que $L$ est réduit.

\smallskip
La proposition \ref{P-A-reduit-2} dit que l'on dispose d'un monomorphisme de $\mathrm{P}$-$\mathrm{A}$-modules instables $M\hookrightarrow\theta J\oplus(P\otimes I)$ avec $J$ et $I$ deux $\mathcal{U}$-injectifs réduits. On note $Q$ le conoyau d'un tel monomorphisme~; on observe que $Q$ est réduit (Proposition \ref{critnilf-0}).

\bigskip
Ceci posé, la démonstration de \ref{tens-P-nilf-2} se décompose en six  (petits) pas~:

\bigskip
 (1) On montre que  $K$ est $\mathcal{N}\hspace{-1,5pt}il$-fermé.
 
\medskip
On considère la suite exacte dans la catégorie $\mathrm{P}$-$\mathcal{U}$
 $$
 \hspace{24pt}
 0\to K\to M\to L\to 0
  \hspace{23pt};
 \leqno{S_{(1)}}
 $$
$K$ est $\mathcal{N}\hspace{-1,5pt}il$-fermé parce que $ M$ et $L$ sont respectivement $\mathcal{N}\hspace{-1,5pt}il$-fermé et réduit (invoquer à nouveau la Proposition \ref{critnilf-0}).
 \hfill$\square$
 
 \bigskip
 (2) On montre que $\mathbb{F}_{2}\otimes_{\mathrm{P}}L$ est réduit.
 
 \medskip
 On applique le foncteur $\mathbb{F}_{2}\otimes_{\mathrm{P}}-$ à la suite exacte $S_{(1)}$~; comme $L$ est $\mathrm{P}$-libre on obtient encore une suite exacte dans la catégorie~$\mathcal{U}$~:
 $$
 \hspace{24pt}
 0\to K\to \mathbb{F}_{2}\otimes_{\mathrm{P}}M\to
 \mathbb{F}_{2}\otimes_{\mathrm{P}}L\to 0
  \hspace{23pt}.
 $$
 Comme $\mathbb{F}_{2}\otimes_{\mathrm{P}}M$ est réduit (puisque $\mathcal{N}\hspace{-1,5pt}il$-fermé) et que $K$ est $\mathcal{N}\hspace{-1,5pt}il$-fermé d'après (1),  la longue suite exacte des $\mathrm{Ext}^{*}_{\mathcal{U}}( \mathrm{nilpotent},-)$ montre que $\mathbb{F}_{2}\otimes_{\mathrm{P}}L$ est réduit.
 \hfill$\square$
  
 \bigskip
(3) On montre que $\mathbb{F}_{2}\otimes_{\mathrm{P}}Q$ est réduit.

\medskip
Par définition de $Q$ on dispose d'une seconde suite exacte dans la catégorie $\mathrm{P}$-$\mathcal{U}$, à savoir
$$
\hspace{24pt}
0\to M\to\theta J\oplus(\mathrm{P}\otimes  I)\to Q\to 0
\hspace{23pt}.
\leqno{S_{(3)}}
$$
Puisque $\mathbb{F}_{2}\otimes_{\mathrm{P}}M$ est réduit on peut invoquer le lemme \ref{suspred}~; on a donc une $\mathcal{U}$-suite exacte~:
$$
\hspace{24pt}
0\to\mathbb{F}_{2}\otimes_{\mathrm{P}}M\to J\oplus I\to\mathbb{F}_{2}\otimes_{\mathrm{P}}Q\to 0
\hspace{23pt}.
$$
Comme $J\oplus I$ est réduit  et que $\mathbb{F}_{2}\otimes_{\mathrm{P}}M$ est $\mathcal{N}\hspace{-1,5pt}il$-fermé, $\mathbb{F}_{2}\otimes_{\mathrm{P}}Q$ est réduit (même argument que précédemment).
\hfill$\square$

\bigskip
(4) Soit $R$ un $\mathrm{P}$-$\mathrm{A}$-module instable réduit,  on montre que $L\otimes_{P}R$ est réduit. On en déduit en particulier que $L\otimes_{P}Q$ et $L\otimes_{P}M$ sont réduits.

\medskip
D'après \ref{P-A-reduit-2}, il existe un monomorphisme de  $\mathrm{P}$-$\mathrm{A}$-modules instables
$$
\hspace{24pt}
i_{R}:R\to\theta J_{R}\oplus(\mathrm{P}\otimes I_{R})
\hspace{23pt},
$$
avec $J_{R}$ et $I_{R}$ deux $\mathcal{U}$-injectifs réduits. Puisque $L$ est $\mathrm{P}$-libre, $L\otimes_{P}i_{R}$ est encore un monomorphisme. Il suffit donc de vérifier que le $\mathrm{A}$-module instable $L\otimes_{P}(\theta J_{R}\oplus(\mathrm{P}\otimes I_{R}))$ est réduit. Or celui-ci est somme directe de deux termes, $(\mathbb{F}_{2}\otimes_{\mathrm{P}}L)\otimes J_{R}$ et $L\otimes I_{R}$, qui sont tous deux réduits parce que les $\mathrm{A}$-modules instables $(\mathbb{F}_{2}\otimes_{\mathrm{P}}L)$, $J_{R}$, $L$ et  $I_{R}$ le sont  (le premier d'après (2)).

\bigskip
(5) Soit $R$ un $\mathrm{P}$-$\mathrm{A}$-module instable tel que $R$ et $\mathbb{F}_{2}\otimes_{\mathrm{P}}R$ sont réduits~; on montre que $M\otimes_{P}R$ est réduit. On en déduit en particulier que $M\otimes_{P}Q$ et  $M\otimes_{P}M$ sont réduits.

\medskip
On applique le foncteur $-\otimes_{\mathrm{P}}R$ à la suite exacte $S_{(1)}$. Comme $L$ est $\mathrm{P}$-libre, on obtient une suite exacte
$$
\hspace{24pt}
0\to K\otimes(\mathbb{F}_{2}\otimes_{\mathrm{P}}R)\to M\otimes_{P}R\to  L\otimes_{P}R\to 0
\hspace{23pt};
$$
$K\otimes(\mathbb{F}_{2}\otimes_{\mathrm{P}}R)$ est réduit, car $K$ et $\mathbb{F}_{2}\otimes_{\mathrm{P}}R$ le sont, et $L\otimes_{P}R$ est réduit d'après (4) du coup $M\otimes_{P}R$ est aussi réduit.
\hfill$\square$

\bigskip
(6) On montre (enfin~!) que $M\otimes_{P}M$ est $\mathcal{N}\hspace{-1,5pt}il$-fermé.

\medskip
On applique le foncteur $M\otimes_{\mathrm{P}}-$ à la suite exacte $S_{(3)}$. Puisque $M\otimes_{\mathrm{P}}M$ est réduit on peut encre invoquer le lemme \ref{suspred}~; on obtient une $\mathcal{U}$-suite exacte
$$
\hspace{24pt}
0\to M\otimes_{P}M\to((\mathbb{F}_{2}\otimes_{\mathrm{P}}M)\otimes J)\oplus(M\otimes I)\to M\otimes_{P}Q\to 0
\hspace{23pt}.
$$

On conclut en observant que $(\mathbb{F}_{2}\otimes_{\mathrm{P}}M)\otimes J$ et $M\otimes I$ sont $\mathcal{N}\hspace{-1,5pt}il$-fermés d'après la proposition \ref{tensnf} (un  $\mathcal{U}$-injectif réduit est $\mathcal{N}\hspace{-1,5pt}il$-fermé~!) et que  $M\otimes_{P}Q$  est réduit d'après (5).
\hfill$\square\square$

\bigskip
\textsc{Complément}

\medskip
Bien que la démonstration de \ref{tens-P-nilf-2} soit terminée on se propose de faire un pas de plus~:

\medskip
(7) On montre que $L$ est $\mathcal{N}\hspace{-1,5pt}il$-fermé.

\medskip
Soient $i':K\to\theta J$ et $i'':L\to\mathrm{P}\otimes I$ des $(\mathrm{P}\text{-}\mathcal{U})$-monomorphismes avec $J$ et $I$ des $\mathcal{U}$-injectif réduits (pour se convaincre de l'existence de $i'$ on peut par exemple invoquer \ref{reduit} et \ref{ab initio}). Un prolongement de $i'$ à $M$ induit un $(\mathrm{P}\text{-}\mathcal{U})$-monomorphisme $i:M\to\theta J\oplus(\mathrm{P}\otimes I)$ qui prend place dans le $(\mathrm{P}\text{-}\mathcal{U})$-diagramme commutatif ci-dessous

\footnotesize
$$
\begin{CD}
& & 0 & & 0 & & 0 \\
& & @VVV @VVV @VVV \\
0@>>>K@>>>M@>>>L@>>>0 \\
& & @VVi'V @VViV @VVi''V \\
0@>>>\theta J@>>>\theta J\oplus(\mathrm{P}\otimes I)@>>>\mathrm{P}\otimes I@>>>0 \\
& & @VVV @VVV @VVV \\
0@>>>Q'@>>>Q@>>>Q''@>>>0 \\
& & @VVV @VVV @VVV \\
& & 0 & & 0 & & 0 
\end{CD}
$$
\normalsize
dans lequel les flèches horizontales de part et d'autre du terme $\theta J\oplus(\mathrm{P}\otimes I)$ sont l'inclusion et la projection canoniques et les termes $Q'$, $Q$, $Q''$ sont les conoyaux respectifs des flèches $i'$, $i$, $i''$. Par construction les lignes et les colonnes de ce diagramme sont exactes.

\smallskip
Compte tenu du point (2) de la démonstration de \ref{tens-P-nilf-2} on peut appliquer le lemme \ref{suspred} à la suite exacte $0\to L\to\mathrm{P}\otimes I\to Q''\to 0$. Il en résulte en particulier $\mathrm{Tor}^{\mathrm{P}}_{1}(\mathbb{F}_{2},Q'')=0$ ce qui implique que $Q''$ est $\mathrm{P}$-libre (voir \ref{sigmatau} et \ref{P-lib-2}).

\smallskip
L'hypothèse $``M$ est $\mathcal{N}\hspace{-1.5pt}il$-fermé'' équivaut à ``$Q$ est réduit''. Comme $Q'$ est $\mathrm{P}$-trivial (puisque quotient de $\theta J$) et que $Q''$ est $\mathrm{P}$-libre, on a $Q'\cong \tau Q$ (le~foncteur $\tau$ est exact à gauche~!). Le point (c) de \ref{P-A-reduit-1} dit que $Q''$ est réduit ce qui entraîne que $L$ est $\mathcal{N}\hspace{-1.5pt}il$-fermé.
\hfill$\square$

\medskip
\begin{rem}\label{helas} En général $\mathbb{F}_{2}\otimes_{\mathrm{P}}L$ n'est pas $\mathcal{N}\hspace{-1.5pt}il$-fermé. Voir \ref{contrexples}.
\end{rem}

\sect{Sur la cohomologie modulo $2$ des groupes alternés}\label{Alt}

L'objet de cette section est de montrer que la cohomologie modulo  $2$ d'un $2$-sous-groupe de Sylow d'un groupe alterné est $\mathcal{N}\hspace{-1.5pt}il$-fermée.

\subsect{Stratégie et énoncés}

Soit  $n\geq 2$ un entier. On note $\epsilon_{n}:\mathfrak{S}_{n}\to\{\pm 1\}\cong\mathbb{Z}/2$ l'homomorphisme signature dont le groupe alterné $\mathfrak{A}_{n}$ est le noyau. Il est clair que la restriction de $\epsilon_{n}$ au $2$-Sylow $\mathrm{S}_{n}$ (décrit en \ref{2-Sylow-1}) 
de $\mathfrak{S}_{n}$ est non-triviale~; on la note toujours $\epsilon_{n}:\mathrm{S}_{n}\to\mathbb{Z}/2$. On pose
$$
\hspace{24pt}
\mathrm{A}_{n}
\hspace{4pt}:=\hspace{4pt}
\ker\hspace{2pt}(\epsilon_{n}:\mathrm{S}_{n}\to\mathbb{Z}/2)
\hspace{23pt};
$$
$\mathrm{A}_{n}$ est un $2$-Sylow de $\mathfrak{A}_{n}$.

\medskip
Incidemment une remarque~:

\smallskip
\begin{rem}\label{Klein}
On constate que l'on a $\mathrm{A}_{4}=[\mathfrak{A}_{4},\mathfrak{A}_{4}]$ (l'abélianisé de $\mathfrak{A}_{4}$ est isomorphe à $\mathbb{Z}/3$), en d'autres termes que $\mathrm{A}_{4}$ est le {\em groupe de Klein}  (sous-groupe de $\mathfrak{S}_{4}$). Nous adopterons  la notation $\mathrm{A}_{4}$ pour le groupe de Klein dans l'appendice \ref{appen-C}. 
\end{rem}

 \smallskip
 Le résultat que nous avons en vue est donc le suivant~:

\begin{theo}\label{alt-1} Le $\mathrm{A}$-module instable
$\mathrm{H}^{*}\mathrm{A}_{n}$ est $\mathcal{N}\hspace{-1,5pt}il$-fermé pour tout entier~$n\geq 2$.
\end{theo}

\medskip
Compte tenu de \ref{Sylow-2} ce théorème implique~:

\begin{cor}\label{quillenalt}
L'application de Quillen $\mathrm{q}_{\mathfrak{A}_{n}}:\mathrm{H}^{*}(\mathfrak{A}_{n}\;\mathbb{F}_{2})\to\mathrm{L}(\mathfrak{A}_{n})$ est un isomorphisme pour tout entier $n\geq 2$.
\end{cor}

\medskip
On en revient maintenant à l'énoncé \ref{alt-1}. L'espace classifiant $\mathrm{B}\mathrm{A}_{n}$ a le type d'homotopie du revêtement double de $\mathrm{B}\mathrm{S}_{n}$ dont la classe caractéristique est l'élément de $\mathrm{H}^{1}\mathrm{S}_{n}:=\mathrm{H}^{1}(\mathrm{S}_{n};\mathbb{F}_{2})=\mathrm{H}^{1}(\mathrm{S}_{n};\mathbb{Z}/2)=\mathrm{Hom}(\mathrm{S}_{n},\mathbb{Z}/2)$ donné par l'homomorphisme $\epsilon_{n}:\mathrm{S}_{n}\to\mathbb{Z}/2$.

\medskip
Notre stratégie pour démontrer \ref{alt-1} est d'utiliser les observations faites en section \ref{Gysin}.

\medskip
Le $\mathcal{K}$-morphisme $\mathrm{P}=\mathrm{H}^{*}\mathbb{Z}/2\to\mathrm{H}^{*}\mathrm{S}_{n}$ munit $\mathrm{H}^{*}\mathrm{S}_{n}$ d'une structure canonique de $\mathrm{P}$-$\mathrm{A}$-algèbre instable et \textit{a fortiori} de $\mathrm{P}$-$\mathrm{A}$-module instable. Comme l'on sait déjà que le $\mathrm{A}$-module instable $\mathrm{H}^{*}\mathrm{S}_{n}$ est $\mathcal{N}\hspace{-1,5pt}il$-fermé (Théorème \ref{symnf}), l'énoncé \ref{cokere} (et la remarque \ref{A-Gysin-bis}) montrent que l'énoncé \ref{alt-1} est équivalent aux suivants~:

\begin{theo}\label{alt-2} Le $\mathrm{A}$-module instable
$\mathbb{F}_{2}\otimes_{\mathrm{P}}\mathrm{H}^{*}\mathrm{S}_{n}$ est $\mathcal{N}\hspace{-1,5pt}il$-fermé pour tout entier~$n\geq 2$.
\end{theo}

\begin{theo}\label{alt-3} Le $\mathrm{P}$-$\mathrm{A}$-module instable
$\mathrm{H}^{*}\mathrm{S}_{n}$ est $\mathrm{P}$-$\mathcal{N}\hspace{-1.5pt}il$-fermé pour tout entier~$n\geq 2$.
\end{theo}

\medskip
C'est ce dernier énoncé que l'on va démontrer  en remplaçant dans les arguments utilisés en \ref{Sn} les $\mathrm{A}$-modules instables par les $\mathrm{P}$-$\mathrm{A}$-modules instables.

\bigskip
On reprend les notations de  \ref{2-Sylow-1}~; 
on observera que si $n$ est impair alors on a $\mathrm{S}_{n}=\mathrm{S}_{n-1}$ si bien que l'on peut supposer $n$ pair et ou encore $m_{r}\geq 1$. L'isomorphisme de $\mathrm{A}$-modules instables
$$
\hspace{24pt}
\mathrm{H}^{*}\mathrm{S}_{n}
\hspace{4pt}\cong\hspace{4pt}
\mathrm{H}^{*}\mathrm{S}_{2^{m_{1}}}\otimes\mathrm{H}^{*}\mathrm{S}_{2^{m_{2}}}\otimes\ldots\otimes\mathrm{H}^{*}\mathrm{S}_{2^{m_{r}}}
\hspace{24pt}
$$
est en fait un isomorphisme de $\mathrm{P}$-$\mathrm{A}$-modules instables. En effet l'homomorphisme
$\epsilon_{n}:\mathrm{S}_{n}=\mathrm{S}_{2^{m_{1}}}\times\mathrm{S}_{2^{m_{1}}}\times\ldots\times\mathrm{S}_{2^{m_{r}}}\to\mathbb{Z}/2$
est le composé de l'homomorphisme
$$
\epsilon_{2^{m_{1}}}\times\epsilon_{2^{m_{2}}}\times\ldots\times\epsilon_{2^{m_{r}}}:\mathrm{S}_{2^{m_{1}}}\times\mathrm{S}_{2^{m_{1}}}\times\ldots\times\mathrm{S}_{2^{m_{r}}}\to\mathbb{Z}/2\times\mathbb{Z}/2\times\ldots\times\mathbb{Z}/2
$$
et de l'homorphisme ``addition'' $\mathbb{Z}/2\times\mathbb{Z}/2\times\ldots\times\mathbb{Z}/2\to\mathbb{Z}/2$.

\medskip
La proposition \ref{tens-P-nilf} montre donc qu'il suffit  de démontrer le théorème \ref{alt-3} pour $n=2^{m}$ avec $m\geq 1$.

\subsect{Constructions quadratiques dans les catégories $\mathrm{P}\backslash\mathcal{K}$ et  $\mathrm{P}$-$\mathcal{U}$}\label{cq-3}

\medskip
On rappelle que l'on a $\mathrm{S}_{2^{m+1}}=\mathfrak{S}_{2}\wr\mathrm{S}_{2^{m}}:=\mathfrak{S}_{2}\rtimes(\mathrm{S}_{2^{m}}\times\mathrm{S}_{2^{m}})$, $\mathfrak{S}_{2}$~agissant sur $\mathrm{S}_{2^{m}}\times\mathrm{S}_{2^{m}}$ par l'échange des facteurs)~; la diagonale $\mathrm{S}_{2^{m}}\to\mathrm{S}_{2^{m}}\times\mathrm{S}_{2^{m}}$ induit une inclusion de groupes
$\mathfrak{S}_{2}\times\mathrm{S}_{2^{m}}\subset\mathrm{S}_{2^{m+1}}$.

\smallskip
On constate que pour $m\geq 1$ la restriction de $\epsilon_{2^{m+1}}$ à $\mathrm{S}_{2^{m}}\times\mathrm{S}_{2^{m}}$ est l'homomorphisme composé
$$
\begin{CD}
\hspace{24pt}
\mathrm{S}_{2^{m}}\times\mathrm{S}_{2^{m}}
@>\epsilon_{2^{m}}\times\epsilon_{2^{m}}>>
\mathbb{Z}/2\times\mathbb{Z}/2
@>\text{addition}>>
\mathbb{Z}/2
\hspace{23pt}.
\end{CD}
$$
 et que sa restriction à $\mathfrak{S}_{2}\times\mathrm{S}_{2^{m}}$ est triviale.
 
 \smallskip
 Il en résulte que l'image de $\epsilon_{2^{m+1}}$ par le $\mathcal{K}$-isomorphisme de \ref{cqgroupe} 
 $$
 \mathrm{H}^{*}\mathrm{S}_{2^{m+1}}\cong
 \mathfrak{S}_{2}\mathrm{H}^{*}\mathrm{S}_{2^{m}}:=
 \lim\hspace{2pt}(\hspace{2pt}(\mathrm{H}^{*}\mathrm{S}_{2^{m}}\otimes\mathrm{H}^{*}\mathrm{S}_{2^{m}})^{\mathfrak{S}_{2}}\overset{\nu}{\longrightarrow}\Phi\mathrm{H}^{*}\mathrm{S}_{2^{m}}\overset{\rho}{\longleftarrow}\mathrm{R}_{1}\mathrm{H}^{*}\mathrm{S}_{2^{m}}\hspace{2pt})
 $$
 est le couple $(\epsilon_{2^{m}}\otimes 1+1\otimes\epsilon_{2^{m}},0)$ de $(\mathrm{H}^{*}\mathrm{S}_{2^{m}}\otimes\mathrm{H}^{*}\mathrm{S}_{2^{m}})^{\mathfrak{S}_{2}}\times\mathrm{R}_{1}\mathrm{H}^{*}\mathrm{S}_{2^{m}}$.
 
 \bigskip
 L'analyse ci-dessus conduit aux définitions \ref{cq-4} et \ref{cq-5} ci-après.
 
 \medskip
 On commence par décrire une variante pour les  $\mathrm{A}$-algèbres instables de la définition \ref{theta}.

\medskip
\begin{defi}\label{theta-bis} Soit $K$ une $\mathrm{A}$-algèbre instable, l'homomorphisme $\mathrm{P}\to K$ composé de l'augmentation $\mathrm{P}\to\mathbb{F}_{2}$ et de l'unité $\mathbb{F}_{2}\to K$ permet de faire de $K$ une $\mathrm{P}$-$\mathrm{A}$-algèbre instable, une telle $\mathrm{P}$-$\mathrm{A}$-algèbre instable, sera dit {\em$\mathrm{P}$-triviale}. On note enore $\theta:\mathcal{K}\to\mathrm{P}\backslash\mathcal{K}$ le foncteur ainsi défini.
\end{defi}

\medskip
Il est clair que les foncteurs $\theta:\mathcal{K}\to\mathrm{P}\backslash\mathcal{K}$ et $\theta:\mathcal{U}\to\mathrm{P}\text{-}\mathcal{U}$ sont en un sens évident ``oubli-compatibles''.

\bigskip
\begin{defi}\label{cq-4} Soit $(K;f)$ une $\mathrm{P}$-$\mathrm{A}$-algèbre instable $K$ c'est-à-dire une $\mathrm{A}$-algèbre instable $K$ munie d'un $\mathcal{K}$-morphisme $f:\mathrm{P}\to K$.

\smallskip
On munit la $\mathrm{A}$-algèbre instable $(K\otimes K)^{\mathfrak{S}_{2}}$ du $\mathcal{K}$-morphisme $\mathrm{P}\to(K\otimes K)^{\mathfrak{S}_{2}}$ composé de la diagonale $\psi:\mathrm{P}\to(\mathrm{P}\otimes\mathrm{P})^{\mathfrak{S}_{2}}$ (induite par la loi de groupe de~$\mathbb{Z}/2$) et du $\mathcal{K}$-morphisme $(f\otimes f)^{\mathfrak{S}_{2}}$ induit par $f$.

\smallskip
On note encore $\nu:(K\otimes K)^{\mathfrak{S}_{2}}\to\theta\hspace{1pt}\Phi K$ le $\mathrm{P}\backslash\mathcal{K}$-morphisme induit par le $\mathcal{K}$-morphisme $\nu$ de \ref{cqalg} (observer que tout $\mathcal{K}$-morphisme de $\mathrm{P}$ dans $\Phi K$ est forcément trivial).

\smallskip
 On pose
$$
\hspace{24pt}
\mathfrak{S}_{2}(K;f):=
\lim\hspace{2pt}(\hspace{2pt}(K\otimes K)^{\mathfrak{S}_{2}}\overset{\nu}{\longrightarrow}\theta\hspace{1pt}\Phi K\overset{\theta\rho}{\longleftarrow}\theta\hspace{1pt}\mathrm{R}_{1}K\hspace{2pt})
\hspace{23pt},
$$
limite (produit fibré) dans la catégorie $\mathrm{P}\backslash\mathcal{K}$ (soit $e:=f(u)$, nous avons vu en \ref{P-A-notions} 
que la donnée de $e$ est équivalente à celle de $f$
aussi $\mathfrak{S}_{2}(K;f)$ peut être aussi notée $\mathfrak{S}_{2}(K;e)$).

\smallskip
L'endofoncteur $\mathfrak{S}_{2}:\mathrm{P}\backslash\mathcal{K}\circlearrowleft$ ainsi défini et  l'endofoncteur $\mathfrak{S}_{2}:\mathcal{K}\circlearrowleft$  de \ref{cqalg} commutent (en un sens évident) avec le foncteur oubli $\mathrm{P}\backslash\mathcal{K}\to\mathcal{K}$.
\end{defi}

\bigskip
On donne ci-dessous la version ``linéaire'' de la définition ci-dessus. Au préalable deux observations~:

\smallskip
On note $\mathcal{O}:\mathrm{P}\text{-}\mathcal{U}\to\mathcal{U}$ le foncteur oubli. Soit $M$ un un $\mathrm{P}$-$\mathrm{A}$-module instable.
 
\smallskip
1) Le $\mathcal{U}$-automorphisme $\sigma$, $x\otimes y\mapsto y\otimes x$, définissant l'action de $\mathfrak{S}_{2}$ sur $\mathcal{O}M\otimes\mathcal{O}M$ est en fait un $(\mathrm{P}\text{-}\mathcal{U})$-automorphisme du $\mathrm{P}$-$\mathrm{A}$-module instable $M\otimes M$. Il en résulte que $(M\otimes M)^{\mathfrak{S}_{2}}$ est un sous-$\mathrm{P}$-$\mathrm{A}$-module instable de $M\otimes M$.

\smallskip
2) Soit $z$ un élément de  $(M\otimes M)^{\mathfrak{S}_{2}}$~; on constate que $uz$ apartient à l'image de $1+\sigma$. Le $\mathcal{U}$-morphisme $(\mathcal{O}M\otimes\mathcal{O}M)^{\mathfrak{S}_{2}}\to\Phi\mathcal{O}M$ de \ref{cqalg} est donc en fait un $(\mathrm{P}\text{-}\mathcal{U})$-morphisme que l'on notera encore $\nu:(M\otimes M)^{\mathfrak{S}_{2}}\to\theta\hspace{1pt}\Phi\mathcal{O}M$ (ou $\nu_{M}$ si cette précision peut être utile). On observera incidemment que comme $\Phi\mathcal{O}M$ est nul en degrés impairs toute $(\mathrm{P}\text{-}\mathcal{U})$-structure sur $\Phi\mathcal{O}M$ est forcément triviale.

\bigskip
\begin{defi}\label{cq-5} Soit $M$ un $\mathrm{P}$-$\mathrm{A}$-module instable instable~; on pose
$$
\hspace{24pt}
\mathfrak{S}_{2}M:=
\lim\hspace{2pt}((M\otimes M)^{\mathfrak{S}_{2}}\overset{\nu}{\longrightarrow}\theta\hspace{1pt}\Phi\mathcal{O}M\overset{\theta\rho}{\longleftarrow}\theta\hspace{1pt}\mathrm{R}_{1}\mathcal{O}M\hspace{2pt})
\hspace{23pt},
$$
limite (produit fibré) dans la catégorie $\mathrm{P}\text{-}\mathcal{U}$. 

\smallskip
L'endofoncteur $\mathfrak{S}_{2}:\mathrm{P}\text{-}\mathcal{U}\circlearrowleft$ ainsi défini et  l'endofoncteur $\mathfrak{S}_{2}:\mathcal{U}\circlearrowleft$  de \ref{cqalg} commutent (en un sens évident) avec le foncteur oubli $\mathcal{O}:\mathrm{P}\text{-}\mathcal{U}\to\mathcal{U}$. De même l'endofoncteur $\mathfrak{S}_{2}:\mathrm{P}\backslash\mathcal{K}\circlearrowleft$ et l'endofoncteur $\mathfrak{S}_{2}:\mathrm{P}\text{-}\mathcal{U}\circlearrowleft$ commute avec le foncteur oubli $\mathrm{P}\backslash\mathcal{K}\to\mathrm{P}\text{-}\mathcal{U}$.
\end{defi}

\bigskip
A nouveau nous avons tout fait  pour avoir~:

\medskip
\begin{pro}\label{cg-6} Pour tout entier $m\geq 1$\ on a un isomorphisme canonique de $\mathrm{P}$-$\mathrm{A}$-algèbres instables et {\em a fortiori} de $\mathrm{P}$-$\mathrm{A}$-modules instables
$$
\hspace{24pt}
\mathrm{H}^{*}\mathrm{S}_{2^{m+1}}
\hspace{4pt}\cong\hspace{4pt}
\mathfrak{S}_{2}\mathrm{H}^{*}\mathrm{S}_{2^{m}}
\hspace{23pt}.
$$
\end{pro}

Disposant de cet isomorphisme on achève la démonstration du théorème \ref{alt-3} dans le cas $n=2^{m}$ par récurrence sur $m$ ($\mathrm{H}^{*}\mathrm{S}_{2^{1}}=\mathrm{P}$ est un  $\mathrm{P}$-$\mathrm{A}$-module instable $\mathrm{P}$-$\mathcal{N}\hspace{-1.5pt}il$-fermé), à l'aide de la proposition suivante~:

\begin{pro}\label{quad-P-nilf} Soit $M$ un $\mathrm{P}$-$\mathrm{A}$-module instable~; si $M$ est $\mathrm{P}$-$\mathcal{N}\hspace{-1.5pt}il$-fermé  alors il en est de même pour $\mathfrak{S}_{2}M$.
\end{pro}

\medskip
Puisque les endofoncteurs $\mathfrak{S}_{2}:\mathrm{P}\text{-}\mathcal{U}\circlearrowleft$ et $\mathfrak{S}_{2}:\mathcal{U}\circlearrowleft$  de \ref{cqalg} sont``oubli-compatibles'', la proposition \ref{quadnilf} dit que le $\mathrm{A}$-module instable $\mathfrak{S}_{2}M$ est $\mathcal{N}\hspace{-1.5pt}il$-fermé~; la proposition ci-dessus est donc impliquée par la suivante~:

\begin{pro}\label{quad-P-nilf-bts} Soit $M$ un $\mathrm{P}$-$\mathrm{A}$-module instable~; si $M$ est $\mathrm{P}$-$\mathcal{N}\hspace{-1.5pt}il$-fermé  alors le $\mathrm{A}$-module instable $\mathbb{F}_{2}\otimes_{\mathrm{P}}\mathfrak{S}_{2}M$ est $\mathcal{N}\hspace{-1.5pt}il$-fermé.
\end{pro}

\medskip
Nous allons ramener la démonstration de  \ref{quad-P-nilf-bts} à celle de \ref{proclef-bis}. Les énoncés intermédiaires \ref{S1} et \ref{S1-bis}  ne font pas intervenir l'hypothèse ``$M$ est $\mathrm{P}$-$\mathcal{N}\hspace{-1.5pt}il$-fermé''.

\medskip
Soit $M$ un  $\mathrm{P}$-$\mathrm{A}$-module instable. On a par définition même de la construction quadratique $\mathfrak{S}_{2}M$, une suite exacte de $\mathrm{P}$-$\mathrm{A}$-modules instables 
$$
\hspace{24pt}
0\to\mathfrak{S}_{2}M\to(M\otimes M)^{\mathfrak{S}_{2}}\oplus\theta\hspace{1pt}\mathrm{R}_{1}\mathcal{O}M\to\theta\hspace{1pt}\Phi \mathcal{O}M\to 0
\hspace{24pt}
$$
(l'avant dernière flèche est surjective parce que, par exemple, le $\mathcal{U}$-morphisme $\mathrm{R}_{1}\mathcal{O}M\to\Phi \mathcal{O}M$ est surjectif). En appliquant le foncteur $\mathbb{F}_{2}\otimes_{\mathrm{P}}-$ à cette suite exacte on obtient la suite exacte de $\mathrm{A}$-modules instables suivante~:
\begin{multline*}
\hspace{24pt}
0\leftarrow
\Phi\mathcal{O}M
\leftarrow
\mathbb{F}_{2}\otimes_{\mathrm{P}}(M\otimes M)^{\mathfrak{S}_{2}}
\oplus\mathrm{R}_{1}\mathcal{O}M
\leftarrow
\mathbb{F}_{2}\otimes_{\mathrm{P}}\mathfrak{S}_{2}M
\overset{\partial}{\leftarrow}\Sigma\hspace{1pt}\Phi \mathcal{O}M
\leftarrow
\\
\Sigma\hspace{1pt}\tau\hspace{1pt}((M\otimes M)^{\mathfrak{S}_{2}})
\oplus\Sigma\hspace{1pt}\mathrm{R}_{1}\mathcal{O}M
\leftarrow
\Sigma\hspace{1pt}\tau(\mathfrak{S}_{2}M)
\leftarrow 0
\hspace{24pt}
\end{multline*}
(se rappeler de \ref{sigmatau}). Comme l'homomorphisme $\Sigma\mathrm{R}_{1}\mathcal{O}M\to\Sigma\Phi\mathcal{O}M$ est surjectif le connectant $\partial$ qui apparaît ci-dessus est trivial. On dispose donc de deux suites exactes courtes de $\mathrm{A}$-modules instables~: 
$$
\leqno{(\mathcal{S}_{1})}
\hspace{12pt}
0\to\mathbb{F}_{2}\otimes_{\mathrm{P}}\mathfrak{S}_{2}M
\to\mathbb{F}_{2}\otimes_{\mathrm{P}}(M\otimes M)^{\mathfrak{S}_{2}}
\oplus\mathrm{R}_{1}\mathcal{O}M\to\Phi\mathcal{O}M\to 0
\hspace{12pt}
$$
et
$$
\leqno{(\mathcal{S}_{2})}
\hspace{24pt}
0\to\tau(\mathfrak{S}_{2}M)
\to\tau\hspace{1pt}((M\otimes M)^{\mathfrak{S}_{2}})
\oplus\mathrm{R}_{1}\mathcal{O}M\to\Phi\mathcal{O}M\to 0
\hspace{23pt}.
$$
Nous reviendrons sur $(\mathcal{S}_{2})$ dans l'appendice \ref{appen-A}~; l'exactitude de $(\mathcal{S}_{1})$ implique~:

\bigskip
\begin{pro}\label{S1} Soit $M$ un $\mathrm{P}$-$\mathrm{A}$-module instable. Le $\mathrm{A}$-module instable $\mathbb{F}_{2}\otimes_{\mathrm{P}}\mathfrak{S}_{2}M$ est naturellement isomorphe à la limite (produit fibré) du $\mathcal{U}$-~diagramme
$$
\hspace{24pt}
\begin{CD}
\mathbb{F}_{2}\otimes_{\mathrm{P}}(M\otimes M)^{\mathfrak{S}_{2}}
@>\mathbb{F}_{2}\otimes_{\mathrm{P}}\nu>>
\Phi\mathcal{O}M
@<\mathcal{O}\rho<<
\mathrm{R}_{1}\mathcal{O}M
\end{CD}
\hspace{23pt}.
$$
\end{pro}

\medskip
Comme la flèche $\mathbb{F}_{2}\otimes_{\mathrm{P}}\nu$ est surjective (parce que le produit tensoriel $\mathbb{F}_{2}\otimes_{\mathrm{P}}-$  préserve les surjections) on a~:

\medskip
\begin{scho}\label{S1-bis} Soit $M$ un $\mathrm{P}$-$\mathrm{A}$-module instable. On a  une $\mathcal{U}$-suite exacte naturelle
$$
\hspace{8pt}
0\to\ker\hspace{1pt}(\hspace{1pt}\mathbb{F}_{2}\otimes_{\mathrm{P}}(M\otimes M)^{\mathfrak{S}_{2}}\to\Phi\mathcal{O}M\hspace{1pt})
\to\mathbb{F}_{2}\otimes_{\mathrm{P}}\mathfrak{S}_{2}M
\to\mathrm{R}_{1}\mathcal{O}M\to 0
\hspace{7pt}.
$$
\end{scho}

\begin{rem}\label{Phi-P-U} On peut définir un endofoncteur $\Phi:\mathrm{P}\text{-}\mathcal{U}\circlearrowleft$ de la même façon que l'on a défini l'endofoncteur $\Phi:\mathcal{U}\circlearrowleft$.

\smallskip
Soit $M$ un $\mathrm{P}$-$\mathrm{A}$-module instable~; on pose à nouveau $\Phi M:=\widehat{\mathrm{H}}^{0}(\mathfrak{S}_{2};M\otimes M)$. Le  $\mathrm{P}$-$\mathrm{A}$-module instable $\Phi M$ est $\mathrm{P}$-trivial, en formule $\Phi M=\theta\hspace{1pt}\Phi\mathcal{O}M$. Le $\mathcal{U}$-morphisme $\mathbb{F}_{2}\otimes_{\mathrm{P}}(M\otimes M)^{\mathfrak{S}_{2}}\to\Phi\mathcal{O}M$ est un avatar du  $\mathrm{P}\text{-}\mathcal{U}$-morphisme  (entre $\mathrm{P}$-$\mathrm{A}$-modules instables $\mathrm{P}$-triviaux) $\mathbb{F}_{2}\otimes_{\mathrm{P}}(M\otimes M)^{\mathfrak{S}_{2}}\to\Phi M$. Nous privilégierons cet avatar dans la suite de cette section car nous aurons besoin de la notation $\mathcal{O}$ pour un autre foncteur oubli.
\end{rem}

\bigskip
Compte tenu de \ref{R1nilf} et \ref{critnilf-2}, le scholie \ref{S1-bis} entraîne le suivant~:

\medskip
\begin{scho}\label{S1-ter} Soit $M$ un $\mathrm{P}$-$\mathrm{A}$-module instable. Si l'on suppose que  $M$ est $\mathcal{N}\hspace{-1,5pt}il$-fermé alors les deux propriétés suivantes sont équivalentes~:
\begin{itemize}
\item[(i)] Le $\mathrm{A}$-module instable $\mathbb{F}_{2}\otimes_{\mathrm{P}}\mathfrak{S}_{2}M$ est $\mathcal{N}\hspace{-1,5pt}il$-fermé.
\item[(ii)] Le $\mathrm{A}$-module instable $\ker\hspace{1pt}(\hspace{1pt}\mathbb{F}_{2}\otimes_{\mathrm{P}}(M\otimes M)^{\mathfrak{S}_{2}}\to\Phi M\hspace{1pt})$ est $\mathcal{N}\hspace{-1,5pt}il$-fermé.
\end{itemize}
\end{scho}

Démontrer la proposition \ref{quad-P-nilf-bts} revient donc à démontrer celle-ci~:

\bigskip
\begin{pro}\label{proclef-bis} Soit $M$ un $\mathrm{P}$-$\mathrm{A}$-module instable~; si $M$ est $\mathrm{P}$-$\mathcal{N}\hspace{-1,5pt}il$-fermé alors le $\mathrm{A}$-module instable
$$
\ker\hspace{2pt}(\hspace{2pt}\mathbb{F}_{2}\otimes_{\mathrm{P}}(M\otimes M)^{\mathfrak{S}_{2}}\to\Phi M\hspace{2pt})
$$
est $\mathcal{N}\hspace{-1,5pt}il$-fermé.
\end{pro}

\bigskip
\begin{rem}
La suite exacte $(\mathcal{S}_{1})$ et les propositions  \ref{R1nilf} et \ref{critnilf-0} montrent que si $M$ et  $\mathbb{F}_{2}\otimes_{\mathrm{P}}(M\otimes M)^{\mathfrak{S}_{2}}$ sont $\mathcal{N}\hspace{-1,5pt}il$-fermés alors $\mathbb{F}_{2}\otimes_{\mathrm{P}}\mathfrak{S}_{2}M$ l'est aussi (même argument que celui utilisé pour \ref{quadnilf}).  Hélas il est faux en général que $\mathbb{F}_{2}\otimes_{\mathrm{P}}(M\otimes M)^{\mathfrak{S}_{2}}$ soit $\mathcal{N}\hspace{-1,5pt}il$-fermé si $M$ est $\mathrm{P}$-$\mathcal{N}\hspace{-1,5pt}il$-fermé~:  voir ci-dessous. 
\end{rem}

\medskip
\begin{exple}\label{contrexples}
On prend $M=\mathrm{H}^{*}\mathrm{S}_{2}=\mathrm{P}$~; on constate que l'on a un $\mathcal{U}$-isomorphisme $\mathbb{F}_{2}\otimes_{\mathrm{P}}(\mathrm{P}\otimes\mathrm{P})^{\mathfrak{S}_{2}}\cong\Phi \mathrm{P}$ (observer que $(\mathrm{P}\otimes\mathrm{P})^{\mathfrak{S}_{2}}$ est l'algèbre de polynômes $\mathbb{F}_{2}[u\otimes 1+1\otimes u,u\otimes u]$). Le $\mathrm{A}$-module instable $\Phi\mathrm{P}$ est réduit mais n'est pas $\mathcal{N}\hspace{-1,5pt}il$-fermé, voir la proposition \ref{phinf} ci-après. Cependant on constate également  que l'isomorphisme $\mathbb{F}_{2}\otimes_{\mathrm{P}}(\mathrm{P}\otimes\mathrm{P})^{\mathfrak{S}_{2}}\cong\Phi \mathrm{P}$ est réalisé par la flèche de gauche de \ref{S1}.

\smallskip
En fait cet exemple correspond au premier pas de récurrence de la démonstra\-tion du théorème \ref{alt-3} dans le cas $n=2^{m}$ ~:  la propriété ``$\mathrm{H}^{*}(\mathrm{S}_{2};\mathbb{F}_{2})$ est $\mathrm{P}$-$\mathcal{N}\hspace{-1,5pt}il$-fermé'' implique la propriété ``$\mathrm{H}^{*}(\mathrm{S}_{4};\mathbb{F}_{2})$ est $\mathrm{P}$-$\mathcal{N}\hspace{-1,5pt}il$-fermé''.

\smallskip
Cet exemple permet aussi de tenir la promesse faite en \ref{helas}. Compte tenu de \ref{cokere} (et de la remarque \ref{A-Gysin-bis}), on sait, sans attendre la démonstration de \ref{quad-P-nilf}, que $\mathrm{H}^{*}\mathrm{S}_{4}\cong\mathfrak{S}_{2}\mathrm{H}^{*}\mathrm{S}_{2}=\mathfrak{S}_{2}\mathrm{P}$ est $\mathrm{P}$-$\mathcal{N}\hspace{-1,5pt}il$-fermé puisque $\mathrm{A}_{4}$ est isomorphe à $(\mathbb{Z}/2)^{2}$ (voir \ref{Klein})~; d'autre part on constate que l'on a $\mathfrak{S}_{2}\mathrm{P}/\tau\mathfrak{S}_{2}\mathrm{P}=(\mathrm{P}\otimes\mathrm{P})^{\mathfrak{S}_{2}}$ (pour une généralisation voir l'appendice \ref{appen-A}).
\end{exple}

\smallskip
\footnotesize
\begin{pro}\label{phinf} Soit $M$ un $\mathrm{A}$-module instable~; les deux conditions suivantes sont équivalentes~:
\begin{itemize}
\item[(i)] le $\mathrm{A}$-module instable $\Phi M$ est $\mathcal{N}\hspace{-1.5pt}il$-fermé~;
\item[(ii)] on a $M^{n}=0$ pour $n>0$ ($M$ est ``concentré'' en degré zéro).
\end{itemize}
\end{pro}

\textit{Démonstration de $(i)\Rightarrow(ii)$.} On note $\lambda_{M}$  l'application $\Phi M\to M, \Phi x\mapsto\mathrm{Sq}_{0}x$~; on constate que $\lambda_{M}$ est un $\mathcal{U}$-morphisme (les $\lambda_{M}$ fournissent une transformation naturelle de l'endofoncteur $\Phi:\mathcal{U}\circlearrowleft$ dans l'endofoncteur identité de $\mathcal{U}$). Si $\Phi M$ est $\mathcal{N}\hspace{-1,5pt}il$-fermé, il est \textit{a fortiori} réduit et il en est de même pour $M$~;  $\lambda_{M}$ est donc injective. On pose $Q:=\mathop{\mathrm{coker}}\lambda_{M}$. La proposition  \ref{critnilf-0} montre que $Q$ est réduit c'est-à-dire que $\lambda_{Q}$ est injective. D'autre part, par construction, l'application  $\lambda_{Q}$ est triviale. On  a donc $\Phi Q=0$, $Q=0$ et $M\cong\Phi M$. La~condition $(ii)$ en résulte aisément.
\hfill$\square$

\normalsize

\bigskip
\textsc{Démonstration de la proposition \ref{proclef-bis}}

\medskip
Cette proposition est au niveau technique le point-clef de cette section~; sa démonstration va nous occuper jusqu'à la fin de celle-ci.

\medskip
La proposition \ref{tens-P-nilf-1}, la proposition \ref{chdp} et le scholie \ref{invnf} montrent que le $\mathrm{A}$-module instable  $(\mathbb{F}_{2}\otimes_{\mathrm{P}}(M\otimes M))^{\mathfrak{S}_{2}}$ est $\mathcal{N}\hspace{-1,5pt}il$-fermé. Notre stratégie de démonstration est de comparer $\mathbb{F}_{2}\otimes_{\mathrm{P}}(M\otimes M)^{\mathfrak{S}_{2}}$ et $(\mathbb{F}_{2}\otimes_{\mathrm{P}}(M\otimes M))^{\mathfrak{S}_{2}}$ (attention au parenthésage~!)

\medskip
Le $(\mathrm{P}\text{-}\mathcal{U})$-morphisme $M\otimes M\to\mathrm{F}_{2}\otimes_{P}(M\otimes M)$ induit un $(\mathrm{P}\text{-}\mathcal{U})$-diagramme commutatif
$$
\begin{CD}
\mathrm{H}^{0}(\mathfrak{S}_{2};M\otimes M)
@>>>
\mathrm{H}^{0}(\mathfrak{S}_{2};\mathbb{F}_{2}\otimes_{\mathrm{P}}(M\otimes M)) \\
@VVV @VVV \\
\widehat{\mathrm{H}}^{0}(\mathfrak{S}_{2};M\otimes M)
@>>>
\widehat{\mathrm{H}}^{0}(\mathfrak{S}_{2};\mathbb{F}_{2}\otimes_{P}(M\otimes M))
\end{CD}
$$
dans lequel les $\mathrm{P}$-$\mathrm{A}$-modules instables de droite sont $\mathrm{P}$-triviaux, il induit donc  un diagrame commutatif de $\mathrm{P}$-$\mathrm{A}$-modules instables $\mathrm{P}$-triviaux
$$
\begin{CD}
\mathbb{F}_{2}\otimes_{\mathrm{P}}\mathrm{H}^{0}(\mathfrak{S}_{2};M\otimes M)
@>>>
\mathrm{H}^{0}(\mathfrak{S}_{2};\mathbb{F}_{2}\otimes_{\mathrm{P}}(M\otimes M)) \\
@VVV @VVV \\
\mathbb{F}_{2}\otimes_{\mathrm{P}}\widehat{\mathrm{H}}^{0}(\mathfrak{S}_{2};M\otimes M)
@>>>
\widehat{\mathrm{H}}^{0}(\mathfrak{S}_{2};\mathbb{F}_{2}\otimes_{P}(M\otimes M))
\end{CD}
$$
que l'on peut considérer comme un $\mathcal{U}$-diagramme commutatif. On allège un peu la notation et on nomme les flèches qui y apparaissent~:
$$
\hspace{5pt}
\begin{CD}
\mathbb{F}_{2}\otimes_{\mathrm{P}}(M\otimes M)^{\mathfrak{S}_{2}}
@>\kappa_{M}>>
(\mathbb{F}_{2}\otimes_{\mathrm{P}}(M\otimes M))^{\mathfrak{S}_{2}} \\
@V\nu_{M}VV @V\widetilde{\nu}_{M}VV \\
\Phi M
@>\widehat{\kappa}_{M}>>
\widehat{\Phi}M
\end{CD}
\hspace{23pt}.
\leqno{(\mathcal{D}_{M})}
$$

\smallskip
Archivons la définition de $\widehat{\Phi}M$~:

\medskip
\begin{defi}\label{phihat}
Soit $M$ un $\mathrm{P}$-$\mathrm{A}$-module instable~; le $\mathrm{A}$-module~instable $\widehat{\mathrm{H}}^{0}(\mathfrak{S}_{2};\mathbb{F}_{2}\otimes_{P}(M\otimes M))$ est noté $\widehat{\Phi}M$.
\end{defi}

\medskip
La commutativité de $(\mathcal{D}_{M})$ fait que
$\kappa_{M}$  induit un $\mathcal{U}$-morphisme naturel, disons $\iota_{M}$, de $\ker\nu_{M}$ sur $\ker\widetilde{\nu}_{M}$.

\begin{pro}\label{isoclef} Soit $M$ un $\mathrm{P}$-$\mathrm{A}$-module instable~; si $M$ est réduit\footnote{L'hypothèse que $M$ est réduit n'est pas nécessaire, cependant elle est satisfaite dans notre contexte (la démonstration de \ref{proclef-bis}) et elle simplifiera notre démonstration.} alors le $\mathcal{U}$-morphisme naturel
$$
\iota_{M}:\ker\hspace{1pt}(\hspace{1pt}\mathbb{F}_{2}\otimes_{\mathrm{P}}(M\otimes M)^{\mathfrak{S}_{2}}\overset{\nu_{M}}{\longrightarrow}\Phi M\hspace{1pt})
\to
\ker\hspace{1pt}(\hspace{1pt}\mathbb{F}_{2}\otimes_{\mathrm{P}}(M\otimes M))^{\mathfrak{S}_{2}}\overset{\widetilde{\nu}_{M}}{\longrightarrow}\widehat{\Phi} M\hspace{1pt})
$$
est un isomorphisme.
\end{pro}

\medskip
\textit{Démonstration.} On rappelle que l'on note $\mathcal{E}$ la catégorie des  $\mathbb{F}_{2}$-espaces vectoriels $\mathbb{N}$-gradués et $\mathrm{P}\text{-}\mathcal{E}$ la catégorie des $\mathbb{F}_{2}$-espaces vectoriels $\mathbb{N}$-gradués munis d'une structure de $\mathrm{P}$-module (au sens gradué). On note $\mathcal{O}:\mathrm{P}\text{-}\mathcal{U}\to\mathrm{P}\text{-}\mathcal{E}$ le~foncteur oubli évident. Il est clair qu'il suffit de montrer que le $\mathcal{E}$-morphisme sous-jacent à $\iota_{M}$ est un isomorphisme.

\medskip
On pose, comme dans la démonstration de \ref{tens-P-nilf-2}, $K=\tau M$ et $L:=M/K$~; le~point (c) de \ref{P-A-reduit-1} dit que $L$ est $\mathrm{P}$-libre ce qui implique que la suite exacte $0\to\mathcal{O}K\to\mathcal{O}M\to\mathcal{O}L\to 0$ est scindable. On peut donc supposer dans la démonstration de la $\mathcal{E}$-version de \ref{isoclef} que l'on a $\mathcal{O}M=\mathcal{O}K\oplus\mathcal{O}L$, ce que nous ferons ci-après. Dégageons au préalable un énoncé quasiment évident.

\medskip
Soit $\mathcal{C}$ l'une des deux catégories $\mathrm{P}\text{-}\mathcal{E}$ ou $\mathrm{P}\text{-}\mathcal{U}$. Soit $M$ un objet de $\mathcal{C}$. On pose $\Theta(M):=M\otimes M$ ;  $\Theta(M)$ est  un $\mathcal{C}$-objet muni d'une action de $\mathfrak{S}_{2}$. En clair  $\Theta(M)$ est muni du $\mathcal{C}$-automorphisme involutif défini par $\sigma(x\otimes y)=y\otimes x$.

\begin{pro}\label{cross-effect} Soient $M_{1}$ et $M_{2}$ deux $\mathcal{C}$-objets~; on a un isomorphisme canonique de $\mathcal{C}$-objets  munis d'une action de~$\mathfrak{S}_{2}$~:
$$
\hspace{24pt}
\Theta(M_{1}\oplus M_{2})
\hspace{4pt}\cong\hspace{4pt}
\Theta(M_{1})\oplus\Theta(M_{2})
\oplus((M_{1}\otimes M_{2})\oplus(M_{1}\otimes M_{2}))
\hspace{23pt},
$$
l'action de $\mathfrak{S}_{2}$ sur $(M_{1}\otimes M_{2})\oplus(M_{1}\otimes M_{2})$ étant donnée par l'échange des deux termes $M_{1}\otimes M_{2}$.
\end{pro}

\medskip
\footnotesize
\textit{Démonstration.} Le $\mathcal{C}$-objet $\Theta(M)$ est somme directe de trois sous-objets invariants sous l'action de $\mathfrak{S}_{2}$~:
$$
\hspace{24pt}
\Theta(M_{1}\oplus M_{2})
\hspace{4pt}=\hspace{4pt}
\Theta(M_{1})\oplus\Theta(M_{2})
\oplus((M_{1}\otimes M_{2})\oplus(M_{2}\otimes M_{1}))
\hspace{23pt}.
$$
Soit $\phi:M_{1}\otimes M_{2}\to M_{2}\otimes M_{1}$ le $\mathcal{C}$-isomorphisme défini par $\phi (x_{1}\otimes y_{2})=y_{2}\otimes x_{1}$.

\smallskip
La matrice du $\mathcal{C}$-automorphisme
$$
\hspace{24pt}
\sigma:(M_{1}\otimes M_{2})\oplus(M_{2}\otimes M_{1})\to
(M_{1}\otimes M_{2})\oplus(M_{2}\otimes M_{1})
\hspace{23pt},
$$
la matrice du $\mathcal{C}$-isomorphisme
$$
\mathrm{id}\oplus\phi:(M_{1}\otimes M_{2})\oplus(M_{1}\otimes M_{2})\to
(M_{1}\otimes M_{2})\oplus(M_{2}\otimes M_{1})
$$
et celle de son inverse
(matrices données par la décomposition en somme directe de la source et du but) sont respectivement
$$
\hspace{12pt}
\begin{bmatrix}
0 & \phi^{-1} \\
\phi & 0
\end{bmatrix}
\hspace{12pt},\hspace{12pt}
\begin{bmatrix}
1 & 0 \\
0 & \phi
\end{bmatrix}
\hspace{12pt}\text{et}\hspace{12pt}
\begin{bmatrix}
1 & 0 \\
0 & \phi^{-1}
\end{bmatrix}
\hspace{11pt}.
$$
On conclut en constatant que l'on a l'identité
$$
\hspace{24pt}
\begin{bmatrix}
1 & 0 \\
0 & \phi^{-1}
\end{bmatrix}
\begin{bmatrix}
0 & \phi^{-1} \\
\phi & 0
\end{bmatrix}
\begin{bmatrix}
1 & 0 \\
0 & \phi
\end{bmatrix}
\hspace{4pt}=\hspace{4pt}
\begin{bmatrix}
0 & 1 \\
1 & 0
\end{bmatrix}
\hspace{23pt}.
$$
\hfill$\square$
\normalsize

\medskip
\begin{cor}\label{add} Le foncteur $\widehat{\Phi}:M\mapsto\widehat{\mathrm{H}}^{0}(\mathfrak{S}_{2};\mathbb{F}_{2}\otimes_{P}(M\otimes M))$, défini sur $\mathrm{P}\text{-}\mathcal{E}$ ou $\mathrm{P}\text{-}\mathcal{U}$ et à valeurs dans $\mathcal{E}$ ou  $\mathcal{U}$, est additif.
\end{cor}

\medskip
\textit{Démonstration.} On a $\widehat{\mathrm{H}}^{0}(\mathfrak{S}_{2};\Gamma)=0$ si $\Gamma$ est un $\mathbb{F}_{2}[\mathfrak{S}_{2}]$-module libre.
\hfill$\square$

\medskip
On revient maintenant à la démonstration de \ref{isoclef}. Soit $M$ un $\mathcal{C}$-objet, on pose~:

\smallskip
--$\hspace{18pt}\mathrm{U}_{0}(M):=\mathbb{F}_{2}\otimes_{\mathrm{P}}(M\otimes M)^{\mathfrak{S}_{2}}\hspace{9pt},\hspace{9pt}\mathrm{U}_{1}(M):=(\mathbb{F}_{2}\otimes_{\mathrm{P}}(M\otimes M))^{\mathfrak{S}_{2}}\hspace{18pt};$

\smallskip
--$\hspace{18pt}\mathrm{V}_{0}(M):=\ker\hspace{1pt}(\mathrm{U}_{0}(M)\to\Phi M)
\hspace{9pt},\hspace{9pt}
\mathrm{V}_{1}(M):=\ker\hspace{1pt}(\mathrm{U}_{1}(M)\to\widehat\Phi M)\hspace{18pt}.$

\medskip
Soit $M$ un $\mathrm{P}$-$\mathrm{A}$-module instable réduit~; on choisit un $\mathrm{P}\text{-}\mathcal{E}$-isomorphisme $\mathcal{O}M\cong\mathcal{O}K\oplus\mathcal{O}L$. La proposition \ref{cross-effect} implique~:

\medskip
On a, pour  $i=0,1$, des $\mathcal{E}$-isomorphismes
$$
\hspace{24pt}
\mathrm{U}_{i}(\mathcal{O}M)\cong\mathrm{U}_{i}(\mathcal{O}K)\oplus\mathrm{U}_{i}(\mathcal{O}L)\oplus\mathbb{F}_{2}\otimes_{\mathrm{P}}(\mathcal{O}K\otimes\mathcal{O}L)
\hspace{24pt}
$$
et le $\mathcal{E}$-morphisme $\mathcal{O}\kappa_{M}:\mathrm{U}_{0}(\mathcal{O}M)\to\mathrm{U}_{1}(\mathcal{O}M)$  s'identifie à la somme directe de $\mathcal{O}\kappa_{K}$, $\mathcal{O}\kappa_{L}$ et l'identité.

\medskip
Le $\mathcal{E}$-morphisme $\mathcal{O}\nu_{M}:\mathrm{U}_{0}(\mathcal{O}M)\to\Phi\hspace{1pt}\mathcal{O}M$ (resp. $\mathcal{O}\widetilde{\nu}_{M}:\mathrm{U}_{1}(\mathcal{O}M)\to\widehat{\Phi}\hspace{1pt}\mathcal{O}M$) est nul sur $\mathbb{F}_{2}\otimes_{\mathrm{P}}(\mathcal{O}K\otimes\mathcal{O}L)$ et sa restriction à $\mathrm{U}_{0}(\mathcal{O}K)\oplus\mathrm{U}_{0}(\mathcal{O}L)$ (resp. $\mathrm{U}_{1}(\mathcal{O}K)\oplus\mathrm{U}_{1}(\mathcal{O}L)$) s'identifie à la somme directe de $\mathcal{O}\nu_{K}$ et $\mathcal{O}\nu_{L}$ (resp.~$\mathcal{O}\widetilde{\nu}_{K}$ et $\mathcal{O}\widetilde{\nu}_{L}$).

\medskip
Ceci montre que le $\mathcal{E}$-morphisme 
$\mathcal{O}\iota_{M}:V_{0}(\mathcal{O}M)\to V_{1}(\mathcal{O}M)$ s'identifie à somme directe de $\mathcal{O}\iota_{K}$ et $\mathcal{O}\iota_{L}$. Il est évident que $\iota_{K}$ est un isomorphisme car $K$ est par définition un $\mathrm{P}$-$\mathrm{A}$-module instable $\mathrm{P}$-trivial~; la démonstration de \ref{isoclef} sera donc achevée une fois démontrée la proposition suivante~:

\begin{pro}\label{isoclef-lib} Soit $L$ un $\mathrm{P}$-$\mathrm{A}$-module instable $\mathrm{P}$-libre~; le $\mathcal{U}$-morphisme naturel
$$
\iota_{L}:\ker\hspace{1pt}(\hspace{1pt}\mathbb{F}_{2}\otimes_{\mathrm{P}}(L\otimes L)^{\mathfrak{S}_{2}}\overset{\nu_{L}}{\longrightarrow}\Phi L\hspace{1pt})
\to
\ker\hspace{1pt}(\hspace{1pt}\mathbb{F}_{2}\otimes_{\mathrm{P}}(L\otimes L))^{\mathfrak{S}_{2}}\overset{\widetilde{\nu}_{L}}{\longrightarrow}\widehat{\Phi} L\hspace{1pt})
$$
est un isomorphisme.
\end{pro}

\medskip
\textit{Démonstration.} On introduit deux complexes de chaînes dans la catégorie $\mathrm{P}\text{-}\mathcal{U}$, $C_{\bullet}$ et $D_{\bullet}$,  et un homomorphisme $\gamma:C_{\bullet}\to D_{\bullet}$, définis comme suit~:

\medskip
-- On pose $C_{k}:=L\otimes L$ pour $k=0,1,2$ et $C_{k}:=0$ sinon. L'opérateur de bord  $\mathrm{d}_{k}:C_{k+1}\to C_{k}$ est $1+\sigma$ pour $k=0,1$ (et $0$ sinon).

\medskip
-- On pose $D_{k}:=L\otimes L$ pour $k=0,1,2,3$  et $D_{k}:=0$ sinon. L'opérateur de bord  $\mathrm{d}_{k}:D_{k+1}\to D_{k}$ est $1+\sigma$ pour $k=0,1,2$ (et $0$ sinon).

\medskip
-- L'homomorphisme $\gamma_{k}:C_{k}\to D_{k}$ est l'identité de $L\otimes L$ pour $k=0,1,2$ (et~$0$ sinon).

\medskip
Le $\mathrm{P}$-$\mathrm{A}$-module instable $L\otimes L$, est libre comme $\mathrm{P}$-module, par exemple parce que $\mathrm{P}\otimes\mathrm{P}$, vu comme un $\mathrm{P}$-module \textit{via} la diagonale $\psi:\mathrm{P}\to\mathrm{P}\otimes\mathrm{P}$, est libre. Comme l'anneau $\mathrm{P}$ est principal, on peut appliquer le théorème des coef-ficients universels (dans sa version fonctorielle). On obtient un $\mathcal{U}$-diagramme commutatif

\footnotesize
$$
\leqno{(\mathcal{D})}\label{D2}
\begin{CD}
0 @>>>
\mathbb{F}_{2}\otimes_{\mathrm{P}}\mathrm{H}_{2}\hspace{1pt}C_{\bullet}
@>>>
\mathrm{H}_{2}(\mathbb{F}_{2}\otimes_{\mathrm{P}}C_{\bullet})
@>>>
\mathrm{Tor}_{1}^{\mathrm{P}}(\mathbb{F}_{2},\mathrm{H}_{1}\hspace{1pt}C_{\bullet})
@>>>0 \\
&& @V\gamma VV @V\gamma VV@ V\gamma VV \\
0 @>>>
\mathbb{F}_{2}\otimes_{\mathrm{P}}\mathrm{H}_{2}\hspace{1pt}D_{\bullet}
@>>>
\mathrm{H}_{2}(\mathbb{F}_{2}\otimes_{\mathrm{P}}D_{\bullet})
@>>>
\mathrm{Tor}_{1}^{\mathrm{P}}(\mathbb{F}_{2},\mathrm{H}_{1}\hspace{1pt}D_{\bullet})
@>>>0
\end{CD}
$$
\normalsize
dont les lignes sont exactes.

\smallskip
\footnotesize
Pour le confort du lecteur, rappelons comment développer  \textit{ab initio} la théorie du théorème des coefficients universels dans le cas qui nous occupe.

\smallskip
Soit $\Gamma_{\bullet}$ un $(\mathrm{P}\text{-}\mathcal{U})$-complexe de chaînes tel que $\Gamma_{n}$ est $\mathrm{P}$-libre pour tout $n$. On~considère à  nouveau la $(\mathrm{P}\text{-}\mathcal{U})$-suite exacte $0\to\widetilde{\mathrm{P}}\to\mathrm{P}\to
\mathbb{F}_{2}\to 0$ qui apparaît dans la démonstration de \ref{sigmatau}, suite exacte qui fournit une résolution de $\mathbb{F}_{2}$ dans la catégorie $\mathrm{P}\text{-}\mathcal{U}$ qui est une résolution libre dans la catégorie $\mathrm{P}\text{-}\mathcal{E}$. Comme les $\Gamma_{n}$ sont $\mathrm{P}$-libres cette suite exacte induit une suite exacte courte de $(\mathrm{P}\text{-}\mathcal{U})$-complexes
$$
\hspace{24pt}
0\to\widetilde{\mathrm{P}}\otimes_{\mathrm{P}}\Gamma_{\bullet}\to
\Gamma_{\bullet}
\to\mathbb{F}_{2}\otimes_{\mathrm{P}}\Gamma_{\bullet}
\to 0
\hspace{23pt};
$$
la longue suite exacte d'homologie associée donne pour tout $n$ une $(\mathrm{P}\text{-}\mathcal{U})$-suite exacte
$$
\widetilde{\mathrm{P}}\otimes_{\mathrm{P}}\mathrm{H}_{n}\Gamma_{\bullet}
\to\mathrm{H}_{n}\Gamma_{\bullet}
\to\mathrm{H}_{n}(\mathbb{F}_{2}\otimes_{\mathrm{P}}\Gamma_{\bullet})
\overset{\partial}{\to}\widetilde{\mathrm{P}}\otimes_{\mathrm{P}}\mathrm{H}_{n-1}\Gamma_{\bullet}
\to\mathrm{H}_{n-1}\Gamma_{\bullet}
$$
telle que le conoyau de la première flèche et le noyau de la dernière s'identifient respectivement à $\mathbb{F}_{2}\otimes_{\mathrm{P}}\mathrm{H}_{n}\Gamma_{\bullet}$ et 
$\mathrm{Tor}^{\mathrm{P}}_{1}(\mathbb{F}_{2},\mathrm{H}_{n-1}\Gamma_{\bullet})$ (invoquer \ref{res-P-lib}).  D'où le théorème des coefficients universels~; la fonctorialité est évidente.
\normalsize

\medskip
On explicite maintenant le diagramme $(\mathcal{D})$. On observe que l'on a tout fait pour avoir $\mathrm{H}_{2}(\mathbb{F}_{2}\otimes_{\mathrm{P}}C_{\bullet})=(\mathbb{F}_{2}\otimes_{\mathrm{P}}(L\otimes L))^{\mathfrak{S}_{2}}$ et $\mathrm{H}_{2}(\mathbb{F}_{2}\otimes_{\mathrm{P}}D_{\bullet})=\widehat{\Phi}L$, puis on fait les constatations suivantes~:

\medskip
--\hspace{4pt}$\mathrm{H}_{2}\hspace{1pt}C_{\bullet}=\mathrm{H}^{0}(\mathfrak{S}_{2}; L\otimes L)=(L\otimes L)^{\mathfrak{S}_{2}}$~;

\medskip
--\hspace{4pt}$\mathrm{H}_{1}\hspace{1pt}C_{\bullet}=\mathrm{H}_{1}\hspace{1pt}D_{\bullet}=\mathrm{H}_{2}\hspace{1pt}D_{\bullet}=\widehat{\mathrm{H}}^{0}(\mathfrak{S}_{2}; L\otimes L)\cong\Phi L$ (voir \ref{Phi-P-U})~;

\medskip
--\hspace{4pt}$\mathrm{Tor}_{1}^{\mathrm{P}}(\mathbb{F}_{2},\mathrm{H}_{1}\hspace{1pt}C_{\bullet})=\mathrm{Tor}_{1}^{\mathrm{P}}(\mathbb{F}_{2}, \mathrm{H}_{1}\hspace{1pt}D_{\bullet})=\Sigma\hspace{1pt}\Phi L$ (se rappeler que le $\mathrm{P}$-$\mathrm{A}$-module instable $\Phi L$ est $\mathrm{P}$-trivial et invoquer \ref{sigmatau}).

\medskip
On en déduit  facilement la version explicite du diagramme $(\mathcal{D})$~:

\footnotesize
$$
\begin{CD}
0@>>>
\mathbb{F}_{2}\otimes_{\mathrm{P}}(L\otimes L)^{\mathfrak{S}_{2}}
@>\kappa_{L}>>
(\mathbb{F}_{2}\otimes_{\mathrm{P}}(L\otimes L))^{\mathfrak{S}_{2}}
@>>>\Sigma\hspace{1pt}\Phi L
@>>> 0
\\
& & @V\nu_{L}VV @V\widetilde{\nu}_{L}VV @V\mathrm{id}VV \\
0@>>> \Phi L @>\widehat{\kappa}_{L}>> \widehat{\Phi}L
@>>>\Sigma\hspace{1pt}\Phi L
@>>> 0
\end{CD}
\leqno{(\mathcal{D}^{\hspace{1pt}+}_{L})}
$$
\normalsize
(le carré commutatif constitué des flèches $\kappa_{L}$, $\widehat{\kappa}_{L}$, $\nu_{L}$ et $\widetilde{\nu}_{L}$ est le diagramme $(\mathcal{D}_{L})$ ce qui explique la notation $(\mathcal{D}^{\hspace{1pt}+}_{L})$).

\medskip
L'exactitude des deux suites horizontales implique que $\kappa_{L}$ induit un isomorphisme entre les noyaux des deux flèches verticales de gauche  ce qui achève la démonstration de \ref{isoclef-lib} et donc de \ref{isoclef}.
\hfill$\square\square$

\medskip
Compte tenu de \ref{isoclef}, démontrer la proposition \ref{proclef-bis} revient à démontrer la suivante~:

\begin{pro}\label{proclef-ter} Soit $M$ un $\mathrm{P}$-$\mathrm{A}$-module instable~; si $M$ est $\mathrm{P}$-$\mathcal{N}\hspace{-1,5pt}il$-fermé alors le $\mathrm{A}$-module instable
$$
\ker\hspace{2pt}(\hspace{2pt}(\mathbb{F}_{2}\otimes_{\mathrm{P}}(M\otimes M))^{\mathfrak{S}_{2}}\to\widehat{\Phi} M\hspace{2pt})
$$
est $\mathcal{N}\hspace{-1,5pt}il$-fermé.
\end{pro}

\medskip
\textit{Démonstration.} Comme $(\mathbb{F}_{2}\otimes_{\mathrm{P}}(M\otimes M))^{\mathfrak{S}_{2}}$ est $\mathcal{N}\hspace{-1,5pt}il$-fermé il faut montrer que $\widehat{\Phi} M$ est réduit (Proposition \ref{critnilf-0}). Or on a~:

\begin{lem}\label{phihat-reduit} Soit $M$ un $\mathrm{P}$-$\mathrm{A}$-module instable~; si $M$ est réduit alors il en est de même pour $\widehat{\Phi} M$.
\end{lem}

\medskip
\textit{Démonstration.} On reprend les notations de la  démonstration de \ref{isoclef}~; on note en outre $i:K\to M$ et $p:M\to L$ l'inclusion et la surjection canoniques (morphismes dans la catégorie $\mathrm{P}\text{-}\mathcal{U}$). Dans la démonstration de \ref{isoclef} on a noté $\mathcal{O}$ le foncteur oubli $\mathrm{P}\text{-}\mathcal{U}\to\mathrm{P}\text{-}\mathcal{E}$ le foncteur oubli ; on note $\mathcal{O}'$ le foncteur oubli $\mathcal{U}\to\mathcal{E}$. Par définition on a $\mathcal{O}'\widehat{\Phi}=\widehat{\Phi}\mathcal{O}$, $\widehat{\Phi}$ désignant à gauche le foncteur $\widehat{\Phi}:\mathrm{P}\text{-}\mathcal{U}\to\mathcal{U}$ et à droite  le foncteur $\widehat{\Phi}:\mathrm{P}\text{-}\mathcal{E}\to\mathcal{E}$.

\begin{lem}\label{scindable} Soit $M$ un $\mathrm{P}$-$\mathrm{A}$-module instable~; si $M$ est réduit alors la suite de $\mathrm{A}$-modules instables
$$
\begin{CD}
0@>>>\widehat{\Phi}K@>\widehat{\Phi}i>>\widehat{\Phi}M@>\widehat{\Phi}p>>
\widehat{\Phi}L@>>>0
\end{CD}
$$
est exacte.
\end{lem}

\medskip
\textit{Démonstration.} Soit $(S)$ la $(\mathrm{P}\text{-}\mathcal{U})$-suite
$$
\begin{CD}
0@>>>K@>i>>M@>p>>L@>>>0
\end{CD}
$$
(qui est exacte~!). Il est clair que la $\mathcal{U}$-suite $\widehat{\Phi}\hspace{1pt}(S)$ est exacte si et seulement la $\mathcal{E}$-suite $\mathcal{O}'\widehat{\Phi}\hspace{1pt}(S)=\widehat{\Phi}\mathcal{O}\hspace{1pt}(S)$ est exacte. La $\mathcal{E}$-suite $\widehat{\Phi}\mathcal{O}\hspace{1pt}(S)$ est exacte parce que le foncteur $\widehat{\Phi}:\mathrm{P}\text{-}\mathcal{E}\to\mathcal{E}$ est additif et que
la  $(\mathrm{P}\text{-}\mathcal{E})$-suite $\mathcal{O}\hspace{1pt}(S)$ est une suite exacte scindable (voir par exemple  \cite[4.6, Theorem 1]{Par}).
\hfill$\square$

\medskip
Comme $K$ est $\mathrm{P}$-trivial, $\widehat{\Phi}K$ s'identifie à $\Phi K$ et est donc réduit~; du coup pour montrer que $\widehat{\Phi}M$ est réduit il suffit de montrer que $\widehat{\Phi}L$ est réduit.

\begin{lem}\label{phihat-reduit-bis} Soit $M$ un $\mathrm{P}$-$\mathrm{A}$-module instable~; si $M$ est réduit alors il en est de même pour $\widehat{\Phi} L$ (avec $L:=M/\tau M$).
\end{lem}

\medskip
\textit{Démonstration.} Il découle de \ref{chdp-1} que le $\mathrm{A}$-module instable $\widehat{\Phi} L$  est sous-jacent au $\mathrm{P}$-$\mathrm{A}$-module instable $\widehat{\mathrm{H}}^{0}(\mathfrak{S}_{2};L\otimes_{\mathrm{P}}L)$, $\mathfrak{S}_{2}$ agissant sur $L\otimes_{\mathrm{P}}L$ par ``échange des facteurs''~; ce point de vue nous sera commode ci-après.

\smallskip
La proposition \ref{P-A-reduit-P-lib} dit qu'il existe un  $(\mathrm{P}$-$\mathcal{U})$-monomorphisme $e:L\to\mathrm{P}\otimes I$ avec $I$ un $\mathcal{U}$-injectif réduit. En utilisant cette information et la ``fonctorialité en $L$'' de la $\mathcal{U}$-suite exacte
$$
\begin{CD}
0@>>> \Phi L @>\widehat{\kappa}_{L}>> \widehat{\Phi}L
@>\partial_{L}>>\Sigma\hspace{1pt}\Phi L
@>>> 0
\end{CD}
\leqno{(\mathcal{S}_{L})}
$$
(la seconde ligne du diagramme ($\mathcal{D}_{L}^{+})$, l'introduction de la notation $\partial_{L}$ est transparente), on obtient un $\mathcal{U}$-diagramme commutatif
$$
\begin{CD}
\hspace{12pt}
0\hspace{-5pt}@>>> \Phi L @>\widehat{\kappa}_{L}>> \widehat{\Phi}L
@>\partial_{L}>>\Sigma\hspace{1pt}\Phi L
@>>>\hspace{-12pt}0
\hspace{12pt}
\\
& & @V\Phi eVV @V\widehat{\Phi}eVV @V\Sigma\Phi eVV
\\
\hspace{20pt}
0@>>> \Phi(\mathrm{P}\otimes I)@>\widehat{\kappa}_{\hspace{1pt}\mathrm{P}\otimes I}>>\widehat{\Phi}(\mathrm{P}\otimes I)
@>\partial_{\hspace{1pt}\mathrm{P}\otimes I}>>\Sigma\hspace{1pt}\Phi(\mathrm{P}\otimes I)
@>>> 0 
\hspace{19pt}.
\end{CD}
$$
Puisque $\Phi e$ et $\Sigma\Phi e$ sont injectifs il en est de même pour $\widehat{\Phi}e$. On constate que $\widehat{\Phi}(\mathrm{P}\otimes I)$ est isomorphe à $\mathrm{P}\otimes\Phi I$~; puisque $\mathrm{P}\otimes\Phi I$ est réduit il en est de même pour $\widehat{\Phi}L$.
\hfill$\square\square$

\bigskip
\textsc{Complément}

\medskip
On démystifie ci-dessous le $\mathrm{P}$-$\mathrm{A}$-module instable $\widehat{\Phi}L$ (et la suite exacte $(\mathcal{S}_{L})$).

\medskip
L'application $\Phi\mathrm{P}\to\mathrm{P},x\mapsto x^{2}$ est un $\mathcal{K}$-morphisme (on peut remplacer dans cette assertion $\mathrm{P}$ par une $\mathrm{A}$-algèbre instable arbitraire)~; ce $\mathcal{K}$-morphisme est injectif, il permet d'identifier $\Phi\mathrm{P}$ avec une sous-$\mathrm{A}$-algèbre instable de $\mathrm{P}$. Un~$\mathrm{P}$-$\mathrm{A}$-module instable est donc naturellement un $\Phi\mathrm{P}$-$\mathrm{A}$-module instable et l'endofoncteur $\Phi:\mathcal{U}\circlearrowleft$ induit un foncteur $\mathrm{P}\text{-}\mathcal{U}\to\Phi\mathrm{P}\text{-}\mathcal{U}$ que l'on notera encore~$\Phi$.

\footnotesize
\begin{lem}\label{phiPlin-1} Le $\mathcal{U}$-morphisme $\widehat{\kappa}_{L}:\Phi L \to\widehat{\Phi}L$ est un $(\Phi\mathrm{P}\text{-}\mathcal{U})$-morphisme
\end{lem}

\smallskip
\textit{Démonstration.}  Par définition $\widehat{\kappa}_{L}$ est le monomorphisme $\mathbb{F}_{2}\otimes_{\mathrm{P}}\mathrm{H}_{2}\hspace{1pt}D_{\bullet}\to\mathrm{H}_{2}(\mathbb{F}_{2}\otimes_{\mathrm{P}}D_{\bullet})$~; comme l'on a $\mathrm{H}_{2}\hspace{1pt}D_{\bullet}=\Phi L$ et que le $\mathrm{P}$-$\mathrm{A}$-module instable $\Phi L$ est forcément $\mathrm{P}$-trivial on a $\mathbb{F}_{2}\otimes_{\mathrm{P}}\mathrm{H}_{2}\hspace{1pt}D_{\bullet}\cong\mathrm{H}_{2}\hspace{1pt}D_{\bullet}$, si bien que $\widehat{\kappa}_{L}$ s'identifie aussi au morphisme $\mathrm{H}_{2}\hspace{1pt}D_{\bullet}\to\mathrm{H}_{2}(\mathbb{F}_{2}\otimes_{\mathrm{P}}D_{\bullet})$. Rappelons que le complexe $\mathbb{F}_{2}\otimes_{\mathrm{P}}D_{\bullet}$ peut être vu comme le $(\mathrm{P}\text{-}\mathcal{U})$-complexe
$$
\begin{CD}
L\otimes_{\mathrm{P}}L@<1+\sigma<< L\otimes_{\mathrm{P}}L@<1+\sigma<< L\otimes_{\mathrm{P}}L @<1+\sigma<< L\otimes_{\mathrm{P}}L\leftarrow 0\leftarrow 0\ldots
\end{CD}
$$
(compte tenu de \ref{chdp-1}). Soit $x$ un élément de $L$, $\Phi x$ est représenté par le cycle $x\otimes x$ de $\mathrm{Z}_{2}\hspace{1pt}\mathrm{D}_{\bullet}$ dont l'image dans $\mathrm{Z}_{2}(\mathbb{F}_{2}\otimes_{P}\mathrm{D}_{\bullet})$ est le cycle $x\otimes_{\mathrm{P}}x$. Les  égalités $\Phi u\hspace{1pt}\Phi x=\Phi(ux)$ et $ux\otimes_{\mathrm{P}}ux=u^{2}(x\otimes_{\mathrm{P}}x)$ permettent de conclure
\hfill$\square$

\begin{lem}\label{phiPlin-2} Le $\mathcal{U}$-morphisme $\partial_{L}:\widehat{\Phi}L\to\Sigma\Phi L$ est un $(\Phi\mathrm{P}\text{-}\mathcal{U})$-morphisme
\end{lem}

\smallskip
\textit{Démonstration.} Par définition, le $\mathcal{U}$-morphisme $\partial_{L}:\widehat{\Phi}L\to\Sigma\Phi L$ est induit par le connectant $\partial:\mathrm{H}_{2}(\mathbb{F}_{2}\otimes_{\mathrm{P}}D_{\bullet})\to\widetilde{\mathrm{P}}\otimes_{\mathrm{P}}\mathrm{H}_{1}D_{\bullet}$ associée à la  suite exacte courte de $(\mathrm{P}\text{-}\mathcal{U})$-complexes
$0\to\widetilde{\mathrm{P}}\otimes_{\mathrm{P}}D_{\bullet}\to
D_{\bullet}\to\mathbb{F}_{2}\otimes_{\mathrm{P}}D_{\bullet}\to 0$.

\smallskip
Soit $B$ une base (homogène) de $L$ en tant que $\mathrm{P}$-$\mathbb{F}_{2}$-espace vectoriel $\mathbb{N}$-gradué~; $B\times B$ est une base du $(\mathrm{P}\otimes\mathrm{P})$-$\mathbb{F}_{2}$-ev\footnote{ev est une abréviation pour espace vectoriel} $\mathbb{N}$-gradué $L\otimes L$ et $\{x\otimes_{P}y;\{x,y)\in B\times B\}$ est une base du $\mathrm{P}$-$\mathbb{F}_{2}$-ev $\mathbb{N}$-gradué $L\otimes_{\mathrm{P}}L=\mathbb{F}_{2}\otimes_{\mathrm{P}}D_{\bullet}$. On fait les constatations suivantes~:

\smallskip
-- $\{x\otimes_{P}x;x\in B\}\coprod\{x\otimes_{P}y:\{x,y\}\in\mathcal{P}_{2}(B)\}$ ($\mathcal{P}_{2}(B)$ désignant l'ensemble des parties à deux éléments de $B$) est une base du $\mathrm{P}$-$\mathbb{F}_{2}$-ev $\mathbb{N}$-gradué $\mathrm{Z}_{2}(\mathbb{F}_{2}\otimes_{P}\mathrm{D}_{\bullet})$~;

\smallskip
-- $\{x\otimes_{P}y:\{x,y\}\in\mathcal{P}_{2}(B)\}$ est une base du $\mathrm{P}$-$\mathbb{F}_{2}$-ev $\mathbb{N}$-gradué $\mathrm{B}_{2}(\mathbb{F}_{2}\otimes_{P}\mathrm{D}_{\bullet})$.

\smallskip
Il en résulte que les classes des cycles $x\otimes_{P}x$, $x$ parcourant $B$, constituent une base du $\mathrm{P}$-$\mathbb{F}_{2}$-ev $\mathbb{N}$-gradué  $\mathrm{H}_{2}(\mathbb{F}_{2}\otimes_{P}\mathrm{D}_{\bullet})=:\widehat{\Phi}L$.

\smallskip
Soit $x$ un élément de $B$~; on note $[x\otimes_{P}x]$ dans $\mathrm{H}_{2}(\mathbb{F}_{2}\otimes_{P}\mathrm{D}_{\bullet})$. du cycle $x\otimes_{P}x$. Soit~$n$ un entier naturel~; si $n$ est pair on a $\partial_{L}(u^{n}[x\otimes_{P}x])=0$ pour des raisons de degré, on vérifie ci-après, par un \textit{diagram chasing} (fastidieux~!), que l'on a $\partial_{L}(u^{2m+1}[x\otimes_{P}x])=\Sigma\Phi(u^{m}x)$, égalité dont le lemme découle~:

\smallskip
-- Le cycle $u^{2m+1}(x\otimes_{\mathrm{P}}x)=u^{m+1}x\otimes_{\mathrm{P}}u^{m}x$ est l'image par le morphisme $L\otimes L=:D_{2}\to\mathbb{F}_{2}\otimes_{\mathrm{P}}D_{2}:=L\otimes_{\mathrm{P}}L$ de l'élément $u^{m+1}x\otimes u^{m}x$~;

\smallskip
-- l'image par l'opérateur de bord $D_{2}\to D_{1}:=L\otimes L$ de $u^{m+1}x\otimes u^{m}x$ est  le cycle $u^{m+1}x\otimes u^{m}x+u^{m}x\otimes u^{m+1}x$ $=u.(u^{m}x\otimes u^{m}x)$ (par définition de la structure de $\mathrm{P}$-$\mathrm{A}$-module instable de $L\otimes L$)~;

\smallskip
-- $u.(u^{m}x\otimes u^{m}x)$ est l'image  par le morphisme $\widetilde{\mathrm{P}}\otimes_{\mathrm{P}}D_{1}\to D_{1}$ du cycle $u\otimes_{\mathrm{P}}(u^{m}x\otimes u^{m}x)$~;

\smallskip
-- la classe du cycle $u\otimes_{\mathrm{P}}(u^{m}x\otimes u^{m}x)$ dans $\mathrm{H}_{1}(\widetilde{\mathrm{P}}\otimes_{\mathrm{P}}D_{\bullet})=\widetilde{\mathrm{P}}\otimes_{\mathrm{P}}\mathrm{H}_{1}D_{\bullet}=\widetilde{\mathrm{P}}\otimes_{\mathrm{P}}\Phi L$ est $u\otimes_{\mathrm{P}}\Phi(u^{m}x)$.
\normalsize

\medskip
On note $\mathrm{e}\text{-}\widehat{\kappa}_{L}:\mathrm{P}\otimes_{\Phi\mathrm{P}}\Phi L\to\widehat{\Phi}L$ le $(\mathrm{P}\text{-}\mathcal{U})$-morphisme obtenu par extension des scalaires à partir du $(\Phi\mathrm{P}\text{-}\mathcal{U})$-morphisme $\widehat{\kappa}_{L}$.

\begin{pro}\label{phihat-final} Le $(\mathrm{P}\text{-}\mathcal{U})$-morphisme $\mathrm{e}\text{-}\widehat{\kappa}_{L}:\mathrm{P}\otimes_{\Phi\mathrm{P}}\Phi L\to\widehat{\Phi}L$ est un isomorphisme.
\end{pro}

\medskip
\textit{Démonstration.} On note $\mathrm{d}: \mathrm{P}\to\Sigma\mathrm{P}$ le $\mathcal{U}$-morphisme composé de la diagonale $\psi:\mathrm{P}\to\mathrm{P}\otimes\mathrm{P}$ et du produit tensoriel $\mathrm{id}\otimes\mathrm{j}:\mathrm{P}\otimes\mathrm{P}\to\mathrm{P}\otimes\Sigma\mathbb{F}_{2}=\Sigma\mathrm{P}$, $\mathrm{j}$~désignant l'unique application non triviale de $\mathrm{P}$ dans $\Sigma\mathbb{F}_{2}$ \footnote{L'application $\mathrm{d}$ est un avatar de la différentielle de Kähler $\mathbb{F}_{2}[u]\to\Omega_{\mathbb{F}_{2}[u]\slash\mathbb{F}_{2}}$.} On constate que $\mathrm{d}$ est un $(\Phi\mathrm{P}\text{-}\mathcal{U})$-morphisme\footnote{L'endofoncteur $\Sigma:\mathcal{U}\circlearrowleft$ induit un endofoncteur $\Sigma:\Phi\mathrm{P}\text{-}\mathcal{U}\circlearrowleft$.} et que le noyau et  l'image de $\mathrm{d}$ sont respectivement $\Phi\mathrm{P}$ et $\Sigma\Phi\mathrm{P}$~; on dispose donc d'une $(\Phi\mathrm{P}\text{-}\mathcal{U})$-suite exacte ``canonique''~:
$$
\hspace{48pt}
0\to\Phi\mathrm{P}\to\mathrm{P}\to\Sigma\Phi\mathrm{P}\to 0
\hspace{47pt}.
\leqno{(\mathcal{T})}
$$
Comme $\Phi L$ est un $\Phi\mathrm{P}$-module libre la $(\Phi\mathrm{P}\text{-}\mathcal{U})$-suite $(\mathcal{T})\otimes_{\Phi\mathrm{P}}\Phi L$, à savoir
$$
\hspace{24pt}
0\to\Phi L\to\mathrm{P}\otimes_{\Phi\mathrm{P}}\Phi L\to\Sigma\Phi L\to 0
\hspace{23pt},
$$
est encore exacte. On considère le $(\Phi\mathrm{P}\text{-}\mathcal{U})$-diagramme
$$
\hspace{24pt}
\begin{CD}
0@>>>\Phi L@>>>\mathrm{P}\otimes_{\Phi\mathrm{P}}\Phi L@>>>\Sigma\Phi L@>>>0
\\
&& @V\mathrm{id}VV @V\mathrm{e}\text{-}\widehat{\kappa}_{L}VV @V\mathrm{id}VV
\\
0@>>> \Phi L @>\widehat{\kappa}_{L}>> \widehat{\Phi}L
@>\partial_{L}>>\Sigma\hspace{1pt}\Phi L
@>>> 0
\end{CD}
\hspace{23pt}.
$$
Le carré de gauche est commutatif par définition même de $\mathrm{e}\text{-}\widehat{\kappa}_{L}$. Les calculs effectués lors des démonstrations de \ref{phiPlin-1} et \ref{phiPlin-2} montrent que le carré de droite l'est aussi. Précisons un peu.

\footnotesize

\smallskip
Soit $x$ un élément de $B$ (base homogène du $\mathbb{F}_{2}$-ev gradué $L$). On a montré lors de la démonstration de \ref{phiPlin-1} l'égalité $\widehat{\kappa}_{L}(\Phi x)=[x\otimes_{\mathrm{P}}x]$, d'où, puisque $\mathrm{e}\text{-}\widehat{\kappa}_{L}$ est $\mathrm{P}$-linéaire, l'égalité $\mathrm{e}\text{-}\widehat{\kappa}_{L}(u^{2m+1}\otimes_{\Phi\mathrm{P}}\Phi x)=u^{2m+1}[x\otimes_{\mathrm{P}}x]$ et, compte tenu de ce que l'on vu lors de la démonstration de \ref{phiPlin-2}, l'égalité $(\partial_{L}\circ\mathrm{e}\text{-}\widehat{\kappa}_{L})(u^{2m+1}\otimes_{\Phi\mathrm{P}}\Phi x)=\Sigma\Phi(u^{m}x)$. On conclut en observant que l'égalité $\mathrm{d}(u^{2m+1})=\Sigma u^{2m}$ implique que l'image de $u^{2m+1}\otimes_{\Phi\mathrm{P}}\Phi x)$ par le morphisme $\mathrm{P}\otimes_{\Phi\mathrm{P}}\Phi L\to\Sigma\Phi L$ est aussi $\Sigma\Phi(u^{m}x)$.

\normalsize

\smallskip
La commutativité du diagramme ci-dessus, dont les lignes sont exactes, entraîne~que $\mathrm{e}\text{-}\widehat{\kappa}_{L}$ est un isomorphisme.
\hfill$\square$

\sect{Sur la cohomologie modulo $2$ des groupes de Coxeter finis}

\begin{theo}\label{Coxeter-1} L'application de Quillen pour la cohomologie modulo $2$ d'un groupe de Coxeter~fini est un isomorphime.
\end{theo}

\medskip
L'énoncé ci-desus est conséquence de l'énoncé suivant, d'apparence technique mais en fait plus fort (voir \ref{Sylow-2})~:

\begin{theo}\label{Coxeter-2} La cohomologie modulo $2$ d'un $2$-Sylow d'un groupe de Coxeter fini est $\mathcal{N}\hspace{-1,5pt}il$-fermée.
\end{theo}

\medskip
\textit{Démonstration.} Le cas des groupes diédraux doit être traité séparément, voir Appendice \ref{appen-B}. On passe ensuite en revue la liste des autres groupes de Coxeter finis irréductibles \cite[Chap. VI, \S 4]{Bo} et on constate que leur 2-Sylow sont produits de $2$-Sylow de groupes symétriques ou alternés, voir les points (1) à (10) ci-après. On conclut en invoquant \ref{tensnf}.

\bigskip
(1) Le groupe $\mathrm{W}(\mathbf{A}_{n})$ ($n\geq 1$) est isomorphe à $\mathfrak{S}_{n+1}$.

\medskip
(2) Les $2$-Sylow de $\mathrm{W}(\mathbf{B}_{n})$ ($n\geq 2$) et $\mathfrak{S}_{2n}$ sont isomorphes ($\mathrm    {W}(\mathbf{B}_{n})$ est isomorphe à $\mathfrak{S}_{n}\wr\{\pm 1\}=\mathfrak{S}_{n}\wr\mathfrak{S}_{2}$ qui s'identifie à un sous-groupe d'indice impair de $\mathfrak{S}_{2n}$ (utiliser \ref{alpha})).

\medskip
(3) Les $2$-Sylow de $\mathrm{W}(\mathbf{D}_{n})$ ($n\geq 3$) et $\mathfrak{A}_{2n}$ sont  isomorphes ($\mathrm{W}(\mathbf{D}_{n})$ est isomorphe au noyau de l'homomorphisme composé $\mathfrak{S}_{n}\wr\mathfrak{S}_{2}\hookrightarrow\mathfrak{S}_{2n}\overset{\epsilon_{2n}}{\to}\mathbb{Z}/2$ qui s'identifie à un sous-groupe d'indice impair de $\mathfrak{A}_{2n}$).

\medskip
(4) Les $2$-Sylow de $\mathrm{W}(\mathbf{E}_{8})$ et $\mathfrak{A}_{16}$ sont isomorphes ($\mathrm{W}(\mathbf{D}_{8})$ s'identifie à un sous-groupe de $\mathrm{W}(\mathbf{E}_{8})$ d'indice $135$).

\medskip
(5) Les $2$-Sylow de $\mathrm{W}(\mathbf{E}_{7})$ et $\mathfrak{A}_{12}\times\mathfrak{S}_{2}$ sont isomorphes ($\mathrm{W}(\mathbf{E}_{7})$ contient un sous-groupe isomorphe à $\mathrm{W}(\mathbf{D}_{6})\times\mathfrak{S}_{2}$, d'indice $63$).

\medskip
(6) Les $2$-Sylow de $\mathrm{W}(\mathbf{E}_{6})$ et $\mathfrak{A}_{10}$ sont isomorphes ($\mathrm{W}(\mathrm{E}_{6})$ contient un sous-groupe isomorphe à $\mathrm{W}(\mathbf{D}_{5})$, d'indice $27$).

\medskip
(7) Les $2$-Sylow de $\mathrm{W}(\mathbf{F}_{4})$ et $\mathfrak{S}_{8}$ sont isomorphes ($\mathrm{W}(\mathbf{F}_{4})$ contient un sous-groupe d'indice $3$ isomorphe à $\mathfrak{S}_{4}\wr\{\pm 1\}=\mathfrak{S}_{4}\wr\mathfrak{S}_{2}$).

\medskip
(8) Le $2$-Sylow de $\mathrm{W}(\mathbf{G}_{2})$ est isomorphe à $\mathfrak{S}_{2}\times\mathfrak{S}_{2}$.

\medskip
(9) Le groupe $\mathrm{W}(\mathbf{H}_{3})$ est isomorphe à $\mathfrak{A}_{5}\times\mathfrak{S}_{2}$.

\medskip
(10) Les $2$-Sylow de $\mathrm{W}(\mathbf{H}_{4})$ et $\mathfrak{A}_{8}$ sont isomorphes (voir Appendice \ref{appen-C}).

\vspace{1.5cm}
\begin{center}
\large
\textbf{Appendices }
\end{center}
\normalsize

\renewcommand{\thesect}{\Alph{sect}}
\setcounter{sect}{0}

\sect{Sur la série de Poincaré de $\mathrm{H}^{*}\mathrm{A}_{2^{m}}$}
\phantomsection\label{appen-A}

\medskip
On se propose dans cet appendice de décrire une méthode de calcul par récurrence sur $m$ de la série de Poincaré de $\mathrm{H}^{*}\mathrm{A}_{2^{m}}$ en utilisant la formule de la remarque \ref{sdP}.

\medskip
Cette formule se spécialise de la façon suivante~:
$$
\mathrm{S}(\mathrm{H}^{*}\mathrm{A}_{2^{m}};t)
\hspace{4pt}=\hspace{4pt}
(1-t)\hspace{2pt}\mathrm{S}(\mathrm{H}^{*}\mathrm{S}_{2^{m}};t)
+(1+t)\hspace{2pt}\mathrm{S}(\tau\mathrm{H}^{*}\mathrm{S}_{2^{m}};t)
$$
(on rappelle que la série de Poincaré d'un  $\mathbb{F}_{2}$-espace vectoriel $\mathbb{N}$-gradué $E$ de dimension finie en chaque degré est notée $\mathrm{S}(E;t)$).

\medskip
Le calcul de $\mathrm{S}(\mathrm{H}^{*}\mathrm{A}_{2^{m}};t)$ se fait par récurrence sur $m$ en utilisant la proposition ci-dessous~:

\begin{pro}\phantomsection\label{sdPcq} Soit $M$ un $\mathrm{A}$-module instable avec $\dim_{\mathbb{F}_{2}}M^{n}<\infty$ pour tout $n$~; on a~:
$$
\hspace{24pt}
\mathrm{S}(\mathfrak{S}_{2}M;t)
\hspace{4pt}=\hspace{4pt}
\frac{1}{2}\hspace{2pt}(\hspace{2pt}\mathrm{S}(M;t)^{2}+\mathrm{S}(M;t^{2})\hspace{2pt})+\frac{t}{1-t}\hspace{2pt}\mathrm{S}(M;t^{2})
\hspace{23pt}.
$$
\end{pro}

\textit{Démonstration.} On considère la $\mathcal{U}$-suite exacte donnée par la définition du $\mathrm{A}$-module instable $\mathfrak{S}_{2}M$~:
$$
\hspace{24pt}
0\to\mathfrak{S}_{2}M\to(M\otimes M)^{\mathfrak{S}_{2}}\oplus\mathrm{R}_{1}M\to\Phi M\to 0
\hspace{24pt}
$$
et on utilise~:

\smallskip
-- le lemme \ref{sdPinv} ci-après, dont la démonstration est laissée au lecteur~;

\smallskip
-- le $\mathcal{E}$-isomorphisme $\mathrm{R}_{1}M\cong\Phi M$ du point (f) de \ref{R1}~;

\smallskip
-- l'égalité $\mathrm{S}(\Phi M;t)=\mathrm{S}(M;t^{2})$.
\hfill$\square$

\medskip
\begin{lem}\phantomsection\label{sdPinv} Soient $k$ un corps et $E$ un $k$-espace vectoriel $\mathbb{N}$-gradué avec $\dim_{k}M^{n}<\infty$ pour tout $n$. Soit $S(t)$ la série de Poincaré de $E$, alors la série de Poincaré de $(E\otimes E)^{\mathfrak{S}_{2}}$ est $(S(t)^2+S(t^{2}))/2$.
\end{lem}

\medskip
Le calcul de $\mathrm{S}(\tau\mathrm{H}^{*}\mathrm{S}_{2^{m}};t)$ se fait par récurrence sur $m$ en utilisant celui de $\mathrm{S}(\mathrm{H}^{*}\mathrm{S}_{2^{m}};t)$, le lemme \ref{sdPinv} et la proposition ci-dessous~:

\begin{pro}\phantomsection\label{M/tauM} Soit $M$ un $\mathrm{P}$-$\mathrm{A}$-module instable~; si $M$ est réduit alors on a un isomorphisme (canonique) de $\mathrm{P}$-$\mathrm{A}$-modules instables
$$
\hspace{24pt}
\mathfrak{S}_{2}M\slash\tau\mathfrak{S}_{2}M
\hspace{4pt}\cong\hspace{4pt}
(M\otimes M)^{\mathfrak{S}_{2}}
\slash
(\tau M\otimes \tau M)^{\mathfrak{S}_{2}}
\hspace{23pt}.
$$
\end{pro}

\textit{Démonstration.} On considère la $(\mathrm{P}\text{-}\mathcal{U})$-suite exacte donnée par la définition du $\mathrm{P}$-$\mathrm{A}$-module instable $\mathfrak{S}_{2}M$~:
$$
\hspace{24pt}
0\to\mathfrak{S}_{2}M\to(M\otimes M)^{\mathfrak{S}_{2}}\oplus\theta\hspace{1pt}\mathrm{R}_{1}\mathcal{O}M\to\theta\hspace{1pt}\Phi \mathcal{O}M\to 0
\hspace{24pt}
$$
(la notation $\mathcal{O}$ désigne à nouveau le foncteur oubli $\mathrm{P}\text{-}\mathcal{U}\to\mathcal{U}$). On applique l'endofoncteur $\tau$ à cette suite exacte~;  comme $\tau$ est exact à gauche et que l'homomorphisme  $\theta\hspace{1pt}\mathrm{R}_{1}\mathcal{O}M\to\theta\hspace{1pt}\Phi\mathcal{O}M$ est surjectif on obtient encore une $(\mathrm{P}\text{-}\mathcal{U})$-suite exacte~:
$$
\hspace{24pt}
0\to\tau\mathfrak{S}_{2}M\to\tau((M\otimes M)^{\mathfrak{S}_{2}})\oplus\theta\hspace{1pt}\mathrm{R}_{1}\mathcal{O}M\to\theta\hspace{1pt}\Phi \mathcal{O}M\to 0
\hspace{24pt}
$$
qui coïncide avec la suite exacte $(\mathcal{S}_{2})$ qui apparaît dans la démonstration de la proposition \ref{quad-P-nilf-bts}\footnote{On observera que l'argument employé ici pour obtenir cette suite exacte évite de faire appel à \ref{sigmatau}.}. La contemplation du $(\mathrm{P}\text{-}\mathcal{U})$-diagramme commutatif
\footnotesize
$$
\begin{CD}
0@>>>\mathfrak{S}_{2}M
@>>>(M\otimes M)^{\mathfrak{S}_{2}}\oplus\theta\hspace{1pt}\mathrm{R}_{1}\mathcal{O}M
@>>>\theta\hspace{1pt}\Phi \mathcal{O}M@>>>0 \\
&& @AAA @AAA @AAA \\
0@>>>\tau\mathfrak{S}_{2}M
@>>>\tau((M\otimes M)^{\mathfrak{S}_{2}})\oplus\theta\hspace{1pt}\mathrm{R}_{1}\mathcal{O}M
@>>>\theta\hspace{1pt}\Phi \mathcal{O}M@>>>0
\end{CD}
$$
\normalsize
donne un $(\mathrm{P}\text{-}\mathcal{U})$-isomorphisme $\mathfrak{S}_{2}M\slash\tau\mathfrak{S}_{2}M\cong(M\otimes M)^{\mathfrak{S}_{2}}\slash\tau((M\otimes M)^{\mathfrak{S}_{2}})$. Si~$M$ est réduit on a en outre $\tau((M\otimes M)^{\mathfrak{S}_{2}})=(\tau M\otimes \tau M)^{\mathfrak{S}_{2}}$. Pour se convaincre de cette égalité, on procède comme dans la démonstration de \ref{tens-P-nilf-2}~: on pose $K=\tau M$, $L=M/K$ et on observe que $L$ est $\mathrm{P}$-libre ce qui implique que l'on a $M\approx K\oplus L$ dans la catégorie $\mathrm{P}\text{-}\mathcal{E}$.
\hfill$\square$

\bigskip
\begin{exples}  Soit $m$ un entier naturel.

\smallskip
-- Pour $m\geq 2$, on a $\dim_{\mathbb{F}_{2}}\mathrm{H}^{1}\mathrm{A}_{2^{m}}=m$.

\smallskip
-- Pour $m\geq 3$, on a $\dim_{\mathbb{F}_{2}}\mathrm{H}^{2}\mathrm{A}_{2^{m}}=\frac{m^3-m+18}{6}$.
\end{exples}

\bigskip
\begin{exples}\phantomsection\label{H1-H2} Soit $X$ un espace topologique avec $\mathrm{H}^{*}X$ de dimension (sur~$\mathbb{F}_{2}$) finie en chaque degré~; on considère l'espace topologique
$$
\hspace{24pt}
\mathrm{A}_{4}X
\hspace{4pt}:=\hspace{4pt}
\mathrm{E}\mathrm{A}_{4}\times_{\mathrm{A}_{4}}X^{4}
\hspace{23pt},
$$
$\mathrm{A}_{4}$ agissant sur $X^{4}$ \textit{via} son inclusion dans~$\mathfrak{S}_{4}$. On dispose pour cet espace d'un énoncé analogue à celui de \ref{cqev}~:

\smallskip
{\em On a un isomorphisme canonique de $\mathbb{F}_{2}$-espaces vectoriels gradués
$$
\hspace{24pt}
\mathrm{H}^{*}\mathrm{A}_{4}X
\hspace{4pt}\cong\hspace{4pt}
\mathrm{H}^{*}(\mathrm{A}_{4};\mathrm{H}^{*}X\otimes\mathrm{H}^{*}X\otimes\mathrm{H}^{*}X\otimes\mathrm{H}^{*}X)
\hspace{23pt},
$$
$\mathrm{A}_{4}$ agissant sur $\mathrm{H}^{*}X\otimes\mathrm{H}^{*}X\otimes\mathrm{H}^{*}X\otimes\mathrm{H}^{*}X$ {\em via}~son inclusion dans $\mathfrak{S}_{4}$.}

\smallskip
On en déduit facilement la série de Poincaré de $\mathrm{H}^{*}\mathrm{A}_{4}X$ en fonction de celle de $\mathrm{H}^{*}X$, disons $S(t)$~:
\begin{multline*}
\mathrm{S}(\mathrm{H}^{*}\mathrm{A}_{4}X;t)
\hspace{4pt}=\hspace{4pt} \\
\frac{1}{(1-t)^{2}}S(t^4)
+\frac{3}{2(1-t)}(S(t^{2})^{2}-S(t^{4}))
+\frac{1}{4}(S(t)^{4}-3\hspace{1pt}S(t^{2})^{2}+2\hspace{1pt}S(t^{4}))\hspace{4pt}.
\end{multline*}
On constate que l'on a $\mathrm{S}(\mathrm{H}^{*}\mathrm{A}_{8};t)=\mathrm{S}(\mathrm{H}^{*}\mathrm{A}_{4}\mathrm{B}\hspace{0,5pt}\mathbb{Z}/2;t)$, égalité  en accord avec l'isomorphisme de groupes $\mathrm{A}_{8}\cong\mathrm{A}_{4}\ltimes(\mathbb{Z}/2)^4$ fourni par le cas $m=3$~de \ref{sign-2}. Par contre la formule ci-dessus donne $\dim_{\mathbb{F}_{2}}\mathrm{H}^{2}(\mathrm{A}_{4}\ltimes(\mathrm{S}_{4})^{4})=15$ alors que l'on~a $\dim_{\mathbb{F}_{2}}\mathrm{H}^{2}\mathrm{A}_{16}=13$~; on ne peut donc avoir $\mathrm{A}_{16}\cong\mathrm{A}_{4}\ltimes(\mathrm{S}_{4})^{4}$ \ldots isomorphisme auquel l'énoncé \ref{sign-1} nous a un temps trompeusement laissé croire.
\end{exples}

\sect{Sur l'application de Quillen pour la cohomologie modulo $2$ \\ des groupes diédraux}
\phantomsection\label{appen-B}

\bigskip
Le calcul de la cohomologie modulo $2$ des groupes diédraux est classique. On le trouvera par exemple dans \cite{AM}~; le fait que l'application de Quillen pour ces groupes soit un isomorphisme est implicite dans cette référence. Nous avons choisi pour notre présentation (loin d'être géodésique~!) du calcul en question de mettre en avant la représentation de Coxeter des groupes diédraux.

\medskip
Soit $n\geq 1$ un entier,  on pose  $\mathrm{D}_{2n}:=\mathbb{Z}/2\ltimes\mathbb{Z}/n$, $\mathbb{Z}/2$ agissant sur $\mathbb{Z}/n$ par multiplication par $-1$~; $\mathrm{D}_{2n}$ est le {\em groupe diédral d'ordre $2n$}\footnote{$\mathrm{D}_{6}$ est isomorphe à $\mathfrak{S}_{3}$, $\mathrm{D}_{8}$ est isomorphe à un $2$-Sylow de $\mathfrak{S}_{4}$ et $\mathrm{D}_{12}$ est isomorphe à $\mathrm{W}(\mathbf{G}_{2}$).}. On écrit $n=2^{m}i$ avec $m\geq 0$ et $i\equiv 1\bmod{2}$~; la suite exacte scindée
$$
\hspace{24pt}
1\to\mathbb{Z}/i\to\mathrm{D}_{2n}\to\mathrm{D}_{2^{m+1}}\to 1
\hspace{24pt}
$$
montre à la fois que le $2$-Sylow de $\mathrm{D}_{2n}$ est $\mathrm{D}_{2^{m+1}}$ et que la restriction, en cohomologie modulo 2, $\mathrm{H}^{*}\mathrm{D}_{2n}\to\mathrm{H}^{*}\mathrm{D}_{2^{m+1}}$ est un isomorphisme.

\bigskip
Le groupe $\mathrm{O}(2)$ est  isomorphe au produit semi-direct  $\mathbb{Z}/2\ltimes\mathrm{SO}(2)$, $\mathbb{Z}/2$ agissant sur le groupe $\mathrm{SO}(2)$ via l'automorphisme involutif $A\mapsto A^{-1}$.

\medskip
On note $\mathrm{s}$ la symétrie orthogonale par rapport à l'axe des $x$. Soit $n\hspace{-3pt}\geq\hspace{-3pt}1$ un entier ; il existe un unique homomorphisme de groupes $\phi_{n}:\mathrm{O}(2)\to\mathrm{O}(2)$ vérifiant $\phi_{n}(\mathrm{s})=\mathrm{s}$ et $\phi_{n}(A)=A^{n}$ pour tout $A$ dans $\mathrm{SO}(2)$. Cet homomorphisme induit une opération sur les fibrés vectoriels   euclidiens (réels) de dimension~$2$ que l'on note encore $\phi_{n}$.

\medskip
On note $\mathrm{r}_{n}$ la rotation de $\mathbb{R}^{2}$ d'angle $\frac{2\pi}{n}$ ($\mathbb{R}^{2}$ est muni de son orientation ``trigonométrique'')~; $\mathrm{r}_{n}\hspace{1pt}\mathrm{s}$ est la réflexion orthogonale par rapport à la droite engendrée par $(\mathrm{cos}\hspace{1pt}\frac{\pi}{n},\mathrm{sin}\hspace{1pt}\frac{\pi}{n})$. On identifie $\mathrm{D}_{2n}$.au sous-groupe de $\mathrm{O}(2)$ engendré par $\{\mathrm{r}_{n}\hspace{1pt}\mathrm{s},\mathrm{s}\}$. L'image réciproque par $\phi_{n}$ du sous-groupe de $\mathrm{O}(2)$ engendré par~$\mathrm{s}$ est $\mathrm{D}_{2n}$.

\begin{pro}\phantomsection\label{diedral-1} Soient $\gamma$ le fibré universel de base $\mathrm{BO}(2)$ et $\mathrm{EO}(2)$ son fibré des repères orthogonaux ($\mathrm{EO}(2)=\mathrm{RO}(\gamma)$). Soit $n\hspace{-3pt}\geq\hspace{-3pt}1$ un entier~; le fibré en sphères $\mathrm{S}(\phi_{n}\gamma)$ s'identifie à l'espace quotient $\mathrm{EO}(2)/\mathrm{D}_{2n}$, le groupe $\mathrm{D}_{2n}$ agissant à droite $\mathrm{EO}(2)$ via son inclusion dans  $\mathrm{O}(2)$.
\end{pro}

\medskip
\textit{Démonstration.} On a par définition $\mathrm{S}(\phi_{n}\gamma)=\mathrm{EO}(2)\times_{\mathrm{O}(2)}\mathrm{S}^{1}$, $\mathrm{O}(2)$ agissant à gauche sur $\mathrm{S}^{1}$ \textit{via} l'homomorphisme  $\phi_{n}:\mathrm{O}(2)\to\mathrm{O}(2)$~; pour la clarté de l'exposition, nous notons ci-après  $\mathrm{S}^{1}_{n}$ l'espace $\mathrm{S}^{1}$ muni de cette action. L'action de $\mathrm{O}(2)$ sur $\mathrm{S}^{1}_{n}$ est transitive et le groupe d'isotropie du point $(1,0)$ est $\mathrm{D}_{2n}$~; on a donc un isomorphisme de $\mathrm{O}(2)$-espaces $\mathrm{S}^{1}_{n}\cong\mathrm{O}(2)/\mathrm{D}_{2n}$ si bien que l'on a $\mathrm{S}(\phi_{n}\gamma)=\mathrm{EO}(2)\times_{\mathrm{O}(2)}\mathrm{O}(2)/\mathrm{D}_{2n}=\mathrm{EO}(2)/\mathrm{D}_{2n}$.
\hfill$\square$

\bigskip
Puisque l'espace $\mathrm{EO}(2)$ est contractile et que l'action de $\mathrm{D}_{2n}$ sur $\mathrm{EO}(2)$ est (topologiquement) libre, la proposition ci-dessus implique~:

\begin{cor}\phantomsection\label{diedral-1-bis} On a un isomorphisme canonique de $\mathrm{H}^{*}\mathrm{BO}(2)$-$\mathrm{A}$-algèbres instables
$$
\hspace{24pt}
\mathrm{H}^{*}\mathrm{D}_{2n}
\hspace{4pt}\cong\hspace{4pt}
\mathrm{H}^{*}\mathrm{S}(\phi_{n}\gamma)
\hspace{23pt}.
$$
\end{cor}

\medskip
Nous allons maintenant étudier $\mathrm{H}^{*}\mathrm{S}(\phi_{n}\gamma)$ à l'aide de la suite exacte de Gysin pour les fibrés en sphères. Pour cela  il nous faut d'abord déterminer la classe d'Euler modulo $2$ de $\phi_{n}\gamma$, c'est-à-dire $\mathrm{w}_{2}(\phi_{n}\gamma)$~; pour préciser la structure de $\mathrm{A}$-module instable de $\mathrm{H}^{*}\mathrm{S}(\phi_{n}\gamma)$ la connaissance de $\mathrm{w}_{1}(\phi_{n}\gamma)$ nous sera aussi utile.

\begin{pro}\phantomsection\label{diedral-2} Soient $\gamma$ le fibré universel de base $\mathrm{BO}(2)$ et $n\hspace{-3pt}\geq\hspace{-3pt}1$ un entier~; on a $\mathrm{w}_{2}(\phi_{n}\gamma)=n\hspace{1pt}\mathrm{w}_{2}(\gamma)$ et $\mathrm{w}_{1}(\phi_{n}\gamma)=\mathrm{w}_{1}(\gamma)$.
\end{pro}

\medskip
\textit{Démonstration.} Soit $\iota:\mathrm{O}(1)\times\mathrm{O}(1)\to\mathrm{O}(2)$ l'inclusion de groupes canonique~;  on va utiliser le fait que l'application
$$
\mathrm{B}\iota:\mathrm{B}(\mathrm{O}(1)\times\mathrm{O}(1))\to\mathrm{BO}(2)
$$
induit un $\mathcal{K}$-isomorphisme de $\mathrm{A}$-algèbres instables $\mathrm{H}^{*}\mathrm{BO}(2)\cong(\mathrm{P}\otimes\mathrm{P})^{\mathfrak{S}_{2}}$.

\medskip
On dispose d'un diagramme commutatif de groupes et d'homorphismes de groupes
$$
\begin{CD}
\mathrm{O}(2)@>\phi_{n}>>\mathrm{O}(2) \\
@A\iota AA @A\iota AA \\
\mathrm{O}(1)\times\mathrm{O}(1)
@>\omega_{n}>>
\mathrm{O}(1)\times\mathrm{O}(1)
\end{CD}
$$

dans lequel l'homomorphisme noté $\omega_{n}$ est défini de la façon suivante~:

\medskip
-- pour $n$ impair, $\omega_{n}$ est l'identité~;

\medskip
-- pour $n$ pair, $\mathrm{pr}_{1}\circ\omega_{n}$ est l'homomorphisme trivial et $\mathrm{pr}_{2}\circ\omega_{n}$ est l'homomorphisme donné par la loi de groupe $\mathrm{O}(1)\times\mathrm{O}(1)\to\mathrm{O}(1)$.

\bigskip
\textit{Vérification.} Soit $(\epsilon_{1} ,\epsilon_{2} )$ un élément de $\{\pm{1}\}\times\{\pm{1}\}=\mathrm{O}(1)\times\mathrm{O}(1)$~; on a~:
$$
\begin{bmatrix}
\epsilon_{1} & 0 \\
0 & \epsilon_{2}
\end{bmatrix}
=
\begin{bmatrix}
1 & 0 \\
0 & \epsilon_{1}\epsilon_{2}
\end{bmatrix}
\begin{bmatrix}
\epsilon_{1} & 0 \\
0 & \epsilon_{1}
\end{bmatrix}
$$
et donc~:
$$
\hspace{24pt}
\phi_{n}
(\begin{bmatrix}
\epsilon_{1} & 0 \\
0 & \epsilon_{2}
\end{bmatrix})
=
\begin{bmatrix}
1 & 0 \\
0 & \epsilon_{1}\epsilon_{2}
\end{bmatrix}
\begin{bmatrix}
\epsilon_{1}^{n} & 0 \\
0 & \epsilon_{1}^{n}
\end{bmatrix}
=
\begin{bmatrix}
\epsilon_{1}^{n} & 0 \\
0 & \epsilon_{1}^{n+1}\epsilon_{2}
\end{bmatrix}
\hspace{23pt}.
$$
\hfill$\square$

\medskip
\begin{rem}\phantomsection\label{diedral-2-bis} L'égalité $\mathrm{w}_{1}(\phi_{n}\gamma)=\mathrm{w}_{1}(\gamma)$ résulte aussi directement de la commutativité du diagramme
$$
\begin{CD}
\mathrm{O}(2)@>\phi_{n}>>\mathrm{O}(2) \\
@V\mathrm{d\acute{e}}tVV @V\mathrm{d\acute{e}}tVV \\
\mathrm{O}(1)
@>=>>
\mathrm{O}(1)
\end{CD}
$$
qui montre que les fibrés $\mathrm{d\acute{e}t}\hspace{1pt}\phi_{n}\gamma$ et $\mathrm{d\acute{e}t}\hspace{1pt}\gamma$ coïncident.
\end{rem}

\medskip
\begin{cor}\phantomsection\label{diedral-3}  Soit $n\hspace{-2pt}\geq\hspace{-2pt}2$ un entier pair~; on a une suite exacte canonique  de $\mathrm{H}^{*}\mathrm{BO}(2)$-$\mathrm{A}$-modules instables
$$
0\to\mathrm{H}^{*}\mathrm{BO}(2)\to\mathrm{H}^{*}\mathrm{D}_{2n}
\to\mathrm{w}_{1}\hspace{1pt}\mathrm{H}^{*}\mathrm{BO}(2)\to 0
$$
($\mathrm{w}_{1}$ étant une abréviation pour $\mathrm{w}_{1}(\gamma)$).
\end{cor}

\bigskip
\textit{Démonstration.} Pour le confort du lecteur, rappelons la théorie, dans notre contexte, de la suite exacte de Gysin pour les fibrés en sphères~:

\medskip
\footnotesize
Soit $\xi$ un fibré vectoriel euclidien de dimension $d$~; on note $\mathrm{B}(\xi)$ sa base, $\mathrm{D}(\xi)$ son fibré en boules, $\mathrm{S}(\xi)$ son fibré en sphères, $\mathrm{p}$ les projections $\mathrm{D}(\xi)\to\mathrm{B}(\xi)$  et $\mathrm{S}(\xi)\to\mathrm{B}(\xi)$, $\mathrm{U}_{\xi}\in\mathrm{H}^{d}(\mathrm{D}(\xi),\mathrm{S}(\xi))$ sa classe de Thom et $\mathrm{e}(\xi)\in\mathrm{H}^{*}\mathrm{B}(\xi)$ sa classe d'Euler. On considère la longue suite exacte de la paire $(\mathrm{D}(\xi),\mathrm{S}(\xi))$ en cohomologie (modulo $2$)~:
$$
\hspace{12pt}
\ldots\to
\mathrm{H}^{*}(\mathrm{D}(\xi),\mathrm{S}(\xi))\to
\mathrm{H}^{*}\mathrm{D}(\xi)\to
\mathrm{H}^{*}\mathrm{S}(\xi)\overset{\partial}{\to}
\mathrm{H}^{*+1}(\mathrm{D}(\xi),\mathrm{S}(\xi))\to\ldots
\hspace{11pt}.
$$
On notera incidemment que $\partial$ peut-être vu comme la désuspension de l'homomorphisme $\Sigma\hspace{1pt}\mathrm{H}^{*}\mathrm{S}(\xi)\to\mathrm{H}^{*}(\mathrm{D}(\xi),\mathrm{S}(\xi))$ induit par l'application naturelle $\mathrm{D}(\xi)/\mathrm{S}(\xi)\to\Sigma\hspace{1pt}\mathrm{S}(\xi)$ ($\mathrm{D}(\xi)$ s'identifie au cylindre de la projection $\mathrm{p}:\mathrm{S}(\xi)\to\mathrm{B}(\xi)$).

\smallskip
Compte tenu des points suivants~:

\smallskip
-- la projection $\mathrm{p}:\mathrm{D}(\xi)\to\mathrm{B}(\xi)$ est une équivalence d'homotopie~;

\smallskip
-- comme $\mathrm{H}^{*}\mathrm{D}(\xi)$-module à droite $\mathrm{H}^{*}(\mathrm{D}(\xi),\mathrm{S}(\xi))$ est libre de dimension $1$ de base $\{\mathrm{U}_{\xi}\}$, ce que équivaut à dire que l'homomorphisme
$\mathrm{H}^{*-d}\mathrm{B}(\xi)\to\mathrm{H}^{*}(\mathrm{D}(\xi),\mathrm{S}(\xi))\hspace{2pt},\hspace{2pt}x\mapsto\mathrm{U}_{\xi}\smile\mathrm{p}^{*}x$ est un isomorphisme (isomorphisme de Thom)~;

\smallskip
-- $\mathrm{p}^{*}\mathrm{e}(\xi)$ est l'image de $\mathrm{U}_{\xi}$ par l'homomorphisme $\mathrm{H}^{*}(\mathrm{D}(\xi),\mathrm{S}(\xi))\to\mathrm{H}^{*}\mathrm{D}(\xi)$~;

\smallskip
la suite exacte ci-dessus peut être réécrite de la façon suivante~:
$$
\hspace{9pt}
\begin{CD}
\ldots\to\mathrm{H}^{*-d}\hspace{1pt}\mathrm{B}(\xi)
@>\mathrm{e}(\xi)\smile\hspace{2pt}->>
\mathrm{H}^{*}\mathrm{B}(\xi)
@>\mathrm{p}^{*}>>
\mathrm{H}^{*}\mathrm{S}(\xi)
@>\partial>>
\mathrm{H}^{*-d+1}\hspace{1pt}\mathrm{B}(\xi)\to\ldots
\end{CD}
\hspace{9pt}.
$$

\medskip
\normalsize
Venons maintenant au cas $\xi=\phi_{n}\gamma$. Comme $n$ est pair la proposition \ref{diedral-2} dit en particulier que la classe $\mathrm{e}(\phi_{n}\gamma)=\mathrm{w}_{2}(\phi\gamma)$ est nulle~; la longue suite exacte de Gysin donne une suite exacte courte :
$$
\hspace{24pt}
0\to\mathrm{H}^{*}\mathrm{BO}(2)\to\mathrm{H}^{*}\mathrm{S}(\phi_{n}\gamma)
\to\mathrm{H}^{*-1}\mathrm{BO}(2)\to 0
\hspace{23pt}.
$$
Il nous reste à préciser la structure de la flèche $\mathrm{H}^{*}\mathrm{S}(\phi_{n}\gamma)\to\mathrm{H}^{*-1}\mathrm{BO}(2)$. Comme nous l'avons rappelé plus haut, celle-ci est la désuspension de la composée de l'homorphisme naturel $\Sigma\hspace{1pt}\mathrm{H}^{*}\mathrm{S}(\phi_{n}\gamma)\to\mathrm{H}^{*}(\mathrm{D}(\phi_{n}\gamma),\mathrm{S}(\phi_{n}\gamma))$ et de l'inverse de l'isomorphisme de Thom
$$
\hspace{24pt}
\Theta:
\mathrm{H}^{*-2}\mathrm{BO}(2)=\mathrm{H}^{*-2}\mathrm{B}(\phi_{n}\gamma)
\to
\mathrm{H}^{*}(\mathrm{D}(\phi_{n}\gamma),\mathrm{S}(\phi_{n}\gamma))
\hspace{23pt}.
$$
La suite exacte courte ci-dessus est donc un avatar de la suivante~:
$$
\hspace{24pt}
0\to\mathrm{H}^{*}\mathrm{BO}(2)\to
\mathrm{H}^{*}\mathrm{S}(\phi_{n}\gamma)
\overset{\partial}{\to}
\Sigma^{-1}\hspace{1pt}\mathrm{H}^{*}(\mathrm{D}(\phi_{n}\gamma),\mathrm{S}(\phi_{n}\gamma))\to 0
\hspace{23pt}.
$$
Comme le $\mathrm{A}$-module $\mathrm{H}^{*}\mathrm{S}(\phi_{n}\gamma)$ est instable et que $\partial$ est surjectif, le $\mathrm{A}$-module $\Sigma^{-1}\hspace{1pt}\mathrm{H}^{*}(\mathrm{D}(\phi_{n}\gamma),\mathrm{S}(\phi_{n}\gamma))$ est lui-aussi instable. Vérifions directement que $\mathrm{H}^{*}(\mathrm{D}(\phi_{n}\gamma),\mathrm{S}(\phi_{n}\gamma))$ est bien une suspension (d'un $\mathrm{A}$-module instable) c'est-à-dire que l'application $\mathrm{Sq}_{0}$ est triviale sur $\mathrm{H}^{*}(\mathrm{D}(\phi_{n}\gamma),\mathrm{S}(\phi_{n}\gamma))$~:

\smallskip
Soit $x$ un élément de la cohomologie modulo $2$ de $\mathrm{BO}(2)$, on a
\begin{multline*}
\mathrm{Sq}_{0}\hspace{1pt}\Theta x=\mathrm{Sq}_{0}\hspace{1pt}(\mathrm{U}_{\phi_{n}\gamma}\smile\mathrm{p}^{*}x)=\mathrm{Sq}_{0}\hspace{1pt}\mathrm{U}_{\phi_{n}\gamma}\smile\mathrm{Sq}_{0}\hspace{1pt}\mathrm{p}^{*}x \\
=\mathrm{Sq}^{2}\hspace{1pt}\mathrm{U}_{\phi_{n}\gamma}\smile\mathrm{Sq}_{0}\hspace{1pt}\mathrm{p}^{*}x 
=(\mathrm{U}_{\phi_{n}\gamma}\smile\mathrm{p}^{*}\mathrm{w}_{2}(\phi_{n}\gamma))\smile\mathrm{Sq}_{0}\hspace{1pt}\mathrm{p}^{*}x
\end{multline*}
(se rappeler que les classes de Stiefel-Whitney d'un fibré réel $\xi$ peuvent être définies par la formule de Thom $\mathrm{w}_{i}(\xi)=\Theta^{-1}(\mathrm{Sq}^{i}\mathrm{U}_{\xi})$)~; d'où $\mathrm{Sq}_{0}\hspace{1pt}\Theta x=0$ puisque la classe $\mathrm{w}_{2}(\phi_{n}\gamma)$ est nulle.

\smallskip
Etudions plus généralement les opérations $\mathrm{Sq}^{i}$ sur $\mathrm{H}^{*}(\mathrm{D}(\phi_{n}\gamma),\mathrm{S}(\phi_{n}\gamma))$, on a
\begin{multline*}
\mathrm{Sq}^{i}\hspace{1pt}(\mathrm{U}_{\phi_{n}\gamma}\smile\mathrm{p}^{*}x)=\\
\mathrm{U}_{\phi_{n}\gamma}\smile\mathrm{p}^{*}\mathrm{Sq}^{i}x \hspace{3pt}+\hspace{3pt}
\mathrm{Sq}^{1}\mathrm{U}_{\phi_{n}\gamma}\smile\mathrm{p}^{*}\mathrm{Sq}^{i-1}x \hspace{3pt}+\hspace{3pt}\mathrm{Sq}^{2}\mathrm{U}_{\phi_{n}\gamma}\smile\mathrm{p}^{*}\mathrm{Sq}^{i-2}x \\
\end{multline*}
soit encore, compte tenu de la formule de Thom et de la proposition \ref{diedral-2}
$$
\mathrm{Sq}^{i}\hspace{2pt}\Theta x
\hspace{4pt}=\hspace{4pt}
\Theta\hspace{2pt}(\mathrm{Sq}^{i}x+\mathrm{w}_{1}\hspace{1pt}\mathrm{Sq}^{i-1}x)
\leqno{(\mathcal{SQ})}
$$
($\mathrm{w}_{1}$ est une abréviation pour $\mathrm{w}_{1}(\gamma)$ et les produits dans $\mathrm{H}^{*}\mathrm{BO}(2)$ sont notés sans le symbole $\smile$).

\smallskip
Soit $\mathrm{w}_{1}\hspace{1pt}\mathrm{H}^{*}\mathrm{BO}(2)$ l'idéal de $\mathrm{H}^{*}\mathrm{BO}(2)$ engendré par $\mathrm{w}_{1}$~; cet idéal est stable sous l'action de $\mathrm{A}$, c'est un sous-$\mathrm{H}^{*}\mathrm{BO}(2)$-$\mathrm{A}$-module instable de $\mathrm{H}^{*}\mathrm{BO}(2)$.

\smallskip
Soit enfin $\nu:\Sigma^{-1}\hspace{1pt}\mathrm{H}^{*}(\mathrm{D}(\phi_{n}\gamma),\mathrm{S}(\phi_{n}\gamma))\to\mathrm{w}_{1}\hspace{1pt}\mathrm{H}^{*}\mathrm{BO}(2)$ l'application de degré zéro définie par $\nu\hspace{1pt}(\Sigma^{-1}\Theta x)=\mathrm{w}_{1}\hspace{1pt}x$. Par définition $\nu$ est un isomorphisme de $\mathrm{H}^{*}\mathrm{BO}(2)$-modules $\mathbb{N}$-gradués, l'égalité $(\mathcal{SQ})$ montre que c'est un isomorphisme de $\mathrm{H}^{*}\mathrm{BO}(2)$-$\mathrm{A}$-modules instables.
\hfill$\square$

\begin{pro}\phantomsection\label{diedral-4} Pour tout entier $n\geq 1$~; le $\mathrm{A}$-module instable $\mathrm{H}^{*}\mathrm{D}_{2n}$ est $\mathcal{N}\hspace{-1,5pt}il$-fermé.
\end{pro}

\bigskip
\textit{Démonstration.} Le cas $n\equiv 1\bmod{2}$ (resp. $n\equiv 2\bmod{4}$) est trivial~: l'homomorphisme canonique $\mathrm{D}_{2n}\to\mathrm{D}_{2}\cong\mathbb{Z}/2$ (resp.  $\mathrm{D}_{2n}\to\mathrm{D}_{4}\cong\mathbb{Z}/2\times\mathbb{Z}/2$) induit un isomorphisme en cohomologie modulo $2$. On traite ci-dessous le cas $n$ pair (qui contient le cas $n\equiv 2\bmod{4}$).

\smallskip
Les $\mathrm{A}$-modules instables $\mathrm{H}^{*}\mathrm{BO}(2)$ et $\mathrm{w}_{1}\hspace{1pt}\mathrm{H}^{*}\mathrm{BO}(2)$ qui apparaissent de part et d'autre de
 $\mathrm{H}^{*}\mathrm{D}_{2n}$ dans la suite exacte de \ref{diedral-3} sont $\mathcal{N}\hspace{-1,5pt}il$-fermés~:
 
 \smallskip
 -- $\mathrm{H}^{*}\mathrm{BO}(2)\cong(\mathrm{P}\otimes\mathrm{P})^{\mathfrak{S}_{2}}$ est $\mathcal{N}\hspace{-1,5pt}il$-fermé d'après \ref{invnf}~;
 
 \smallskip
 -- $\mathrm{w}_{1}\hspace{1pt}\mathrm{H}^{*}\mathrm{BO}(2)$, qui est le noyau du $\mathcal{U}$-morphisme $(\mathrm{P}\otimes\mathrm{P})^{\mathfrak{S}_{2}}\to\mathrm{P}$ induit par le produit $\varphi:\mathrm{P}\otimes\mathrm{P}\to\mathrm{P}$, est $\mathcal{N}\hspace{-1,5pt}il$-fermé d'après \ref{critnilf-1}.
 
\smallskip
Il en résulte que $\mathrm{H}^{*}\mathrm{D}_{2n}$ est  $\mathcal{N}\hspace{-1,5pt}il$-fermé(voir \ref{critnilf-2}).
\hfill$\square$

\bigskip
\textsc{Détermination de $\mathrm{L}(\mathrm{D}_{2n})$ pour $n\geq 1$}

\medskip
Pour $n\equiv 1\bmod{2}$, on a $\mathrm{L}(\mathrm{D}_{2n})=\mathrm{L}(\mathrm{D}_{2})=\mathrm{L}(\mathbb{Z}/2)=\mathrm{H}^{*}\mathbb{Z}/2$.

\medskip
Pou $n\equiv 2\bmod{4}$, on a $\mathrm{L}(\mathrm{D}_{2n})=\mathrm{L}(\mathrm{D}_{4})=\mathrm{L}(\mathbb{Z}/2\times\mathbb{Z}/2)=\mathrm{H}^{*}(\mathbb{Z}/2\times\mathbb{Z}/2)$.

\medskip
Le cas $n\equiv 0\bmod{4}$ est plus intéressant. Soit $\mathrm{E}_{1}$ (resp. $\mathrm{E}_{2}$) le sous-groupe de $\mathrm{D}_{2n}$ engendré par $-\mathrm{s}$ et $\mathrm{s}$ (resp.  $-\mathrm{r}_{n}\mathrm{s}$ et $\mathrm{r}_{n}\mathrm{s}$). On fait les constatations suivantes (pour alléger la notation on pose $\mathcal{Q}:=\mathcal{Q}_{\mathrm{D}_{2n}}$)~:

\smallskip
--\hspace{8pt}pour $i=1,2$, $\mathrm{E}_{i}$ est un $2$-groupe abélien élémentaire avec $\dim_{\mathbb{Z}/2}\mathrm{E}_{i}=2$~;

\smallskip
--\hspace{8pt} $\mathrm{E}_{1}$ et $\mathrm{E}_{2}$ ne sont pas conjugués (en d'autres termes, ne sont pas $\mathcal{Q}$-isomorphes)~;

\smallskip
--\hspace{8pt} tout $2$-sous-groupe abélien élémentaire $E$ de $\mathrm{D}_{2n}$ vérifie $\dim_{\mathbb{Z}/2}E\leq 2$ et si l'on a $\dim_{\mathbb{Z}/2}E=2$ alors $E$ est conjugué soit de $\mathrm{E}_{1}$ soit de $\mathrm{E}_{2}$ (en d'autres termes, est $\mathcal{Q}$-isomorphe soit à $\mathrm{E}_{1}$ soit à $\mathrm{E}_{2}$)~;

\smallskip
--\hspace{8pt}$\mathrm{E}_{1}\cap\mathrm{E}_{2}$ est égal à $\{\pm\mathrm{id}_{\mathbb{R}^{2}}\}$ (ce sous-groupe est le centre de $\mathrm{D}_{2n}$)~;

\smallskip
--\hspace{8pt}le seul élément non-trivial de $\mathrm{Aut}_{\mathcal{Q}}(\mathrm{E}_{1})$ (resp. $\mathrm{Aut}_{\mathcal{Q}}(\mathrm{E}_{2})$), est induit par la conjugaison par $(\mathrm{r}_{n})^{\frac{n}{4}}$ (la rotation d'angle $\frac{\pi}{2}$), il échange  $-\mathrm{s}$ et $\mathrm{s}$ (resp.  $-\mathrm{r}_{n}\mathrm{s}$ et $\mathrm{r}_{n}\mathrm{s}$)~;

\smallskip
--\hspace{8pt}$\mathrm{E}_{1}\cap\mathrm{E}_{2}$ est la diagonale de $\mathrm{E}_{i}\cong\mathbb{Z}/2\times\mathbb{Z}/2$ pour $i=1,2$~;

\smallskip
--\hspace{8pt}on a $\lim_{\mathcal{Q}^{\mathrm{op}}}\mathrm{H}^{*}E=\lim_{\mathcal{R}^{\mathrm{op}}}\mathrm{H}^{*}E$, $\mathcal{R}$ désignant la sous-catégorie pleine de $\mathcal{Q}$ dont les objets sont $\mathrm{E}_{1}$, $\mathrm{E}_{2}$ et $\mathrm{E}_{1}\cap\mathrm{E}_{2}$ (pour la notation $\lim_{(-)^{\mathrm{op}}}\mathrm{H}^{*}E$, voir~\ref{limite}).

\medskip
On en déduit que $\mathrm{L}(\mathrm{D}_{2n})$ est la limite dans la catégorie $\mathcal{K}$ du diagramme
$$
\begin{CD}
(\mathrm{P}\otimes\mathrm{P})^{\mathfrak{S}_{2}}
@>>> \mathrm{P} @<<<
(\mathrm{P}\otimes\mathrm{P})^{\mathfrak{S}_{2}}
\end{CD}
$$
les deux flèches étant induites par le produit de $\mathrm{P}$ soit encore du diagramme
$$
\hspace{24pt}
\begin{CD}
(\mathrm{P}\otimes\mathrm{P})^{\mathfrak{S}_{2}}
@>>> \Phi\mathrm{P} @<<<
(\mathrm{P}\otimes\mathrm{P})^{\mathfrak{S}_{2}}
\end{CD}
\hspace{23pt}.
$$

\medskip
\begin{rem}\phantomsection\label{diedral-4} Ce qui précède montre que la $\mathcal{U}$-suite exacte de \ref{diedral-3} 
$$
0\to\mathrm{H}^{*}\mathrm{BO}(2)\to\mathrm{H}^{*}\mathrm{D}_{2n}
\to\mathrm{w}_{1}\hspace{1pt}\mathrm{H}^{*}\mathrm{BO}(2)\to 0
$$
est scindée pour $n\equiv 0\bmod{4}$ (la restriction $\mathrm{H}^{*}\mathrm{D}_{2n}\to\mathrm{H}^{*}\mathrm{E}_{1}$ induit un isomorphisme de $\mathrm{H}^{*}\mathrm{BO}(2)$ sur $(\mathrm{P}\otimes\mathrm{P})^{\mathfrak{S}_{2}}$). Par contre pour $n\equiv 2\bmod{4}$ l'extension associée est non-triviale car $(\mathrm{P}\otimes\mathrm{P})^{\mathfrak{S}_{2}}$ n'est pas un $\mathcal{U}$-facteur direct de $\mathrm{P}\otimes\mathrm{P}$.
\end{rem}

\sect{Sur les $2$-Sylow du groupe $\mathrm{W}(\mathbf{H}_{4})$}
\phantomsection\label{appen-C}

\medskip
On commence par la description faite dans \cite{Hu} du groupe de Coxeter  $\mathrm{W}(\mathbf{H}_{4})$ comme sous-groupe de  $\mathrm{SO}(4)$.

\medskip
Au préalable quelques rappels (bien classiques~!) concernant les groupes de~Lie $\mathrm{SO}(4)$ et $\mathrm{SO}(3)$~:

\medskip
On note $\mathbb{H}$ le corps des quaternions et $\mathrm{S}^{3}\subset\mathbb{H}^{\times}$ le sous-groupe constitué des  quaternions de norme $1$.

\medskip
Soient $q_{1}$ et $q_{2}$ deux éléments de $\mathrm{S}^{3}$. On observe que l'application
$$
\mathbb{R}^{4}=\mathbb{H}\to\mathbb{H}=\mathbb{R}^{4}
\hspace{24pt},\hspace{24pt}
x\mapsto q_{1}\hspace{1pt}x\hspace{1pt}\bar{q}_{2}
$$
est un élément de $\mathrm{SO}(4)$~; on définit ainsi un homomorphisme de groupes de~Lie
$$
\pi:\mathrm{S}^{3}\times\mathrm{S}^{3}\to\mathrm{SO}(4)
$$
dont le noyau est $\mu_{2}:=\{\pm 1\}$ diagonalement plongé dans $\mathrm{S}^{3}\times\mathrm{S}^{3}$. L'homomorphisme $\pi$ est ``le'' revêtement universel de $\mathrm{SO}(4)$~; il induit un isomorphisme
$$
\hspace{24pt}
\mathrm{S}^{3}\times_{\mu_{2}}\mathrm{S}^{3}
\hspace{4pt}\cong\hspace{4pt}
\mathrm{SO}(4)
\hspace{23pt};
$$
la restriction de $\pi$ à la diagonale induit un homorphisme de groupes  de Lie $\mathrm{S}^{3}\to\mathrm{SO}(3)$ (on identifie l'espace euclidien $\mathbb{R}^{3}$ avec le sous-espace de $\mathbb{H}$ constitué des $x$ avec $x+\bar{x}=0$) qui est ``le'' revêtement universel de $\mathrm{SO}(3)$.

\medskip
On passe ensuite au groupe $\mathrm{O}(4)$. On note $\mathrm{c}$ l'involution
$$
\hspace{24pt}
\mathbb{R}^{4}=\mathbb{H}\to\mathbb{H}=\mathbb{R}^{4}
\hspace{24pt},\hspace{24pt}
x\mapsto\bar{x}
\hspace{23pt};
$$
on constate que  $\mathrm{c}$ est un élément de $\mathrm{O}(4)$ de déterminant $-1$. On considère la suite exacte de groupes
$$
\hspace{24pt}
\begin{CD}
0@>>>\mathrm{S}^{3}\times_{\mu_{2}}\mathrm{S}^{3}
\cong\mathrm{SO}(4)@>>>\mathrm{O}(4)
@>\mathrm{d\acute{e}t}>>\mu_{2}
@>>>0
\end{CD}
\hspace{23pt};
$$
l'involution $\mathrm{c}$ fournit une section de $\mathrm{d\acute{e}t}$ et donc un isomorphisme
$$
\mathrm{O}(4)
\hspace{4pt}\cong\hspace{4pt}
\mu_{2}\ltimes\mathrm{SO}(4)
$$
l'action de $\mu_{2}$ sur $\mathrm{SO}(4)$ étant donnée par la conjugaison par $\mathrm{c}$. On constate que cette action se traduit \textit{via} l'isomorphisme $\mathrm{SO}(4)\cong\mathrm{S}^{3}\times_{\mu_{2}}\mathrm{S}^{3}$ par l'échange des deux  facteurs $\mathrm{S}^{3}$. On voit donc au bout du compte que l'on dispose d'un isomorphisme canonique
$$
\hspace{24pt}
\mathrm{O}(4)
\hspace{4pt}\cong\hspace{4pt} 
(\mathfrak{S}_{2}\ltimes(\mathrm{S}^{3}\times\mathrm{S}^{3}))/\langle(-1,-1)\rangle
\hspace{23pt}.
$$

\medskip
Soit maintenant $\mathrm{ic}:\mathfrak{A}_{5}\to\mathrm{SO}(3)$ un monomorphisme de groupes dont l'image est le groupe des isométries directes d'un icosaèdre régulier centré à l'origine~; il sera commode par la suite de supposer que l'image par $\mathrm{ic}$ du  groupe de Klein $\mathrm{A}_{4}$ (qui est le $2$-Sylow de $\mathfrak{A}_{4}$ et donc un $2$-Sylow de $\mathfrak{A}_{5}$, voir \ref{Klein}) est le sous-groupe de $\mathrm{SO}(3)$ constitué des matrices diagonales. Soit enfin $\widetilde{\mathfrak{A}}_{5}$ l'image inverse de $\mathfrak{A}_{5}$ par le revêtement $\mathrm{S}^{3}\to\mathrm{SO}(3)$, Le résultat de \cite{Hu} est l'isomorphisme de groupes suivant
$$
\hspace{24pt}
\mathrm{W}(\mathbf{H}_{4})
\hspace{4pt}\cong\hspace{4pt} 
(\mathfrak{S}_{2}\ltimes(\widetilde{\mathfrak{A}}_{5}\times\widetilde{\mathfrak{A}}_{5}))/\langle(-1,-1)\rangle
\hspace{24pt}
$$
qui  va nous permettre d'identifier ``le'' $2$-Sylow de $\mathrm{W}(\mathbf{H}_{4})$.

\medskip
L'hypothèse faite sur $\mathrm{ic}$ assure que le groupe $\mathrm{Q}_{8}:=\{\pm 1,\pm\mathbf{i},\pm\mathbf{j},\pm\mathbf{k}\}$ est un $2$-Sylow de $\widetilde{\mathfrak{A}}_{5}$. Il en résulte que $(\mathfrak{S}_{2}\ltimes(\mathrm{Q}_{8}\times\mathrm{Q}_{8}))/\langle(-1,-1)\rangle$ est un $2$-Sylow de $\mathrm{W}(\mathbf{H}_{4})$ (de cardinal $2^{6}$ et d'indice $15^{2}$).

\begin{pro}\phantomsection\label{2-Sylow-H4-1} Le groupe $(\mathfrak{S}_{2}\ltimes(\mathrm{Q}_{8}\times\mathrm{Q}_{8}))/\langle(-1,-1)\rangle$ est isomorphe au produit semi-direct $\mathrm{A}_{4}\ltimes\{\pm 1\}^{4}$,  le groupe de Klein $\mathrm{A}_{4}$ agissant sur $\{\pm 1\}^{4}$ {\em via} son inclusion dans $\mathfrak{S}_{4}$.
\end{pro}

\medskip
\textit{Démonstration.} On considère le monomorphisme canonique
$$
\hspace{6pt}
\gamma:(\mathfrak{S}_{2}\ltimes(\mathrm{Q}_{8}\times\mathrm{Q}_{8}))/\langle(-1,-1)\rangle
\to
(\mathfrak{S}_{2}\ltimes(\mathrm{S}^{3}\times\mathrm{S}^{3}))/\langle(-1,-1)\rangle
=\mathrm{O}(4)
\hspace{5pt}.
$$
Le groupe de Lie $\mathrm{O}(4)$ peut être vu comme le groupe $\mathrm{O}_{4}(\mathbb{R})$ des points réels du groupe algébrique $\mathrm{O}_{4}$~; on observe que l'image de $\gamma$ est contenue dans $\mathrm{O}_{4}(\mathbb{Z})$. Le groupe $\mathrm{O}_{4}(\mathbb{Z})$ est canoniquement isomorphe au produit semi-direct $\mathfrak{S}_{4}\ltimes\{\pm 1\}^{4}$,  $\mathfrak{S}_{4}$ agissant sur $\{\pm 1\}^{4}$ par permutation des facteurs~; $\gamma $ induit donc un monomorphisme $\gamma^{\mathrm{bis}}:(\mathfrak{S}_{2}\ltimes(\mathrm{Q}_{8}\times\mathrm{Q}_{8}))/\langle(-1,-1)\rangle\to\mathfrak{S}_{4}\ltimes\{\pm 1\}^{4}$. On constate par inspection que l'image de l'homomorphisme composé
$$
\begin{CD}
(\mathfrak{S}_{2}\ltimes(\mathrm{Q}_{8}\times\mathrm{Q}_{8}))/\langle(-1,-1)\rangle
@>\gamma^{\mathrm{bis}}>>
\mathfrak{S}_{4}\ltimes\{\pm 1\}^{4}
@>>>\mathfrak{S}_{4}
\end{CD}
$$
est contenue dans $\mathrm{A}_{4}$~; $\gamma$ induit donc au bout du compte un monomorphisme $\gamma^{\mathrm{ter}}:(\mathfrak{S}_{2}\ltimes(\mathrm{Q}_{8}\times\mathrm{Q}_{8}))/\langle(-1,-1)\rangle\to\mathrm{A}_{4}\ltimes\{\pm 1\}^{4}$. Comme la source et le but de $\gamma^{\mathrm{ter}}$ ont même cardinal, à savoir $2^{6}$, $\gamma^{\mathrm{ter}}$  est un isomorphisme.
\hfill$\square$

\medskip
\begin{cor}\phantomsection\label{2-Sylow-H4-2} Les $2$-sous-groupes de Sylow du groupe de Coxeter $\mathrm{W}(\mathbf{H}_{4})$ et du groupe alterné $\mathfrak{A}_{8}$ sont  isomorphes.
\end{cor}

\medskip
\textit{Démonstration.} Compte tenu de la proposition précédente il s'agit de montrer qu'un $2$-Sylow de $\mathfrak{A}_{8}$ est aussi isomorphe au produit semi-direct $\mathrm{A}_{4}\ltimes\{\pm 1\}^{4}$.

\smallskip
L'existence d'un tel isomorphisme est  conséquence du cas $m=3$ de l'énoncé ci-dessous (et de son  corollaire)~:

\begin{pro}\phantomsection\label{sign-1} Soit $m\geq 2$ un entier~; il existe un isomorphisme de groupes
$$
\hspace{24pt}
\alpha_{m}:\mathrm{S}_{2^{m}}\to\mathrm{S}_{2^{m-1}}\ltimes(\mathbb{Z}/2)^{\hspace{0.5pt}2^{m-1}}
\hspace{24pt}
$$
(\hspace{2pt}$\mathrm{S}_{2^{m-1}}$ agissant sur $(\mathbb{Z}/2)^{\hspace{0.5pt}2^{m-1}}$ {\em via} son inclusion dans $\mathfrak{S}_{2^{m-1}}$), tel que l'homomorphisme composé
$$
\begin{CD}
\mathrm{S}_{2^{m}}
@>\alpha_{m}>>
\mathrm{S}_{2^{m-1}}\ltimes(\mathbb{Z}/2)^{\hspace{0.5pt}2^{m-1}}
@>\mathrm{surjection\hspace{3pt}canonique}>>
\mathrm{S}_{2^{m-1}}
@>\epsilon_{2^{m-1}}>>
\mathbb{Z}/2
\end{CD}
$$
est $\epsilon_{2^{m}}$.
\end{pro}

\begin{cor}\phantomsection\label{sign-2} Soit $m\geq 2$ un entier~; l'isomorphisme de groupes $\alpha_{m}$ induit un isomorphisme de groupes
$$
\hspace{24pt}
\mathrm{A}_{2^{m}}
\hspace{4pt}\cong\hspace{4pt}
\mathrm{A}_{2^{m-1}}\ltimes(\mathbb{Z}/2)^{\hspace{0.5pt}2^{m-1}}
\hspace{23pt},
$$
$\mathrm{A}_{2^{m-1}}$ agissant sur $(\mathbb{Z}/2)^{\hspace{0.5pt}2^{m-1}}$ {\em via} son inclusion dans $\mathfrak{S}_{2^{m-1}}$.
\end{cor}

\bigskip
\textit{Démonstration de la proposition \ref{sign-1}.} En fait le cas $m=2$ entraîne le cas général.

\smallskip
On revient sur \ref{2-Sylow-1}. 
On a par définition $\mathrm{S}_{2^{m}}\cong\mathrm{S}_{2^{m-k}}\ltimes(\mathrm{S}_{2^{k}})^{\hspace{0.5pt}2^{m-k}}$ ($ \mathrm{S}_{2^{m-k}}$ agissant sur $(\mathrm{S}_{2^{k}})^{\hspace{0.5pt}2^{m-k}}$ {\em via} son inclusion dans $ \mathfrak{S}_{2^{m-k}}$) pour tout $m$ et tout entier $k$ avec $0\leq k\leq m$. On constate que pour $k\geq 1$ l'homomorphisme $\epsilon_{2^{m}}:\mathrm{S}_{2^{m}}\to\mathbb{Z}/2$, vu comme un homomorphisme $\mathrm{S}_{2^{m-k}}\ltimes(\mathrm{S}_{2^{k}})^{\hspace{0.5pt}2^{m-k}}\to\mathbb{Z}/2$, est l'homomorphisme dont la restriction à $\mathrm{S}_{2^{m-k}}$ est triviale et dont la restriction à $(\mathrm{S}_{2^{k}})^{\hspace{0.5pt}2^{m-k}}$ est composé de l'homomorphisme $(\mathrm{S}_{2^{k}})^{\hspace{0.5pt}2^{m-k}}\to(\mathbb{Z}/2)^{\hspace{0.5pt}2^{m-k}}$ induit par $\epsilon_{2^{k}}$ et de l'addition $(\mathbb{Z}/2)^{\hspace{0.5pt}2^{m-k}}\to\mathbb{Z}/2$. Les cas $k=2$ et $k=1$ de ce qui précède et ``l'associativité'' du produit semi-direct montrent que le cas $m=2$ de la proposition \ref{sign-1} entraînent la cas général.

\smallskip
Reste à traiter le cas $m=2$. On considère la suite exacte de groupes
$$
\hspace{24pt}
\begin{CD}
1@>>>\mathrm{A}_{4}@>>>\mathrm{S}_{4}@>\epsilon_{4}>>\mathbb{Z}/2@>>>1
\end{CD}
\hspace{23pt}.
$$
La donné d'un élément de $\mathrm{S}_{4}$ d'ordre $2$ et de signature non-triviale, par exemple la transposition $(2,1,3,4)$, fournit un isomorphisme de groupes de $\mathrm{S}_{4}$ sur un produit semi-direct $\mathbb{Z}/2\ltimes\mathrm{A}_{4}$. L'action de $\mathbb{Z}/2$ sur $\mathrm{A}_{4}$ est non-triviale sinon $\mathrm{S}_{4}$ serait isomorphe à $(\mathbb{Z}/2)^{3}$. Comme les involutions non-triviales de $\mathrm{Aut}(\mathrm{A}_{4})$ sont conjuguées (on a $\mathrm{Aut}(\mathrm{A}_{4})\simeq\mathrm{GL}_{2}(\mathbb{Z}/2)\simeq\mathfrak{S}_{3}$) il existe une base du $\mathbb{Z}/2$-espace vectoriel $\mathrm{A}_{4}$ qui est stable sous l'action de $\mathbb{Z}/2$. On dispose donc d'un isomorphisme $\alpha:\mathrm{S}_{4}\to\mathbb{Z}/2\ltimes(\mathbb{Z}/2\times\mathbb{Z}/2)$, $\mathbb{Z}/2$ agissant sur $\mathbb{Z}/2\times\mathbb{Z}/2$ par échange des facteurs,  tel que l'homomorphisme composé
$$
\begin{CD}
\mathrm{S}_{4}
@>\alpha>>
\mathbb{Z}/2\ltimes(\mathbb{Z}/2\times\mathbb{Z}/2)
@>\mathrm{surjection\hspace{3pt}canonique}>>
\mathbb{Z}/2
\end{CD}
$$
est $\epsilon_{4}$. L'homomorphisme $\alpha_{2}$ est l'avatar évident de $\alpha$.
\hfill$\square$

\bigskip
\begin{rem} Une façon plus rapide (mais plus indirecte) de se convaincre de l'existence d'un isomorphisme $\mathrm{A}_{8}\simeq\mathrm{A}_{4}\ltimes(\mathbb{Z}/2)^{4}$ est d'utiliser l'isomorphisme exceptionnel $\mathfrak{A}_{8}\cong\mathrm{GL}_{4}(\mathbb{Z}/2)$~:

\smallskip
\footnotesize
On note $\mathrm{M}_{2}(\mathbb{Z}/2)$ le $\mathbb{Z}/2$-espace vectoriel des matrices $2\times 2$ à coefficients dans $\mathbb{Z}/2$. On introduit les éléments suivants de $\mathrm{M}_{2}(\mathbb{Z}/2)$~:
$$
\hspace{12pt}
\mathrm{S}:=\begin{bmatrix} 0 & 1 \\  1 & 0 \end{bmatrix}
\hspace{12pt},\hspace{12pt}
\mathrm{I}:=\begin{bmatrix} 1 & 0 \\  0 &1 \end{bmatrix}
\hspace{12pt},\hspace{12pt}
\mathrm{O}:=\begin{bmatrix} 0 & 0 \\  0 & 0 \end{bmatrix}
\hspace{12pt}
$$
(on observera que l'on a $\mathrm{S}^{2}=\mathrm{I}$), et le sous-ensemble suivant de $\mathrm{GL}_{4}(\mathbb{Z}/2)$~:
$$
\hspace{18pt}
\mathbb{T}
\hspace{4pt}:=\hspace{4pt}
\{\hspace{4pt}\begin{bmatrix} \mathrm{S}^{\nu_{1}} & \mathrm{M} \\  \mathrm{O} & \mathrm{S}^{\nu_{2}} \end{bmatrix}\hspace{4pt}; \hspace{4pt}(\nu_{1},\nu_{2})\in\mathbb{Z}/2\times\mathbb{Z}/2\,\hspace{4pt},\hspace{4pt}M\in\mathrm{M}_{2}(\mathbb{Z}/2)\hspace{4pt}\}
\hspace{17pt}.
$$
On constate que $\mathbb{T}$ est un sous-groupe de $\mathrm{GL}_{4}(\mathbb{Z}/2)$, que ce sous-groupe est un $2$-Sylow et que le sous-ensemble de $\mathbb{T}$ constitué des éléments avec $(\nu_{1},\nu_{2})=(0,0)$ est un sous-groupe distingué de $\mathbb{T}$ isomorphe à $\mathrm{M}_{2}(\mathbb{Z}/2)$ et donc à $(\mathbb{Z}/2)^{4}$. On constate également que l'on a les formules suivantes~:
$$
\begin{bmatrix} \mathrm{S} & \mathrm{O} \\  \mathrm{O} & \mathrm{I} \end{bmatrix}
\begin{bmatrix} \mathrm{I} & M \\  \mathrm{O} & \mathrm{I} \end{bmatrix}
{\begin{bmatrix} \mathrm{S} & \mathrm{O} \\  \mathrm{O} & \mathrm{I} \end{bmatrix}}^{-1}
=\begin{bmatrix} \mathrm{I} & \mathrm{S}\hspace{1pt}M \\  \mathrm{O} & \mathrm{I} \end{bmatrix}
$$
$$
\begin{bmatrix} \mathrm{I} & \mathrm{O} \\  \mathrm{O} & \mathrm{S} \end{bmatrix}
\begin{bmatrix} \mathrm{I} & M \\  \mathrm{O} & \mathrm{I} \end{bmatrix}
{\begin{bmatrix} \mathrm{I} & \mathrm{O} \\  \mathrm{O} & \mathrm{S} \end{bmatrix}}^{-1}
=\begin{bmatrix} \mathrm{I} & M\hspace{1pt}\mathrm{S} \\  \mathrm{O} & \mathrm{I} \end{bmatrix}
$$
$$
\hspace{24pt}
\mathrm{S}\begin{bmatrix} x_{1} & x_{3} \\  x_{2} & x_{4} \end{bmatrix}
=
\begin{bmatrix} x_{2} & x_{4} \\  x_{1} & x_{3} \end{bmatrix}
\hspace{24pt}\text{et}\hspace{24pt}
\begin{bmatrix} x_{1} & x_{3} \\  x_{2} & x_{4} \end{bmatrix}\mathrm{S}
=
\begin{bmatrix} x_{3} & x_{1} \\  x_{4} & x_{2} \end{bmatrix}
\hspace{23pt}.
$$

\medskip
La litanie de formules ci-dessus montre bien que $\mathbb{T}$ est isomorphe au produit semi-direct $\mathrm{A}_{4}\ltimes(\mathbb{Z}/2)^{4}$, le groupe de Klein $\mathrm{A}_{4}$ agissant sur $(\mathbb{Z}/2)^{4}$ \textit{via} son inclusion dans $\mathfrak{S}_{4}$.

\end{rem}
\normalsize

\renewcommand{\refname}{Références}

\end{document}